\documentclass[]{imag-ms-template}
\usepackage{multirow}
\usepackage{xcolor}
\usepackage{graphicx}
\usepackage{bm}
\usepackage{siunitx}
\usepackage{subcaption}
\usepackage{amsmath}
\usepackage{algorithmic}
\usepackage{multicol}
\usepackage{booktabs}
\usepackage{tabularray}
\usepackage{diagbox}
\DeclareMathOperator{\trace}{trace}
\newcommand{\argmax}[3][\in]{\underset{#2 #1 #3}{\arg\mathrm{max}}}
\newcommand{\argmin}[3][\in]{\underset{#2 #1 #3}{\arg\mathrm{min}}}

\title{Overview of Bayesian Solvers in EEG Distributed Source Models: Prior Selection, Algorithmic Implementation, and Depth Bias Reduction}

\author{Joonas Lahtinen,$^{1\ast}$ Alexandra Koulouri,$^{2}$\\
{\small $^{1}$Faculty of Information Technology and Communication Sciences, Tampere University,}\\
{\small 33720 Tampere, Finland}\\
{\small $^{2}$Institute for Mathematical Innovation, University of Bath,}\\
{\small Bath BA2 7AY, United Kingdom}\\
{\small $^\ast$Correspondence:  joonas.j.lahtinen@tuni.fi}
}

\usepackage[natbib=true,backend=biber,style=apa]{biblatex}
\begin{document} 

\maketitle 

\begin{abstract}
Electroencephalography (EEG) source imaging aims to reconstruct the spatial distribution of neural activity within the brain from non-invasive scalp measurements. This inverse problem is severely ill-posed due to the low spatial resolution of EEG and the presence of measurement noise, necessitating robust regularization techniques. Bayesian approaches provide a principled framework for incorporating prior knowledge into the solution, where regularization naturally arises through prior distributions and their associated hyperparameters.

In this work, we provide an overview of  key Bayesian methods for EEG source imaging based on Gaussian, Laplace, and group Laplace priors, with particular emphasis on hierarchical models that promote sparsity. We analyse the connections between these hierarchical formulations and classical optimization techniques, and provide an analytical description of their implementation using expectation–maximization and alternating optimization algorithms.

To address the issue of depth bias—where deeper sources are systematically underestimated or mislocalized—we extend a statistical signal-to-noise ratio (SNR) framework to derive depth-weighted priors that account for differences in how strongly sources at different depths are reflected in the measurements. Finally, we  illustrate the behaviour of the considered models through simulation studies involving sources at varying depths. The results highlight the impact of prior selection and depth weighting on reconstruction accuracy and demonstrate the importance of informed model design for depth-sensitive EEG source localization.

\end{abstract}


\section{Introduction}
Electroencephalography (EEG) is a non-invasive neuroimaging technique that measures the brain's electrical activity with high temporal resolution \citep{Nunez2006}. Despite its advantages, EEG suffers from low spatial resolution and high sensitivity to noise, making the reconstruction of brain source activity from EEG recordings an ill-posed inverse problem \citep{Kaipio2007,Natterer2001}. This ill-posedness means that many different source configurations can explain the observed measurements, requiring robust inversion techniques to obtain a stable and meaningful solution.

The first step in solving this inverse problem is to formulate the EEG forward model in either a continuous or a discrete setting. While analytic solutions are feasible in simplified geometries (e.g., spherical conductivity models), realistic source localization requires modeling the true geometry \citep{VattaFederica2010realvsSpherical,VanrumsteBart2002realvsspherical} and conductivity of the head \citep{Antonakakis2019SEPheadmodels}. This is often achieved using the finite element method (FEM), which allows for anatomically accurate, subject-specific models derived from MRI scans and accommodates complex properties such as anisotropic conductivity. Although FEM modeling is computationally intensive \citep{miinalainen2019realistic}, it has been shown to significantly improve localization accuracy—by as much as 0.5 cm compared to spherical models \citep{VanrumsteBart2002realvsspherical}—highlighting the importance of accurate forward models in achieving reliable source estimates.

In the FEM framework, the brain volume is discretized into small polyhedral elements, and often each element node is associated with a source basis function whose coefficient reflects the amplitude and possibly orientation of neural activity at that location. This formulation is often referred to as distributed source modeling \citep{Michel2004}. The inverse problem then becomes estimating these coefficients from the observed EEG data, for example in \citep{Haemaelaeinen1994, Uutela1999,DaleAnders1993,calvetti2007,WipfDavid2010,Gramfort2014,FRISTON20081104,Gorodnitsky1995FOCUSS,Ou2008,Liang2023}. Because the problem is underdetermined and susceptible to measurement noise, different regularization strategies, including sparsity or smoothness constraints, are employed to ensure stable estimates\citep{Ghosh2009,Darbas2018}.

In particular, a promising framework for addressing these challenges is the Bayesian approaches \citep{Kaipio2004}, which integrate prior knowledge about sources and noise directly into the inversion process. One of the key advantages is that regularization, often introduced as an explicit penalty term in classical optimization, arises naturally through the specification of prior distributions and their associated parameters, e.g., prior variances, \citep{Engl1996,Kaipio2004}. These parameters, often called hyperparameters, which act as regularization coefficients, can be estimated directly from the data via simple statistical inference methods such as Markov Chain Monte Carlo (MCMC) \citep{Kaipio2004,Lucka2016} . As a result, there is no need to rely on heuristic or deterministic techniques—such as L-curves \citep{Hansen1998}—for tuning regularization parameters.
Several studies have addressed the time-invariant EEG source imaging problem from the Bayesian perspective \citep{hamalainen1993,Uutela1999,Sato2004HBM,Calvetti2009,Wipf2009,Friston2008,Lahtinen2022,Koulouri2015,Kaipio2004}, particularly in the context of distributed source models \citep{Michel2004}
\footnote{Other Bayesian approaches not covered in this review include parametric models based on non-linear dipole fitting \citep{SorrentinoAlberto2014}, Bayesian model evidence methods \citep{Mattout2006}, and time-varying frameworks \citep{Gramfort2012,Bekhti2018} to name but a few.}.
These methods aim to impose informative priors that better reflect underlying neural activity while providing robustness to measurement noise and modeling uncertainties. In this work, we revisit a core set of Bayesian approaches for solving the EEG inverse problem, with particular emphasis on the algorithmic implementations to reconstruct focal brain activity. Furthermore, we focus on the well-known depth bias challenge \citep{Badia1998}, a systematic error in EEG source imaging, wherein neural sources located closer to the scalp (i.e., superficial regions) are more likely to be detected and reconstructed with higher accuracy than those located deeper within the brain, such as in subcortical or medial regions.

To reduce this bias, several studies have introduced sensitivity weighting schemes that modify the prior covariance structure or regularization terms to enhance depth sensitivity. Notable contributions include \citep{Koehler1996, PascualMarqui1994, PascualMarqui, Fuchs1999, Wagner2000, PalmeroSoler2007, Buchner1997, Gramfort2014}, each offering different strategies to balance source visibility across depths. In our current review, we build upon the statistical framework for automatic sensitivity weighting proposed by \citet{Calvetti2019AutomaticDepthWeighting}. We extend this framework to a variety of Bayesian prior models and evaluate their effect on localization accuracy through simulation studies, with a focus on reconstructing sources at varying depths.

While numerous studies have independently proposed various Bayesian methods for EEG source imaging, the literature remains fragmented, with many algorithms presented in isolation and under different modeling assumptions. This work aims to consolidate a coherent overview of major Bayesian approaches—particularly for distributed source models—and to offer practical insights into their implementation. First, we provide a comprehensive overview of EEG source imaging from a Bayesian perspective, focusing on models that incorporate Gaussian, Laplace, and group Laplace priors. We also consider hierarchical formulations that promote sparsity and structured sparsity in the source estimates and detail their algorithmic implementations. Second, we elaborate on how prior variances—and the parameters of their associated hyperpriors—can be systematically derived from the signal-to-noise ratio (SNR), by extending
the statistical framework introduced by Calvetti et al. \citep{Calvetti2019AutomaticDepthWeighting}. This approach eliminates the need for manual tuning of regularization parameters, helping to mitigate the well-known depth bias problem in source localization.
 Finally, we present analytical formulations of the optimization procedures associated with each prior type and evaluate their performance through simulation studies. These experiments highlight the impact of model choice on localization accuracy, depth sensitivity, and spatial focality, providing practical guidance for method selection.

\section{Bayesian Inference in EEG source imaging}
In the Bayesian framework, all the variables are random, and estimating the unknown vector ${\bf x}$  is interpreted as a probabilistic inference problem. The observed data, denoted by ${\bf y}$ is linked to the unknowns through the \emph{likelihood} $p({\bf y} \mid {\bf x})$, while prior information about ${\bf x}$ is encoded through a \emph{prior} distribution $p({\bf x})$. The goal is to compute the \emph{posterior} distribution $p({\bf x} \mid {\bf y})$, which combines these two sources of information.

The first step is to define a numerical tractable observation model and thus to derive the likelihood function.  In practice, for the EEG source imaging problem, the domain is discretized, and following the distributed source modelling\citep{Baillet2001,Nunez2006}, a linear mapping that connects dipole sources with observations is obtained, i.e.
\begin{equation}\label{eq:forwardModel}
    {\bf y}= L{\bf x}+\bm\xi,
\end{equation}
where ${\bf y}\in\mathbb{R}^m$ are the measured electrical potentials on the scalp at a certain time point, $L\in\mathbb{R}^{m\times dn}$ is the lead field matrix, 
and ${\bf x}$ are the coefficients of $n$ dipole sources with degree of freedom of the dipole orientation $d$ (e.g., $d=2$ if radial and tangential sources with respect to the gray matter). Therefore, we can write $\mathbf{x} \in \mathbb{R}^{dn}$ into $n$ blocks of size $d$ as:
\[
\mathbf{x} = \begin{bmatrix} \mathbf{x}_1^\top & \cdots & \mathbf{x}_n^\top \end{bmatrix}^\top, \quad \mathbf{x}_k \in \mathbb{R}^d \mbox{ for }k=1,\ldots,n
\]
and  $\bm{\xi}\in\mathbb{R}^m$ is the additive measurement noise. 

If the noise follows a known distribution $p_{\mathrm{noise}}$, we can write the likelihood as:
\begin{equation}
    p({\bf y} \mid {\bf x}) = \int_{\mathbb{R}^m} p({\bf y} \mid {\bf x}, \bm{\xi}) \, p_{\mathrm{noise}}(\bm{\xi} \mid {\bf x}) \, \mathrm{d}\bm{\xi}.
\end{equation}
Given ${\bf x}$ and $\bm{\xi}$, the conditional likelihood is a Dirac delta function:
\[
p({\bf y} \mid {\bf x}, \bm{\xi}) = \delta({\bf y} - L{\bf x} - \bm{\xi}),
\]
leading to a simplified likelihood, determined by the distribution of the noise:
\[
p({\bf y} \mid {\bf x}) = p_{\mathrm{noise}}({\bf y} - L{\bf x} \mid {\bf x}).
\]
In this problem, the additive noise $\bm \xi$ is modeled as Gaussian, i.e.
is $\bm{\xi}\sim\mathcal{N}(\bm{\xi}^*,\Gamma_{\bm{\xi}})$, where $\bm{\xi}^*\in\mathbb{R}^m$ is the mean and $\Gamma_{\bm{\xi}}\in\mathbb{R}^{m\times m}$ is the measurement noise covariance matrix. Therefore, the likelihood is Gaussian given by
\begin{equation}
    p({\bf y}\mid {\bf x})\propto \exp\left(-\frac{1}{2}({\bf y}-L{\bf x}-\bm{\xi}^*)^\top\Gamma_{\bm{\xi}}^{-1}({\bf y}-L{\bf x}-\bm{\xi}^*)\right).
\end{equation}

The posterior distribution is then obtained via Bayes' theorem:
\[
p({\bf x} \mid {\bf y}) = \frac{p({\bf y} \mid {\bf x}) \, p({\bf x})}{p({\bf y})},
\]
Since $p({\bf y}={\bf y}_{ob})$ for a given set of measurements ${\bf y}_{ob}\in\mathbb{R}^m$, it is often omitted, yielding:
\[
p({\bf x} \mid {\bf y}) \propto p({\bf y} \mid {\bf x}) \, p({\bf x}).
\]

In practice, the full posterior is often intractable or difficult to visualize, so we resort to point estimates. In our analysis, we will employ the  \emph{Maximum a posteriori (MAP)} given by 
    \begin{equation}
    \hat{\bf x} = \arg\max_{{\bf x} \in \mathbb{R}^{dn}} \, p({\bf x} \mid {\bf y}).
    \end{equation}
Based on the Gaussian likelihood and Bayes' theorem, we have 
\begin{equation}\label{eq:MAP}
     \hat{\bf x} = \arg\min_{{\bf x} \in \mathbb{R}^{dn}} \,\left\lbrace\frac{1}{2}({\bf y}-L{\bf x}-\bm{\xi}^*)^\top\Gamma_{\bm{\xi}}^{-1}({\bf y}-L{\bf x}-\bm{\xi}^*)-\log(p({\bf x}))\}\right\rbrace
\end{equation}

As we will see later, when the prior includes unknown parameters $\bm{\gamma}$, known as \textit{hyperparameters},  we can model them as random variables with a distribution $p(\bm{\gamma})$, called the \textit{hyperprior}. This leads to a hierarchical Bayesian model where both ${\bf x}$ and $\bm{\gamma}$ are inferred from the data.

\subsection{Bayesian formulation of standard EEG source imaging}
The most common source modelling approach considers smooth sources represented by a zero-mean Gaussian distribution with isotropic covariance. Specifically, ${\bf x}\sim \mathcal{N}(0,\gamma I_{dn})$, with the corresponding prior density given by
\begin{equation}
    p({\bf x}\mid \gamma)\propto \exp\left(-\frac{1}{2\gamma} \left\|{\bf x}\right\|^2_2\right),
\end{equation}
where $\gamma$ is a constant prior variance controlling the energy of the dipole sources. Here $\|{\bf x}\|_2^2=\sum_{k=1}^{dn}x_k^2$ denotes the $\ell_2$ norm.

From (\ref{eq:MAP}) and $p({\bf x})$ Gaussian, the MAP estimate is
\begin{equation}
    \hat{\bf x}=L^T\left(LL^T+\gamma^{-1}\Gamma_{\bm{\xi}}\right)^{-1}({\bf y}-\bm{\xi}^*).
\end{equation}
This is equivalent to Tikhonov regularization \citep{Tikhonov1963}, which is known as Minimum Norm Estimate (MNE) \citep{HamalainenMNE,DaleAnders1993}. 

Another common prior, promoting sparsity (i.e., only a few active sources), is the Laplace prior, given by
\begin{equation}
    \mathrm{Lap}({\bf x}\mid \bm{\mu},\bm{\gamma} )= \prod_{k=1}^{dn}\frac{\gamma_k}{2}\exp\left(-{\gamma_k} \left|\mu_i-x_i\right|\right),
\end{equation}
where often the mean parameter $\bm{\mu}$ is a zero vector and the scaling parameters $\gamma_k$ are equal to each other.

For $\bm{\mu}=0$ and scaling parameter $\gamma$,
the MAP estimate (\ref{eq:MAP}) becomes
\begin{equation}
    \hat{\bf x}=\argmin{{\bf x}}{\mathbb{R}^{dn}}\left\lbrace ({\bf y}-L{\bf x}-\bm{\xi}^*)^\top\Gamma_{\bm{\xi}}^{-1}({\bf y}-L{\bf x}-\bm{\xi}^*)+\gamma \left\|{\bf x}\right\|_1\right\rbrace,
\end{equation}
where $\left\|{\bf x}\right\|_1=\sum_{k=1}^{dn} \left|x_k\right|$ is $\ell_1$ norm.
This is the well-known Minimum Current Estimate (MCE) \citep{Uutela1999}.

 \subsection{Hierarchical Bayesian methods for EEG source imaging}
In the Bayesian inference process, where hyperpriors 
such as Gamma or inverse-Gamma distributions are employed, sparsity-promoting effects similar to those achieved by $\ell_1$ norm-based regularization in classical frameworks can be achieved. 
In particular, by placing hyperpriors on parameters that control the variance or precision of the primary variables (e.g., source amplitudes), the model gains the flexibility to adaptively shrink irrelevant components toward zero while retaining significant ones. This adaptive shrinkage arises naturally in this hierarchical framework. 
 
 Currently, several hierarchical Bayesian methods have been proposed for source localization (add reference here). 
 The advantage of these types of hyperpriors is the increased focality of the estimation compared to the methods with fixed parametrization \citep{CalvettiSparseHBM2020,RezaeiA2021}.

In this framework, the prior modelling becomes
\begin{align*}
\bf{x} &\sim p(\bf{x} \mid \bm{\gamma}) \quad &\text{(Prior)} \\
\bm{\gamma} &\sim p(\bm{\gamma}) \quad &\text{(Hyperprior)}
\end{align*}

The Bayesian inference is 
\begin{equation}
     p({\bf x},\bm{\gamma}\mid{\bf y})\propto p({\bf y}\mid{\bf x})p({\bf x}\mid\bm{\gamma}) p(\bm{\gamma})
\end{equation}
and the standard MAP estimate is
\begin{equation}
({\bf x}_{\text{MAP}}, \bm{\gamma}_{\text{MAP}}) = 
\argmax[]{{\bf x},\bm{\gamma}}{} \left\lbrace 
\log p({\bf y} \mid {\bf x}) + \log p({\bf x} \mid \bm{\gamma}) + \log p(\bm{\gamma})\right\rbrace
\end{equation}
Two main strategies to obtain MAP estimates have been proposed: the Iterative Alternating Sequential (IAS) or Expectation Maximization.
In IAS \citep{calvetti2007b_IAS}, they solve iteratively
\begin{align}
{\bf x}^{(t+1)} &= \argmax{{\bf x}}{\mathbb{R}^{dn}} \left\lbrace 
\log p({\bf y} \mid {\bf x}) + \log p\left({\bf x} \mid \bm{\gamma}^{(t)}\right) \right\rbrace\\
\bm{\gamma}^{(t+1)} &= \argmax[>]{\bm{\gamma}}{{\bf 0}} \left\lbrace 
\log p\left({\bf x}^{(t+1)} \mid \bm{\gamma}\right) + \log p(\bm{\gamma})\right\rbrace
\end{align}

In  {\em Expectation-Maximization} (EM) algorithm \citep{Figueiredo2003,Caron2008},  we alternate between optimization of the lower bound of the log-posterior (E-step) and maximization of the said lower bound with respect to the model parameters (M-step). In practice, this means that in the E-step, we define the objective function:
\begin{equation}\label{eq:EMobj}
    Q({\bf x};\bar{\bm{\gamma}})=\log p({\bf y}\mid{\bf x})+\mathbb{E}_{\bm{\gamma}\mid \hat{\bf x}}\left[\log p({\bf x}\mid \bm{\gamma})\right],
\end{equation}
and in the M-step, we solve the following optimization problem
\begin{equation}
    \hat{\bf x}=\argmax{{\bf x}}{\mathbb{R}^{dn}}\: Q({\bf x};\bar{\bm{\gamma}}).
\end{equation}
The hyperparameters are updated implicitly inside the expectation in Equation (\ref{eq:EMobj}). 
\subsubsection{Conditionally Gaussian with gamma or inverse gamma hyperprior distributions}\label{sub:CG}

 The method proposed by Calvetti {\em et al.} is based on a conditionally Gaussian (CG) model, where the hyperparameters, the Gaussian prior variances, are either set to be gamma (Ga) or inverse gamma (IG) distributed \citep{Calvetti2009}. As described in the publication, both of the distributions of the hyperpriors fall under the {\em generalized gamma distribution} (add reference here). Hence, 
 \begin{align*}
\bf{x} &\sim \mathcal{N}(0,{\bm \gamma}) \quad &\text{(Prior)} \\
\bm{\gamma} &\sim  \mathrm{GenGamma}(\bm{\gamma}; \alpha,\beta,s)\quad &\text{(Hyperprior)}
\end{align*}
 where
 \begin{equation}\label{eq:GenGamma}
     \mathrm{GenGamma}(\bm{\gamma}; \alpha,\beta,s)\propto \exp\left(-\sum_{k=1}^n\frac{\gamma_k^s}{\beta^s}+(s\alpha-1)\sum_{k=1}^n\log \gamma_k\right),
 \end{equation}
 where $\alpha>0$ is the {\em scaling parameter} and $\beta>0$ is the {\em shape parameter}. The gamma distribution follows from the selection $s=1$ and inverse gamma by setting $s=-1$. Here, we consider one $\gamma_k>0$ to be associated with $d$-dimensional dipole ${\bf x}_k\in\mathbb{R}^d$ for each $k=1,\cdots,n$.

The posterior distribution with the generalized gamma hyperprior reads: 
\begin{equation}
\begin{split}
    p(\mathbf{x}, \bm{\gamma} \mid \mathbf{y}) 
    &\propto p(\mathbf{y} \mid \mathbf{x}) \cdot p(\mathbf{x} \mid \bm{\gamma}) \cdot p(\bm{\gamma}) \\
    &= \mathcal{N}\left(\mathbf{y} \mid L\mathbf{x}, \Gamma_{\bm{\xi}}\right) 
       \cdot \prod_{k=1}^n \mathcal{N}\left(\mathbf{x}_k \mid \mathbf{0}, \gamma_k I_d\right) 
       \cdot \prod_{k=1}^n \mathrm{GenGamma}(\gamma_k \mid \alpha, \beta, s) \\
    &\propto \exp\Bigg( 
        -\frac{1}{2} (\mathbf{y} - L\mathbf{x})^\top \Gamma_{\bm{\xi}}^{-1} (\mathbf{y} - L\mathbf{x}) 
        - \sum_{k=1}^n \frac{\|\mathbf{x}_k\|^2}{2\gamma_k} 
        - \sum_{k=1}^n \frac{\gamma_k^s}{\beta^s} 
        + \left(s\alpha - \frac{d+2}{2} \right) \sum_{k=1}^n \log \gamma_k 
    \Bigg).
\end{split}
\end{equation}
where the measurement noise is assumed zero-mean Gaussian; $\bm{\xi}\sim\mathcal{N}({\bf 0},\Gamma_{\bm{\xi}})$. 
\paragraph{Optimizations using IAS:}\hfill\break
\noindent With IAS, we perform iterative MAP estimation by alternating between updates for $\mathbf{x}$ and $\bm{\gamma}$.
\begin{itemize}
    \item Step 1: Optimize $\mathbf{x}$ for fixed $\bm{\gamma}$: \[
\mathbf{x}^* = \argmin{{\bf x}}{\mathbb{R}^{dn}} \left\lbrace \frac{1}{2} \left( \mathbf{y} - L\mathbf{x} \right)^\top \Gamma_{\bm{\xi}}^{-1}\left( \mathbf{y} - L\mathbf{x} \right) + \frac{1}{2} \mathbf{x}^\top\Gamma_{\mathbf{x}}^{-1}\mathbf{x} \right\rbrace  ,
\]
where $\Gamma_{\mathbf{x}} = \mathrm{diag}(\gamma_1 I_d, \dots, \gamma_n I_d) \in \mathbb{R}^{nd \times nd}$.
\vspace{0.1in}
\item {Step 2: Optimize $\bm{\gamma}$ for fixed $\mathbf{x}$}:\\
We optimize each $\gamma_k$ separately. Since, $
\log p(\gamma_k \mid \mathbf{x}_k) \propto 
- \frac{\|\mathbf{x}_k\|^2}{2\gamma_k} 
- \frac{\gamma_k^s}{\beta^s} 
+ \left( s\alpha - \frac{d+2}{2} \right) \log \gamma_k.
$ by taking the derivative w.r.t.\ $\gamma_k$ and setting to zero we have
\[
\frac{\|\mathbf{x}_k\|^2}{2\gamma_k^2} 
- \frac{s \gamma_k^{s-1}}{\beta^s} 
+ \frac{s\alpha - \frac{d+2}{2}}{\gamma_k} = 0.
\]
This is a nonlinear equation with respect to $\gamma_k$, which can be solved numerically (e.g., Newton-Raphson or bisection). In the following algorithm, we also show special cases for the hyperprior.
\end{itemize}


\clearpage

\begin{algorithm}[H]
\caption{Alternating Optimization for Posterior with Generalized Gamma Hyperprior}
\begin{algorithmic}[1]\label{algo:IAS_CG}
\STATE \textbf{Input:} Data $\mathbf{y}$, forward model $L$, noise covariance $\Gamma_{\bm{\xi}}$, hyperparameters $\alpha$, $\beta$, $s \in [-1, 1]$
\STATE \textbf{Initialize:} $\hat{\mathbf{x}}^{(0)}$, set $t \leftarrow 0$
\STATE Initialize $\gamma_k^{(0)} > 0$ for all $k = 1, \ldots, n$ (e.g., $\gamma_k^{(0)} = 1$)
\REPEAT
    \STATE Update $\mathbf{x}^{(t+1)}$ by solving:  
    \STATE \hspace{1em} $\mathbf{x}^{(t+1)} = \left( L^\top \Gamma_{\bm{\xi}}^{-1} L + \Gamma_{\mathbf{x}}^{-1} \right)^{-1} L^\top \Gamma_{\bm{\xi}}^{-1} \mathbf{y}$
    
    \FOR{$k = 1$ to $n$}
        \IF{$s = 1$ \COMMENT{Gamma hyperprior}}
            \STATE $\gamma_k^{(t+1)} = \dfrac{\beta}{2} \left( \alpha - \dfrac{d+2}{2} + \sqrt{\left(\alpha - \dfrac{d+2}{2} \right)^2 + \dfrac{2 \|\mathbf{x}_k^{(t+1)}\|^2}{\beta}} \right)$
        \ELSIF{$s = -1$ \COMMENT{Inverse Gamma hyperprior}}
            \STATE $\gamma_k^{(t+1)} = \dfrac{ \|\mathbf{x}_k^{(t+1)}\|^2 / 2 + \beta }{ \alpha - \dfrac{d+2}{2} }$
        \ELSE 
            \STATE Solve for $\gamma_k^{(t+1)}$:
            \STATE \hspace{1em} $\dfrac{\|\mathbf{x}_k^{(t+1)}\|^2}{2\gamma_k^2} 
        - \dfrac{s \gamma_k^{s-1}}{\beta^s} 
        + \dfrac{s\alpha - \frac{d+2}{2}}{\gamma_k} = 0$
        \ENDIF
    \ENDFOR
\UNTIL{convergence}
\end{algorithmic}
\end{algorithm}
\paragraph{ Optimization based on EM: }\hfill\break
To apply the EM, we need to derive the conditional $p(\bm{\gamma}\mid {\bf x})$. 
\begin{itemize}
    \item when $s=1$ in (\ref{eq:GenGamma}), i.e. $p(\bm{\gamma})\sim \mathrm{GenGamma}(\bm{ \gamma};\alpha,\beta,1)$ we have a Gamma hyperprior. In this case, we can show that $p(\gamma_k \mid {\mathbf{x}}_k)$ is a generalized inverse Gaussian distribution ($\mathrm{GIG}$)\citep{GoodI.J.GenInvGauss}
\[
p(\gamma_k \mid {\mathbf{x}}_k) = \mathrm{GIG}\left(\gamma_k \mid  \frac{2}{\beta}, \|\hat{\mathbf{x}}_k\|^2,\alpha - \frac{d}{2} \right)
\]
which analytically is given by
\[
p(\gamma_k \mid \hat{\mathbf{x}}_k) = 
\frac{\left( \dfrac{\|\hat{\mathbf{x}}_k\|^2 \cdot \beta}{2} \right)^{\frac{\alpha - \frac{d}{2}}{2}}}
{2 K_{\alpha - \frac{d}{2}}\left( \sqrt{ \dfrac{2\|\hat{\mathbf{x}}_k\|^2}{\beta} } \right)} 
\gamma_k^{\alpha - \frac{d}{2} - 1} 
\exp\left( -\frac{\|\hat{\mathbf{x}}_k\|^2}{2\gamma_k} - \frac{\gamma_k}{\beta} \right).
\]
where  $K_{\alpha - \frac{d}{2}}(\cdot)$ is a modified Bessel function of the second kind \citep{MathHandbook2002}. 

Now, using the properties of the generalized inverse Gaussian distribution \citep{Joergensen1982GenInvGauss}, we get
\begin{equation}
    Q({\bf x};\hat{\bf x})=-\frac{1}{2}({\bf y}-L{\bf x})^\top\Gamma_{\bm{\xi}}^{-1}({\bf y}-L{\bf x})-\sum_{k=1}^n\frac{\left\|{\bf x}_k\right\|^2}{\sqrt{2\beta}\left\|\hat{\bf x}_k\right\|}\frac{K_{\alpha-d/2-1}(\sqrt{2/\beta}\left\|\hat{\bf x}_k\right\|)}{K_{\alpha-d/2}(\sqrt{2/\beta}\left\|\hat{\bf x}_k\right\|)},
\end{equation}
\item When $s=-1$ in (\ref{eq:GenGamma}), i.e. $p(\bm{\gamma})\sim \mathrm{GenGa}(\bm{ \gamma};\alpha,\beta,-1)$ and we deal with an inverse gamma hyperprior, we can just use conjugacy to obtain $\gamma_k\mid{\bf x}_k\sim \mathrm{IG}(\left\|\hat{\bf x}_k\right\|^2+\beta,\alpha+d/2)$ and thus
\begin{equation}
\begin{split}
    Q({\bf x};\hat{\bf x})&=-\frac{1}{2}({\bf y}-L{\bf x})^\top\Gamma_{\bm{\xi}}^{-1}({\bf y}-L{\bf x})-\sum_{k=1}^n\frac{\left\|{\bf x}_k\right\|^2}{2}\mathbb{E}_{\gamma_k\sim \mathrm{IG}}\left[\gamma_k^{-1}\right]\\
    &=-\frac{1}{2}({\bf y}-L{\bf x})^\top\Gamma_{\bm{\xi}}^{-1}({\bf y}-L{\bf x})-\sum_{k=1}^n\frac{\left\|{\bf x}_k\right\|^2}{2}\mathbb{E}_{\eta_k\sim \mathrm{Ga}}\left[\eta_k\right] \\
    &=-\frac{1}{2}({\bf y}-L{\bf x})^\top\Gamma_{\bm{\xi}}^{-1}({\bf y}-L{\bf x})-\sum_{k=1}^n\frac{\left\|{\bf x}_k\right\|^2}{2}\frac{\alpha+d/2}{\left\|\hat{\bf x}_k\right\|^2+\beta},
\end{split}
\end{equation}
which differs from the IAS algorithm only in terms of the factor containing the scale and shape parameters.
\end{itemize}

\clearpage

\begin{algorithm}[H]
\caption{Expectation-Maximization for Conditionally Gaussian Model with Generalized Gamma Hyperprior for $s = \pm1$}\label{algo:EM_CG}
\begin{algorithmic}[1]
\STATE \textbf{Input:} Data $\mathbf{y}$, forward model $L$, noise covariance $\Gamma_{\bm{\xi}}$, hyperparameters $\alpha$, $\beta$, $s \in \{-1, 1\}$, perturbation term $\delta$ for Gamma hyperprior
\STATE \textbf{Initialize:} $\hat{\mathbf{x}}^{(0)}$, set $t \leftarrow 0$
\REPEAT
  \STATE \textbf{E-step:} Compute expectations for each $k = 1, \ldots, n$
    \IF{$s = 1$ \COMMENT{Gamma hyperprior}}
        \STATE Compute:
        \[ \frac{1}{w_k^{(t)}} = \sqrt{2\beta} \|\hat{\mathbf{x}}_k^{(t)}\|
        \cdot 
        \frac{K_{\alpha - \frac{d}{2}}\left( \sqrt{\frac{2}{\beta}} \|\hat{\mathbf{x}}_k^{(t)}\| \right)}{K_{\alpha - \frac{d}{2} - 1}\left( \sqrt{\frac{2}{\beta}} \|\hat{\mathbf{x}}_k^{(t)}\| \right)+\delta} \]
    \ELSIF{$s = -1$ \COMMENT{Inverse Gamma hyperprior}}
        \STATE Compute expected inverse:
        \[
        w_k^{(t)} = \frac{\alpha + \frac{d}{2}}{\|\hat{\mathbf{x}}_k^{(t)}\|^2 + \beta}
        \]
    \ENDIF

    \STATE \textbf{M-step: Update $\hat{\mathbf{x}}^{(t+1)}$ using regularized least squares}
    \STATE Define the prior precision matrix:
    \[
    \Gamma_{\mathbf{x}} = \mathrm{diag} \left( \frac{1}{w_1^{(t)}} I_d, \ldots, \frac{1}{w_n^{(t)}} I_d \right)
    \]
    \STATE Solve:
    \[
\hat{\mathbf{x}} = \Gamma_{\mathbf{x}} L^\top \left( L \Gamma_{\mathbf{x}} L^\top + \Gamma_{\bm{\xi}} \right)^{-1} \mathbf{y}
\]
    \STATE $t \leftarrow t + 1$
\UNTIL{convergence}
\STATE \textbf{Output:} Final estimate $\hat{\mathbf{x}}$
\end{algorithmic}
\end{algorithm}

\paragraph{Marginal distributions: }
Based on the choices of the scaling and shape parameters of the hyperprior of the generalized gamma distribution, we can end up with a closed-form expression for the prior $p({\bf x})$ useful for comprehending and modeling sparsity.

For each ${\bf x}_k$ the conditional distribution:
\[
\mathbf{x}_k \mid \gamma_k \sim \mathcal{N}(\mathbf{0}, \gamma_k I_d).
\]
where the hyperprior is $
    p(\gamma_k) \propto \gamma_k^{s\alpha - 1} \exp\left( -\left( \frac{\gamma_k}{\beta} \right)^s \right)$.
    
    Then, the marginal density of $\mathbf{x}$ becomes
\begin{equation}
    p(\mathbf{x}) = \prod_{k=1}^n \int_0^\infty 
    \frac{1}{(2\pi \gamma_k)^{d/2}} \exp\left( -\frac{\|\mathbf{x}_k\|^2}{2\gamma_k} \right)
    \cdot \gamma_k^{s\alpha - 1} \exp\left(-\left( \frac{\gamma_k}{\beta} \right)^s \right)
    d\gamma_k.
\end{equation}

Each integral is a scale mixture of multivariate Gaussians, resulting in a heavy-tailed marginal distribution for each source $\mathbf{x}_k$.

We now consider specific settings of the generalized gamma distribution that lead to known closed-form marginal distributions for each $\mathbf{x}_k$.

\begin{itemize}
    \item {Case 1: \texorpdfstring{$s = 1$}{s = 1} (Gamma Prior)}

If:
\[
\gamma_k \sim \mathrm{Gamma}(\alpha, \beta),
\]
then the marginal belongs to the {\em variance-gamma} family. No general closed form exists, but special cases are tractable.
For example, as has been shown in \citep{Calvetti2019Magic}, a special case of the previous conditionally Gaussian is the group Laplace prior. In particular, the Laplace  distribution can be written as a mixture of a Gaussian with a Gamma distribution
\begin{equation}\label{eq:SparsityPrior}
\mathrm{Lap}({\bf x}_k) =
\pi({\bf x}_k)\propto\exp{\left(-\lambda\|{\bf x}_k\|_2\right)}=\int_{\gamma}\pi({\bf x}_k|\gamma_k)
\;\pi(\gamma_k)\;d\gamma\end{equation} where $\pi({\bf x}_k|\gamma_k)$ is
Gaussian and $\pi(\gamma_k)$ is Gamma with shape parameter $\alpha = \frac{1+d}{2}$ and scale
$\beta=2/\lambda^2$ and $s=1$ and $\|{\bf x}_k\|_2=\sqrt{\sum_{i=1}^d x^2_{(k-1)d+i}}$. 

\item {Case 2: \texorpdfstring{$\gamma_k \sim \mathrm{InvGamma}(\alpha, \beta)$}{Inverse-Gamma Prior}}\\
For $\gamma_k \sim \mathrm{InvGamma}(\alpha, \beta)$, $p({\bf x_k})$ is a classical heavy-tailed prior. 

 the marginal distribution of $\mathbf{x}_k$ is a multivariate Student-$t$ distribution


\begin{equation}
\mathbf{x}_k \sim \text{Student-}t_{2\alpha}\left( \mathbf{0}, \frac{\beta}{\alpha} I_d \right)
\end{equation}

which is given by:

\begin{equation}
p(\mathbf{x}_k) = 
\frac{\Gamma\left(\alpha + \frac{d}{2}\right)}{\Gamma(\alpha) (\pi \beta)^{d/2}} 
\left(1 + \frac{\|\mathbf{x}_k\|^2}{\beta} \right)^{-\alpha - \frac{d}{2}}
\end{equation}




\item {Summary Table}

\begin{table}[h!]
\centering
\renewcommand{\arraystretch}{1.4}
\begin{tabular}{@{}lllll@{}}
\toprule
$s$ & Prior on $\gamma_k$ & Parameters & Marginal $p(\mathbf{x}_k)$ & Notes \\ \midrule
1   & Gamma                   & $\alpha, \beta$ & Variance-Gamma & No simple closed form \\
1   & Exponential             & $\alpha=\frac{1+d}{2}, \beta=2/\lambda^2$       & Laplace-like   & Scalar case is exact Laplace \\
--1 & Inverse-Gamma           & $\alpha, \beta$   & \text{Student-}$t_{2\alpha}\left( \mathbf{0}, \frac{\beta}{\alpha} I_d \right)$ & Closed form for all $d$ \\
$\neq 1$ & GenGamma           & $\alpha, \beta, s$ & No closed form & Numerical methods required \\
\bottomrule
\end{tabular}
\end{table}
\end{itemize}

\subsubsection{Conditionally Laplace with Gamma Hyperprior}\label{sec:CL}

To mitigate the intensity bias induced by $\ell_1$-type norm priors (i.e., Lasso or conditionally Gaussian) \footnote{These distributions rely on the assumption that the source amplitudes follow a distribution with thin tails, which may not be suitable in cases where there are substantial differences in the strengths of the sources and thus may underestimate strong sources, leading to biased reconstructions} indiscriminately, as well as the impact of measurement noise, we adopt the hierarchical adaptive scheme proposed in \citep{Figueiredo2003,Lee2010hierarchicalEXP}. This approach introduces additional flexibility and adaptivity to the prior, allowing for better discrimination between active and inactive sources and reducing the over-shrinkage commonly associated with fixed sparse regularizers.

In particular, in our previous work \citep{Lahtinen2022}, we proposed the hierarchical prior modeling 
\begin{align*}
\bf{x} &\sim \mathrm{Lap}({\bf x}| 0,{\bm \gamma}) \quad &\text{(Prior)} \\
\bm{\gamma} &\sim \mathrm{Gamma}(\bm{\gamma};\alpha,\beta)\quad &\text{(Hyperprior)}
\end{align*}
which is referred to as  {\em Hierarchical Adaptive $L1$-Regularization} (HAL1R) in \cite{Lahtinen2024SHALpR}, and utilizes the gamma distribution as a hypermodel for the Laplace distribution's parameter. 

 Then, the conditionally Laplace, HAL1R, posterior distribution is
\begin{equation}
    p({\bf x},\bm{\gamma}\mid {\bf y})\propto \exp\left(-\frac{1}{2}({\bf y}-L{\bf x})^\top\Gamma_{\bm{\xi}}^{-1}({\bf y}-L{\bf x})-\sum_{k=1}^{n}\gamma_k \left\| {\bf x}_k\right\|_1-\sum_{k=1}^{n}\beta_k \gamma_k+\sum_{k=1}^{n}\alpha_k\log \gamma_k\right).
\end{equation}

Then, the IAS algorithm solves :
\begin{align}
    \hat{\bf x}&=\argmax{{\bf x}}{\mathbb{R}^{dn}}\left\lbrace \frac{1}{2}({\bf y}-L{\bf x})^\top\Gamma_{\bm{\xi}}^{-1}({\bf y}-L{\bf x})+\sum_{i=1}^{dn}\hat{\gamma}_i\left|x_i\right| \right\rbrace,\\
    \hat{\gamma}_i&=\frac{\alpha}{\beta+\left|\hat{x}_i\right|}\quad\textnormal{for}\quad i=1,\cdots,dn.
\end{align}
Similarly to the conditionally Gaussian model with an inverse gamma hyperprior, we get the EM hyperparameter update rule $\bar{\gamma}_i=(\alpha+1)/(\left|\hat{x}_i\right|+\beta)$ for the HAL1R.

\begin{algorithm}[h!]
\caption{IAS or EM for conditional Laplace with same $\alpha$ and $\beta$}
\begin{algorithmic}[1]
\STATE \textbf{Input:} Data $\mathbf{y}$, forward model $L$, noise covariance $\Gamma_{\bm{\xi}}$, hyperparameters $\alpha, \beta$
\STATE \textbf{Initialize:} ${\bf x}^{(0)}=0 $,  $\bm{\gamma}^{(0)} > 0$ and $
\max_{i} \left| \left(L^\top \Gamma_{\bm{\xi}}^{-1} \mathbf{y} \right)_i \right| > \hat{\gamma}_i^{(0)}
$

iteration count $t = 0$
\REPEAT
    \STATE \textbf{Step 1: Update ${\bf x}$ ($\ell_1$ norm minimization)}
    \[
    \mathbf{x}^{(t+1)} = \argmin{\mathbf{x}}{\mathbb{R}^{dn}}\left\{ 
        \frac{1}{2}(\mathbf{y} - L\mathbf{x})^\top \Gamma_{\bm{\xi}}^{-1}(\mathbf{y} - L\mathbf{x})
        + \sum_{i=1}^{dn} \hat{\gamma}_i^{(t)} |x_i| 
    \right\}
    \]
    \COMMENT{to solve see Algorithm 4}
    \STATE \textbf{Step 2: Update $\bm{\gamma}$}\\
    IAS: 
    \[
    \hat{\gamma}_i^{(t+1)} = \frac{\alpha}{\beta + |x_i^{(t+1)}|}, \quad \text{for } i = 1, \dots, dn
    \]

    or EM: 
    \[\bar{\gamma}^{(t+1)}_i=(\alpha+1)/(\left|\hat{x}_i\right|+\beta) , \quad \text{for } i = 1, \dots, dn\]
    
    \STATE $t \gets t + 1$
\UNTIL{Convergence}

\STATE \textbf{Output:} $\hat{\mathbf{x}}, \hat{\bm{\gamma}}$
\end{algorithmic}
\end{algorithm}

To solve the previous $\ell_1$ norm  (or LASSO) problem, a fast and efficient algorithm is the {\em Majorization-Minimization using Local Quadratic Approximation} (MM-LQA)\footnote{Other solver e.g. ADMM(add reference) or Barrier methods (add ref) can be employed}. \citep{KimBaekjin2018LASSOsolver}. Here we give it:
\begin{algorithm}[H]
\caption{MM-LQA for $\ell_1$-Minimization}
\begin{algorithmic}[1]\label{Alg:MM-LQA }
\STATE \textbf{Input:} $L$, $\mathbf{y}$, $\Gamma_{\bm{\xi}}$, initial $\mathbf{u}^{(0)}={\bf x}^{(t)}$, $\bm{\gamma}^{(t)}$, tolerance $\delta$, small constant $\epsilon$
\REPEAT
    \STATE Compute weights:
    \[
    q_i^{(l)} = \frac{2\gamma_i^{(t)}}{2|u_i^{(l)}| + \epsilon}, \quad \text{for } i = 1,\ldots,dn
    \]
    \STATE Form diagonal matrix $Q^{(l)} = \mathrm{diag}(q_1^{(l)}, \ldots, q_{dn}^{(l)})$
    \STATE Solve:
    \[
    {\bf u}^{(l+1)} = \left( L^\top \Gamma_{\bm{\xi}}^{-1} L + Q^{(l)} \right)^{-1} L^\top \Gamma_{\bm{\xi}}^{-1} \mathbf{y}
    \]
    \STATE $l \gets l + 1$
\UNTIL convergence ($\|{\bf u}^{(t)} - {\bf u}^{(t-1)}\| < \delta$)
\STATE \textbf{Update:} ${\bf x}^{(t+1)}={\bf u}^{(l)}$
\end{algorithmic}
\end{algorithm}


\section{SNR-based weights or prior parameter estimation}
To reduce the depth bias \citep{Badia1998}, different prior parameters have to be assigned for different source locations. 
An effective way to tune parameters or hyperparameters, that is called {\em Sensitivity weighting} proposed in \citep{Calvetti2019AutomaticDepthWeighting}, where a single weight is assigned to each location $k$ and employs the signal-to-noise ratio (SNR), i.e.
\begin{equation}\label{eq:SNR}
\mathrm{SNR}=\frac{\mathbb{E}[\|{\bf y}\|_2^2]}{\mathbb{E}[\|\bm{\xi}\|_2^2]}=
\frac{\mathbb{E}[\|L{\bf x}+\bm{\xi}\|_2^2]}{\mathbb{E}[\|\bm{\xi}\|_2^2]}=\frac{\mathrm{trace}\{\mathbb{E}[(L{\bf x}+\bm{\xi})(L{\bf x}+\bm{\xi})^\top]\}}
{\mathrm{trace}\{\mathbb{E}[\bm{\xi}\bm{\xi}^\top]\}}=\frac{\mathrm{trace}\{L\mathbb{E}[{\bf x}{\bf x}^\top]
L^\top\}}{\mathrm{trace}\{\mathrm{Cov}[\bm{\xi}]\}}+1.
\end{equation}
where  $\mathrm{Cov}[\bm{\xi}]=\Gamma_\xi$ is the
covariance of $\xi$ and $\mathbb{E}[{\bf x}{\bf x}^\top]=\mathrm{Cov}[{\bf x}]$, when $\mathbb{E}[{\bf x}]={\bf 0}$.
Here, we have that $\bf{x}$ and $\bm{\xi}$ are independent. Furthermore, in the current analysis, ${\bf x}_k$ are statistically independent. Now, if $q$ is the number of active sources (usually we consider only $q=1$), and we have the expectation
$\mathbb{E}[{\bf x}{\bf x}^\top]=\mathrm{diag}(\theta_1,\ldots,\theta_{n})\otimes I_d$
where $\otimes$ Kronecker \footnote{here we denote the variances in different locations with $\theta$, we note that for the Gaussian and Conditionally Gaussian priors, this coincides with $\gamma$, however this is not the case for the Laplace prior}, we have based on the analysis in \citep{Calvetti2019AutomaticDepthWeighting}
\begin{equation}\label{eq:varianceWithRespectToSNR}
  \theta_k= \frac{1}{q}\frac{ (\mathrm{SNR}-1) {\mathrm{trace\{\Gamma_\xi\}}}}{\|L_k\|^2_F},
\end{equation}
where ${\bf L}_k=[L_{(k-1)d+1}, \ldots,L_{(k-1)d+d}]\in \mathbb{R}^{m\times d}$ and $\|\cdot\|_F^2$ is the Frobenius norm.

The previous formulation allows us to relate the weights $w_k$ or prior variances $\theta_k$ for different priors—such as Gaussian, Laplace, and hierarchical models—to the signal-to-noise ratio (SNR) and the norms of the lead field matrix columns. To estimate these weights or variances for various prior models, including Gaussian, Laplace, Group Laplace, and hierarchical formulations, we define the corresponding prior distributions and their associated variances as follows:

\begin{itemize}
    \item Weighted Gaussian Prior:
 \begin{equation}
    p({\bf x})\propto \prod_{k=1}^n\exp{\left(-w_k\|{\bf x}_k\|^2_2\right)}
\end{equation}
and since $\mathbb{E}[{\bf x}_k]=0$ and $\mathrm{cov}\left[{\bf x}_k\right]=\theta_k I_d$ for a source at location $k$, $w_k=\frac{1}{2\theta_k}$ from the equation (\ref{eq:varianceWithRespectToSNR}) we get 
$w_k=\frac{q\|{\bf L}_k\|_F^2}{2(\mathrm{SNR}-1)\trace\{\Gamma_\xi\}}$

\item Weighted  Laplace prior:
\begin{equation}\label{eq:LapPrior}
    p({\bf x})= \prod_{k=1}^{n}\frac{w_k}{2}\exp\left(-w_k\|{\bf x}_k\|_1\right),
    \end{equation} 
where $\|{\bf x}_k\|_1=\sum_{i=1}^d|x_{(k-1)d+i}|$.     
Since $\mathbb{E}[{\bf x}_k]=0$ and covariance $\mathrm{cov}\left[{\bf x}_k\right]=\theta_k I_d=\frac{2}{w_k^2}$, based on (\ref{eq:varianceWithRespectToSNR}), we get
$w_k=\sqrt{2\frac{q\|{\bf L}_k\|_F^2}{(\mathrm{SNR}-1)\trace\{\Gamma_\xi\}}}$.


\item Weighted group Laplace prior:
\begin{equation}
 p({\bf x})\propto \prod_{k=1}^{n}\frac{w_k}{2}\exp\left(-w_k\|{\bf x}_k\|_2\right), 
\end{equation}
where $\|{\bf x}_k\|_2=\sqrt{\sum_{i=1}^dx_{(k-1)d+i}^2}$. Now for this zero mean distribution, the covariance is $\mathrm{cov}[\mathbf{x}_k] =\theta_k I_d = \frac{d + 1}{w_k^2} \cdot {I}_d$, thus from (\ref{eq:varianceWithRespectToSNR}), we have $w_k=\sqrt{(d+1)\frac{q\|{\bf L}_k\|_F^2}{(\mathrm{SNR}-1)\trace\{\Gamma_\xi\}}}$.

\item For the  
 hierarchical modelling
\begin{align*}
{\bf x}_k &\sim p({\bf x}_k \mid {\gamma}_k) \quad &\text{(Prior)} \\
{\gamma_k} &\sim p({\gamma_k};\alpha_k,\beta_k) \quad &\text{(Hyperprior)}
\end{align*}
for $k=1,\ldots,n$ and  ${\bf x}_k=(x_{d(k-1)+1},\ldots, x_{d(k-1)+d})$, we are interested in specifying hyperpriors $\alpha_k$ and $\beta_k$. Usually, one parameter is set fixed, for example $\alpha_k=\bar{\alpha}$, and $\beta_k$ is estimated.
\begin{enumerate}
  \item For the CG formulation of section \ref{sub:CG}, as it was proposed in \citep{Calvetti2019SensitivityWeight}, we can $$ \mathrm{cov}[{\bf x_k}]=\theta_k I_d=\mathbb{E}[\gamma_k] I_d,$$ where the mean of the hyperparameter $\mathbb{E}[\gamma_k]=\frac{\beta_k\Gamma\left(\frac{\alpha_k+1}{s}\right)}{\Gamma\left(\frac{\alpha_k}{s}\right)}$ with $s>0$ and $\alpha_k>1$ for a GenGamma distribution. 
Therefore,  based on the hypeprior mean and (\ref{eq:varianceWithRespectToSNR}) for these zero mean distributions, we have that 
\begin{equation}
    \beta_k =\frac{{\Gamma\left(\frac{\bar{\alpha}}{s}\right)}}{\Gamma\left(\frac{\bar{\alpha}+1}{s}\right)}  \frac{1}{q}\frac{ (\mathrm{SNR}-1) {\mathrm{trace\{\Gamma_\xi\}}}}{\|{\bf L}_k\|^2_F}
\end{equation}
For special cases where a close form for $p({\bf x}_k)$ exists, the corresponding covariance $\mathrm{cov}[{\bf x}_k]$ or $\mathbb{E}[{\bf x}_k{\bf x}_k^\mathrm{T}]$ can be estimated.
\begin{enumerate}
    \item When $s=-1$ in (\ref{eq:GenGamma}), the marginal distribution is 
$
\mathbf{x}_k \sim \text{Student-}t_{2\alpha_k}\left( \mathbf{0}, \frac{\beta_k}{\alpha_k} I_d \right),
$ the covariance exists for $\alpha_k > 1$, and is given by:
\[
\mathrm{Cov}[\mathbf{x}_k] =\theta_k I_d = \frac{\beta_k}{\alpha_k - 1} I_d
\]
Therefore, we have that 
\begin{equation}
    \beta_k = (\bar{\alpha}-1)  \frac{1}{q}\frac{ (\mathrm{SNR}-1) {\mathrm{trace\{\Gamma_\xi\}}}}{\|{\bf L}_k\|^2_F}.
\end{equation}
\item When $s=1$ and $\alpha=\frac{d+1}{2}$, the marginal $p({\bf x}_k)\propto\exp{\left(-\lambda_k\|{\bf x}_k\|_2\right)}$ with  $\alpha = \frac{1+d}{2}$ and scale
$\beta_k=2/\lambda_k^2$. The covariance is $\mathrm{cov}[{\bf x}_k]=\theta_k I_d = \frac{d+1}{\lambda^2} I_d=\beta_k\frac{d+1}{2} I_d$, and thus  
\begin{equation}
     \beta_k =  \frac{2}{(d+2)q}\frac{ (\mathrm{SNR}-1) {\mathrm{trace\{\Gamma_\xi\}}}}{\|{\bf L}_k\|^2_F}.
\end{equation}

\end{enumerate}

\item  For the weighted version of the CL  of section \ref{sec:CL} when we introduce different hyperparameter per location $k$, we have that 
\begin{align*}\mathbf{x}_k \mid \gamma_k &\sim \mathrm{Lap}({{\bf x}_k}|\mathbf{0}, \gamma_k) \quad  \\
\gamma_k &\sim \mathrm{Ga}(\gamma_k;\alpha_k, \beta_k) \end{align*}
where $p({\bf {x}_k|\gamma_k})\propto \gamma_k^d \exp{\left(-{\gamma_k\|{\bf x}_k\|_1}\right)}$ and 
\[p({\bf x}_k)
\propto \int_0^\infty 
\gamma_k^d \exp\left(-\gamma_k \|\mathbf{x}_k\|_1\right)
\cdot 
\frac{\beta_k^{\alpha_k}}{\Gamma(\alpha_k)} \gamma_k^{\alpha_k - 1} \exp\left(-\beta_k \gamma_k\right)
\, d\gamma_k
\]

So, 
\[
p(\mathbf{x}_k) = 
\frac{\beta_k^{\alpha_k} \cdot \Gamma(d + \alpha_k)}{2\Gamma(\alpha_k)} \cdot 
{(\|\mathbf{x}_k\|_1 + \beta_k)^{-(d + \alpha_k)}},
\]
that is called the multivariate Lomax distribution \citep{Nayak1987MultiLomax}. 
The hyperparameters are computed using
\begin{equation}
    \mathbb{E}\left[{\bf x}_k{\bf x}_k^\top\right]=\frac{\beta_k^{2}\Gamma\left(\alpha_k-2\right)\Gamma\left(3\right)}{\Gamma\left(\alpha_k\right)}I_d=\frac{2\beta_k^2}{(\alpha_k-1)(\alpha_k-2)}I_d
\end{equation}
with condition $\alpha_k>2$. 
Therefore, given  $\bar{\alpha}$ we estimate
 \begin{equation}
     \beta_k = \sqrt{\frac{(\bar{\alpha}-1)(\bar{\alpha}-2)(\mathrm{SNR}-1)\mathrm{trace\{\Gamma_\xi\}}}{2q\|{\bf L}_k\|^2_F}}
 \end{equation}
Algorithm ~\ref{Algo:WeightedL1} gives the steps to solve the weighted conditional $\ell_1$ Laplace problem described here.  
Moreover, considering the hierachical model \begin{align*}\mathbf{x}_k \mid \gamma_k &\sim \mathrm{GroupLap}({{\bf x}_k}|\mathbf{0}, \gamma_k) \quad  \\
\gamma_k &\sim \mathrm{Ga}(\gamma_k;\alpha_k, \beta_k) \end{align*}
where $p({ {\bf x}_k|\gamma_k})\propto \gamma_k^d \exp{\left(-{\gamma_k\|{\bf x}_k\|_2}\right)}$ and $p({\bf x}_k)\propto (\beta_k+\|{\bf x}_k\|_2)^{-(\alpha_k+d)}$, hence we have

\begin{equation}
    \mathbb{E}\left[{\bf x}_k{\bf x}_k^\top\right]=\frac{(d+1)\beta_k^2}{(\bar{\alpha}-1)(\bar{\alpha}-2)}I_d
\end{equation}
and thus
 \begin{equation}
    \beta_k=\sqrt{\frac{(\bar{\alpha}-1)(\bar{\alpha}-2)(\mathrm{SNR}-1)\mathrm{trace\{\Gamma_\xi\}}}{(d+1)q\|L_k\|^2_F}}
 \end{equation}
 Algorithm~\ref{algo:l1-l2weighted} presents the implementation for this setup.
\end{enumerate}
\end{itemize}

In the following table, we summarize the estimated weights or hyperparameters.
\begin{table}[h!]
\centering
\hspace*{-1.6cm} 
\renewcommand{\arraystretch}{1.2}
\small
\resizebox{\textwidth}{!}{
\begin{tabular}{|l|p{5.7cm}|p{6cm}|}
\hline
\textbf{Prior} & \textbf{Marginal prior $p(\mathbf{x})$} & \textbf{Weight / Hyperparameter} \\
\hline

Weighted Gaussian (wG) & 
$p(\mathbf{x}) \propto \prod_{k=1}^n \exp\left(-w_k \|\mathbf{x}_k\|_2^2\right)$ &
$w_k = \dfrac{q \|{\bf L}_k\|_F^2}{2(\mathrm{SNR}-1)\, \mathrm{trace}(\Gamma_\xi)}$ \\
\hline

Weighted Laplace (wL)& 
$p(\mathbf{x}) = \prod_{k=1}^n \dfrac{w_k}{2} \exp\left(-w_k \|\mathbf{x}_k\|_1\right)$ &
$w_k = \sqrt{2 \dfrac{q \|{\bf L}_k\|_F^2}{(\mathrm{SNR}-1)\, \mathrm{trace}(\Gamma_\xi)}}$ \\
\hline

Weighted Group Laplace (wGL) & 
$p(\mathbf{x}) \propto \prod_{k=1}^n \dfrac{w_k}{2} \exp\left(-w_k \|\mathbf{x}_k\|_2\right)$ &
$w_k = \sqrt{(d+1)\dfrac{q \|{\bf L}_k\|_F^2}{(\mathrm{SNR}-1)\, \mathrm{trace}(\Gamma_\xi)}}$ \\
\hline

weighted Conditional Gaussian (wCG) with no closed form marginal & 
\shortstack{$\mathbf{x}_k \sim \mathcal{N}(0, \gamma_kI_d)$\\ $\gamma_k \sim \mathrm{GenGamma}(\alpha_k, \beta_k, s)$}
&
$\beta_k = \dfrac{\Gamma\left(\tfrac{\bar{\alpha}}{s}\right)}{\Gamma\left(\tfrac{\bar{\alpha}+1}{s}\right)} \cdot \dfrac{(\mathrm{SNR}-1) \, \mathrm{trace}(\Gamma_\xi)}{q \|{\bf L}_k\|_F^2}$ \\
\hline

Student-$t$ (wCG: $s = -1$, $\bar{\alpha}>1$, $\beta_k$) & 
$p(\mathbf{x}_k) \sim t_{2\alpha_k}(0, \frac{\beta_k}{\alpha_k} I_d)$ &
$\beta_k = (\bar{\alpha} - 1) \cdot \dfrac{(\mathrm{SNR}-1) \, \mathrm{trace}(\Gamma_\xi)}{q \|{\bf L}_k\|_F^2}$ \\
\hline

Group Laplace (wCG: $s=1$, $\alpha=\frac{d+1}{2}$, $\beta_k = \frac{2}{w_k^2}$ ) & 
$p(\mathbf{x}_k) \propto \exp(-w_k \|\mathbf{x}_k\|_2)$  &
$\beta_k = \dfrac{2}{(d+2)q} \cdot \dfrac{(\mathrm{SNR}-1)\, \mathrm{trace}(\Gamma_\xi)}{\|{\bf L}_k\|_F^2}$ \\
\hline

Weighted conditional  Laplace (wCL), $\bar{\alpha}>2$ & 
$p(\mathbf{x}_k) \propto (\|\mathbf{x}_k\|_1 + \beta_k)^{-(d + \bar{\alpha})}$ &
 $ \beta_k = \sqrt{\frac{(\bar{\alpha}-1)(\bar{\alpha}-2)(\mathrm{SNR}-1)\mathrm{trace\{\Gamma_\xi\}}}{2q\|{\bf L}_k\|^2_F}}$\\
\hline
Weighted conditional Group Laplace (wCGL), $\bar{\alpha}>2$ & 
$p(\mathbf{x}_k) \propto (\|\mathbf{x}_k\|_2 + \beta_k)^{-(d + \bar{\alpha})}$ &
 $ \beta_k = \sqrt{\frac{(\bar{\alpha}-1)(\bar{\alpha}-2)(\mathrm{SNR}-1)\mathrm{trace\{\Gamma_\xi\}}}{(d+1)q\|{\bf L}_k\|^2_F}}$\\
 \hline
\end{tabular}
}
\caption{Summary of priors, expressions of marginal priors, and corresponding weights or hyperparameters derived based on the SNR formula.}\label{table:weights}
\end{table}

\subsection{Algorithms with weights for EEG source imaging}
In this section, we present the algorithms designed to solve the EEG source imaging problem with structured sparsity using weighted priors. First, we describe the  MM-LQA algorithm \ref{MM-LQA_L_p} for minimizing cost functions with either standard Laplace or group Laplace regularization, where the depth weights $w_k$ (estimated in the previous section) are incorporated into a majorization-minimization framework for efficient updates. We then describe two iterative algorithms, employing either the IAS or EM optimization framework, for solving the EEG source problem under the weighted conditional Laplace (wCL) and weighted conditional group Laplace (wCGL) priors. These methods iteratively update the source estimate and associated hyperparameters to reflect both sparsity (or group sparsity) and SNR-based prior weights (or variances). The algorithms also include mechanisms to avoid degeneracy in the initial iterations, ensuring meaningful reconstructions.

\clearpage

\begin{algorithm}[H]
\caption{MM-LQA for Weighted  Laplace or Group Laplace Minimization}
\begin{algorithmic}[1]\label{MM-LQA_L_p}
\STATE \textbf{Input:} $p=1$ (Laplace) or 2 (Group Laplace), data $\mathbf{y}$, operator $L$, noise covariance $\Gamma_{\bm{\xi}}$, weights $w_k$ (from table~\ref{table:weights}), block size $d$, tolerance $\delta$, small $\epsilon > 0$  
\STATE \textbf{Initialize:} $\mathbf{x}^{(0)} \in \mathbb{R}^{nd}$, set $t \gets 0$
\REPEAT
    \FOR{$k = 1$ to $n$}
        \STATE Extract block $\mathbf{x}_k^{(t)} = \left[ {x}^{(t)}_{(k-1)d+1},\ldots,{x}^{(t)}_{(k-1)d+d} \right]$
        \STATE Compute block weight:
        \[
        q_k^{(t)} = \frac{w_k}{\| \mathbf{x}_k^{(t)} \|_p + \epsilon}
        \]
    \ENDFOR
    \STATE Form block-diagonal matrix:
    \[
    Q^{(t)} = \text{block-diag}\left( q_1^{(t)} I_d, \dots, q_n^{(t)} I_d \right)
    \]
    \STATE Update estimate:
    \[
    \mathbf{x}^{(t+1)} = \left( L^\top \Gamma_{\bm{\xi}}^{-1} L + Q^{(t)} \right)^{-1} L^\top \Gamma_{\bm{\xi}}^{-1} \mathbf{y}
    \]
    \STATE $t \gets t + 1$
\UNTIL{convergence: $\| \mathbf{x}^{(t)} - \mathbf{x}^{(t-1)} \| < \delta$}
\STATE \textbf{Return:} $\hat{\mathbf{x}} = \mathbf{x}^{(t)}$
\end{algorithmic}
\end{algorithm}

\begin{algorithm}[H]
\caption{IAS or EM for weighted conditional Laplace (wCL)}\label{Algo:WeightedL1}
\begin{algorithmic}[1]
\STATE \textbf{Input:} EEG Data $\mathbf{y}$, forward model $L$, noise covariance $\Gamma_{\bm{\xi}}$, 
\STATE \textbf{Hyperparameters $\alpha, \beta$:} for IAS: $\hat{\gamma}_k^{(0)}=\frac{\bar{\alpha}}{\beta_k}$ or for EM:  $\bar{\gamma}_k^{(0)}=\frac{\bar{\alpha}+1}{\beta_k}$ select $\bar{\alpha}>2$
and estimate $\beta_k = \sqrt{\frac{(\bar{\alpha}-1)(\bar{\alpha}-2)(\mathrm{SNR}-1)\mathrm{trace\{\Gamma_\xi\}}}{2p\|{\bf L}_k\|^2_F}}$ for $k=1,\ldots,n$. 
\\
\STATE \textbf{Non-degeneracy condition to avoid } $\mathbf{x}^{(1)} = \mathbf{0}$: ensure that
\[
1 < \max_{k,i}\left\{\frac{1}{\hat{\gamma}_k^{(0)}} \left| \left[L^\top \Gamma_{\bm{\xi}}^{-1} \mathbf{y} \right]_{(k-1)d+i} \right|\right\}, \quad \text{for } k = 1, \dots, n \quad \text{and }\quad i = 1, \dots, d
\]
\textit{If the condition is not satisfied, rescale } $\hat{\gamma}_k^{(0)}$ \textit{as :}
\[
\hat{\gamma}_k^{(0)} \leftarrow \hat{\gamma}_k^{(0)} \cdot \mu \cdot \max_{j,i}\left\{\frac{1}{\hat{\gamma}_j^{(0)}} \left| \left[L^\top \Gamma_{\bm{\xi}}^{-1} \mathbf{y} \right]_{(j-1)d+i} \right|\right\}, \quad \text{with } \mu < 1.
\]
\STATE {\textbf{Initialize:}}${\bf x}^{(0)}=0 $  iteration count $t = 0$
\REPEAT
    \STATE \textbf{Step 1: Update ${\bf x}$}
    \[
    \mathbf{x}^{(t+1)} = \argmin{\mathbf{x}}{\mathbb{R}^{dn}}\left\{ 
        \frac{1}{2}(\mathbf{y} - L\mathbf{x})^\top \Gamma_{\bm{\xi}}^{-1}(\mathbf{y} - L\mathbf{x})
        + \sum_{k=1}^{n} \hat{\gamma}_k^{(t)} \|{\bf{x}}_k\|_1 
    \right\}
    \]
    where $\|{\bf x}_k\|_1=\sum_{i=1}^d|x_{(k-1)d+i}|$. 
    \STATE \textbf{Step 2: Update $\bm{\gamma}$}\\
    IAS: 
    \[
    \hat{\gamma}_k^{(t+1)} = \frac{\bar{\alpha}}{\beta_k + \|{\bf x}_k^{(t+1)}\|_1}, \quad \text{for } k = 1, \dots, n
    \]

    or EM: 
    \[\bar{\gamma}^{(t+1)}_k=(\bar{\alpha}+1)/(\left\|\hat{\bf x}_k\right\|_1+\beta_k) , \quad \text{for } k = 1, \dots, n\]
    
    \STATE $t \gets t + 1$
\UNTIL{Convergence}

\STATE \textbf{Output:} $\hat{\mathbf{x}}, \hat{\bm{\gamma}}$
\end{algorithmic}
\end{algorithm}

\begin{algorithm}[H]
\caption{IAS or EM for weighted conditional group Laplace (wCGL)}\label{algo:l1-l2weighted}
\begin{algorithmic}[1]
\STATE \textbf{Input:} Data $\mathbf{y}$, forward model $L$, noise covariance $\Gamma_{\bm{\xi}}$, $d$
\STATE \textbf{Hyperparameters $\alpha, \beta$:} Set $\bar{\alpha}>2$ and
$
 \beta_k=\sqrt{\frac{(\bar{\alpha}-1)(\bar{\alpha}-2)(\mathrm{SNR}-1)\mathrm{trace\{\Gamma_\xi\}}}{(d+1)p\|L_k\|^2_F}}, \quad \text{for } k = 1, \ldots, n
$
\begin{itemize}
    \item For IAS: $\hat{\gamma}_k^{(0)} = \frac{\bar{\alpha} + d - 1}{\beta_k}$
    \item For EM: $\bar{\gamma}_k^{(0)} = \frac{\bar{\alpha} + d}{\beta_k}$
\end{itemize}
\STATE \textbf{Non-degeneracy condition to avoid } $\mathbf{x}^{(1)} = \mathbf{0}$: ensure that
\[
\max_{k = 1, \dots, n} \left\{\frac{1}{\gamma_k^{(0)}}\left\| \left[ L^\top \Gamma_{\bm{\xi}}^{-1} \mathbf{y} \right]_{(k)} \right\|_2\right\} > 1
\]
where $\left[ L^\top \Gamma_{\bm{\xi}}^{-1} \mathbf{y} \right]_{(k)}=\left[\left[L^\top \Gamma_{\bm{\xi}}^{-1} \mathbf{y}\right]_{(k-1)d + 1},\ldots, \left[L^\top \Gamma_{\bm{\xi}}^{-1} \mathbf{y}\right]_{kd}\right]^\mathrm{T}$ denotes the $k^\text{th}$ block of size $d$. 

\textit{If the condition is not satisfied, rescale the weights as:}
\[
\hat{\gamma}_k^{(0)} \leftarrow 
\hat{\gamma}_k^{(0)}\cdot \mu \cdot 
\max_{j=1,\ldots,n} \left\{\frac{1}{\gamma_k^{(0)}}\left\| \left[ L^\top \Gamma_{\bm{\xi}}^{-1} \mathbf{y} \right]_{(j)} \right\|_2\right\}
, \quad \text{with } \mu < 1 \quad \text{for }k=1,\ldots,n. 
\]

\STATE \textbf{Initialize:} $\mathbf{x}^{(0)} = 0$; set iteration count $t = 0$

\REPEAT
    \STATE \textbf{Step 1: Update $\mathbf{x}$}
    \[
    \mathbf{x}^{(t+1)} = \arg\min_{\mathbf{x} \in \mathbb{R}^{dn}} \left\{ 
        \frac{1}{2} (\mathbf{y} - L \mathbf{x})^\top \Gamma_{\bm{\xi}}^{-1} (\mathbf{y} - L \mathbf{x})
        + \sum_{k=1}^{n} \hat{\gamma}_k^{(t)} \| \mathbf{x}_k \|_2
    \right\}
    \]
    where $\| \mathbf{x}_k \|_2 = \sqrt{ \sum_{i=1}^d x_{(k-1)d + i}^2 }$
    
    \STATE \textbf{Step 2: Update $\bm{\gamma}$}:
    
    \textit{IAS update (MAP):}
    \[
    \hat{\gamma}_k^{(t+1)} = \frac{\bar{\alpha}+d-1}{\beta_k + \| \mathbf{x}_k^{(t+1)} \|_2}, \quad \text{for } k = 1, \dots, n
    \]
    
    \textit{EM update (expectation):}
    \[
    \bar{\gamma}_k^{(t+1)} = \frac{\bar{\alpha} + d}{\| \hat{\mathbf{x}}_k \|_2 + \beta_k}, \quad \text{for } k = 1, \dots, n
    \]
    
    \STATE $t \gets t + 1$
\UNTIL{Convergence}

\STATE \textbf{Output:} $\hat{\mathbf{x}}, \hat{\bm{\gamma}}$
\end{algorithmic}
\end{algorithm}

\section{Comparative Study}
A descriptive comparison of Bayesian EEG source imaging solvers developed over the years is not straightforward, primarily because many of these methods have been evaluated in vastly different EEG experimental or simulation settings.
This variability makes it difficult to draw definitive conclusions from previous findings. To address this, and to enable a more objective and quantitative comparison—particularly with regard to sensitivity weighting strategies— we design and conduct a dedicated simulation study in the following section. To ease the reading of this section, we have included Table~\ref{tab:map-costs} that summarizes the names and the Maximum a Posteriori (MAP)  expressions for the various EEG optimizations that were tested here. 
\begin{table}[h!]
\centering
\caption{Bayesian EEG algorithms and their cost functions}
\label{tab:map-costs}
\resizebox{\textwidth}{!}{%
\begin{tabular}{|l|p{13cm}|}
\hline
\textbf{EEG Solver Name} & \textbf{Cost Function} \\
\hline
Weighted Gaussian (wMNE) & 
$\hat{\mathbf{x}} = \arg\min_{\mathbf{x}} \left\{ 
\frac{1}{2} \| \mathbf{y} - L\mathbf{x} \|^2_{\Gamma_\xi^{-1}} + \sum_{k=1}^n w_k \| \mathbf{x}_k \|_2^2 
\right\}$ \\
\hline
Weighted Laplace (wMCE) & 
$\hat{\mathbf{x}} = \arg\min_{\mathbf{x}} \left\{ 
\frac{1}{2} \| \mathbf{y} - L\mathbf{x} \|^2_{\Gamma_\xi^{-1}} + \sum_{k=1}^n w_k \| \mathbf{x}_k \|_1 
\right\}$ \\
\hline
Weighted Group Laplace (wGL)& 
$\hat{\mathbf{x}} = \arg\min_{\mathbf{x}} \left\{ 
\frac{1}{2} \| \mathbf{y} - L\mathbf{x} \|^2_{\Gamma_\xi^{-1}} + \sum_{k=1}^n w_k \| \mathbf{x}_k \|_2 
\right\}$ \\
\hline
Weighted Conditionally Gaussian (wCG) & 
$\hat{\mathbf{x}}, \hat{\bm{\gamma}} = \arg\min_{\mathbf{x}, \bm{\gamma}} \left\{ 
\frac{1}{2} \| \mathbf{y} - L\mathbf{x} \|^2_{\Gamma_\xi^{-1}} + \sum_{k=1}^n \left( 
\frac{\| \mathbf{x}_k \|_2^2}{2\gamma_k} + \text{penalty}(\gamma_k) 
\right) \right\}$ \\
\hline
Weighted Conditional Laplace (wCL) & 
$\hat{\mathbf{x}}, \hat{\bm{\gamma}} = \arg\min_{\mathbf{x}, \bm{\gamma}} \left\{ 
\frac{1}{2} \| \mathbf{y} - L\mathbf{x} \|^2_{\Gamma_\xi^{-1}} + \sum_{k=1}^n \gamma_k \| \mathbf{x}_k \|_1 + \text{penalty}(\gamma_k) 
\right \}$ \\
\hline
Weighted Conditional Group Laplace (wCGL) & 
$\hat{\mathbf{x}}, \hat{\bm{\gamma}} = \arg\min_{\mathbf{x}, \bm{\gamma}} \left\{ 
\frac{1}{2} \| \mathbf{y} - L\mathbf{x} \|^2_{\Gamma_\xi^{-1}} + \sum_{k=1}^n \gamma_k \| \mathbf{x}_k \|_2 + \text{penalty}(\gamma_k) 
 \right\}$ \\
\hline
\end{tabular}
}
\end{table}

\subsection{Simulation Setup}
For our study, we built one 3D mesh with the help of
the MRI data of the so-called Ernie subject and SimNIBS\footnote{\url{https://simnibs.github.io/simnibs/build/html/index.html}} 4 software \citep{PUONTI2020117044}. The mesh consisted of  743,575 tetrahedral elements
joined in 136,868 nodes. Four different
tissue compartments (scalp, skull, cerebrospinal fluid, and brain) were considered, and 76 electrodes were placed around the head according to the international 10-10 system. 
The lead field matrix used in this study was constructed with the help of custom-made software that exploited the Finite Element Method with linear basis functions, as in \citep{Wolters2004}. The
tissue electric conductivity values were 0.43 S/m for the scalp, 0.0103 S/m for the skull, 1.79 S/m for cerebrospinal
fluid, and 0.33 S/m for the brain (gray matter and white matter) \citep{ram06}.

The dipole source space used for reconstructions consisted of 10,000 sources distributed throughout the gray matter. To avoid the inverse crime, two different forward models of 10,000 sources were constructed so that they do not share a source at exactly the same location; however, the distance to the closest corresponding source in one model to another is at most 3 \unit{\milli\meter}. The source spaces are structured so that the sources on the grey matter layer are equally spread on 0 to 30 \unit{\milli\meter} depth from the surface of the inner skull surface. The average distance between the simulated sources and their nearest neighbors in the reconstruction space was approximately 1.3-2.5 mm.

First, we perform two validation tests. In Experiment(I), the objective is to localize a cortical source located on a sulcal wall. The source is positioned near the bottom of the sulcus to assess whether the estimators can recover activity at the true source location when sensitivity weighting is applied, or whether the estimates remain biased toward the tops of the adjacent gyri. The same forward and inversion models as described previously described, and 5\% noise is added to the simulated observations. The reconstructions are visualized by interpolating the estimated values from the source points nearest to the cutting plane onto the corresponding MRI slice. When hyperparameter updating is required, both algorithms (EM and IAS) are presented only if there is a notable difference between their estimates; otherwise, only the results obtained with the EM algorithm are shown. Experiment (II) aims to show visualize the reconstruction results in slightly deeper location...

While Experiment I and II provides qualitative insight into the spatial characteristics of the reconstructions, we next present a systematic numerical analysis to quantitatively evaluate the performance of the considered methods. This analysis focuses on localization accuracy, focality, and depth bias for dipolar sources placed at different depths relative to the inner skull surface. 
 The evaluation is structured as follows. First, we assess the overall algorithmic performance and robustness to noise using the Earth Mover’s Distance (EMD) as a measure of focality and localization accuracy. Next, we examine how reconstruction performance varies as a function of source depth by analyzing the average EMD across depths. This is followed by a statistical comparison of the methods at two selected depth ranges. Finally, we evaluate depth bias by comparing the depth of the reconstructed maximum with the true source depth. The simulated data include both low (1\%) and high (10\%) measurement noise levels.
The results are presented using histograms, depth-dependent curves, and summary statistics.

\subsection{Visualization of Source Reconstructions}
\begin{figure}
\newcommand{\mywidth}{1}
\newcommand{\TransversalSz}{0.16}
    \centering\begin{minipage}{0.2\textwidth}
        \hspace{1.5cm}\footnotesize{Transversal}
              \vspace{0.2cm}
    \end{minipage}\begin{minipage}{0.2\textwidth}
    \centering
        \footnotesize{Coronal}
              \vspace{0.2cm}
    \end{minipage}\begin{minipage}{0.2\textwidth}
    \hspace{1cm}\footnotesize{Sagittal}
              \vspace{0.2cm}
    \end{minipage}
    \begin{minipage}{0.02\textwidth}
    \rotatebox{90}{\small{\bf wMNE}}
          \end{minipage}\begin{minipage}{\TransversalSz\textwidth}
          \begin{center}
     \includegraphics[trim={2cm 1.6cm 2.1cm 3cm},clip,height=\mywidth\linewidth]{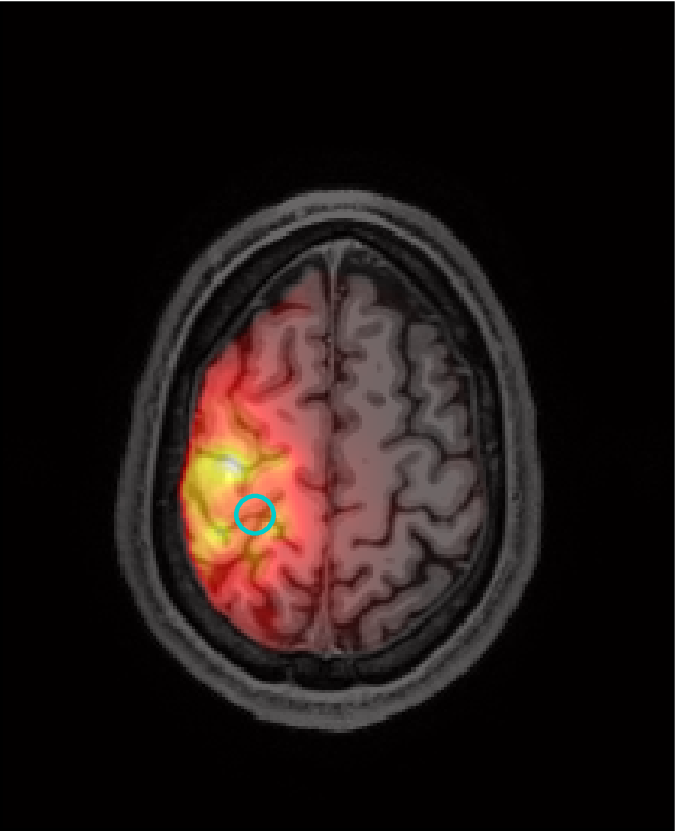}
     \end{center}
      \end{minipage}\begin{minipage}{0.18\textwidth}
          \begin{center}
     \includegraphics[trim={1cm 6cm 1.5cm 1cm},clip,height=0.75\mywidth\linewidth]{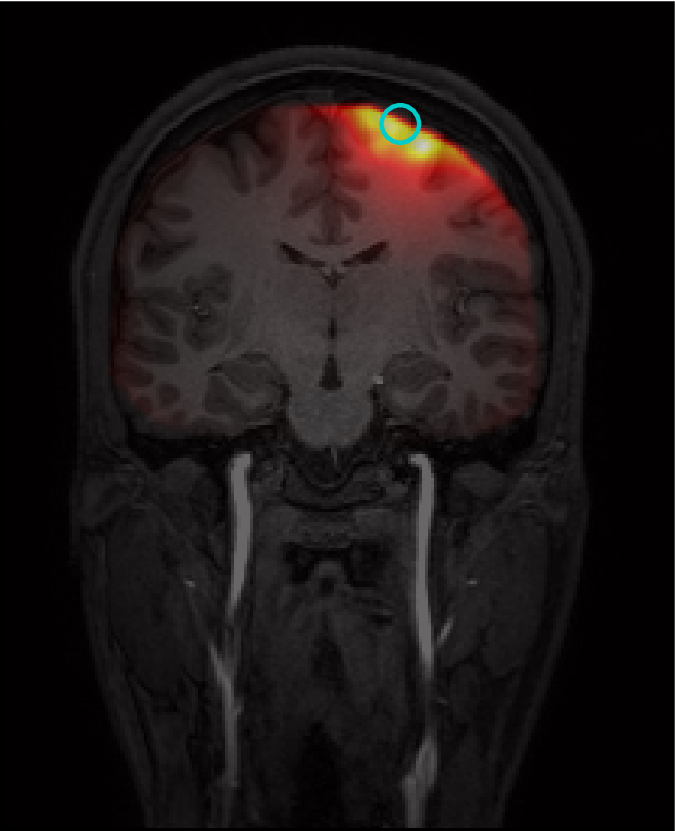}
     \end{center}
      \end{minipage}\vspace{0.5cm}\begin{minipage}{0.18\textwidth}
          \begin{center}
     \includegraphics[trim={3.5cm 5.8cm 1.8cm 2cm},clip,height=0.75\mywidth\linewidth]{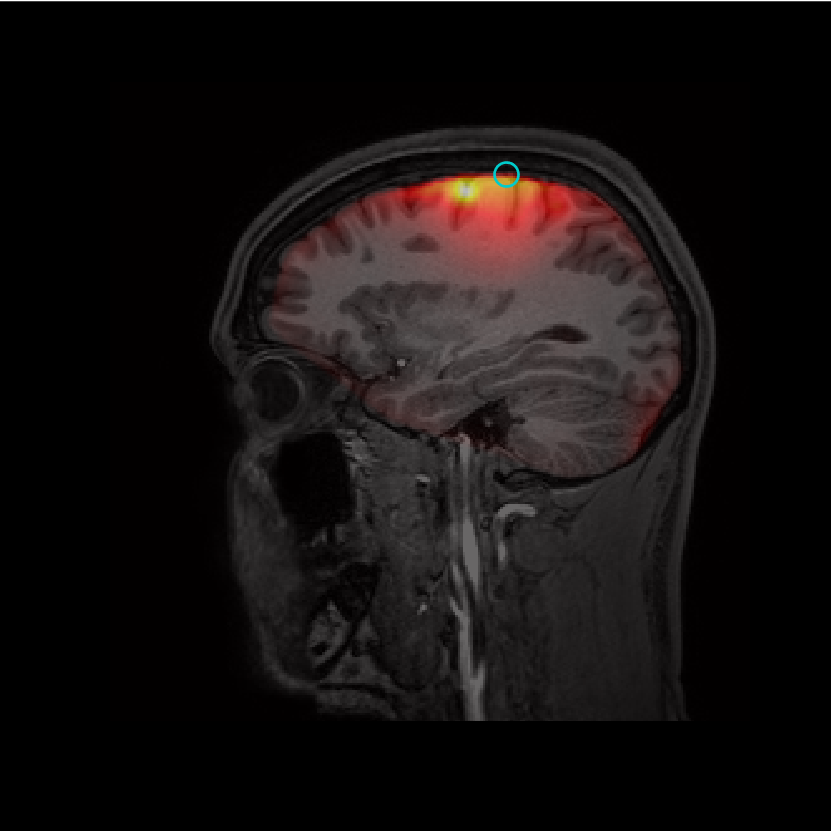}
     \end{center}
      \end{minipage}
      \begin{minipage}{0.02\textwidth}
    \rotatebox{90}{\small{\bf CG-Ga (EM)}}
          \end{minipage}\begin{minipage}{\TransversalSz\textwidth}
          \begin{center}
     \includegraphics[trim={2cm 1.6cm 2.1cm 3cm},clip,height=\mywidth\linewidth]{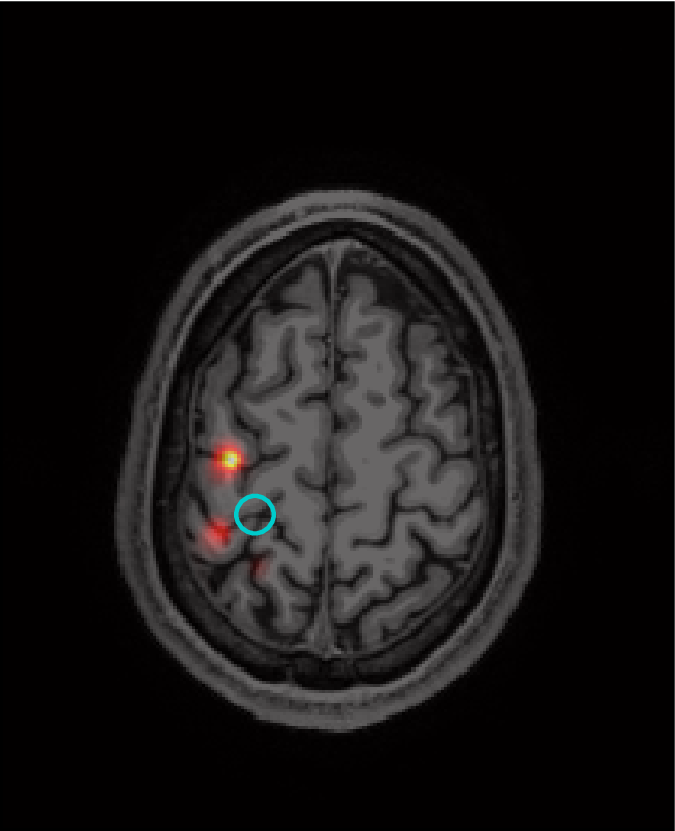}
     \end{center}
      \end{minipage}\begin{minipage}{0.18\textwidth}
          \begin{center}
     \includegraphics[trim={1cm 6cm 1.5cm 1cm},clip,height=0.75\mywidth\linewidth]{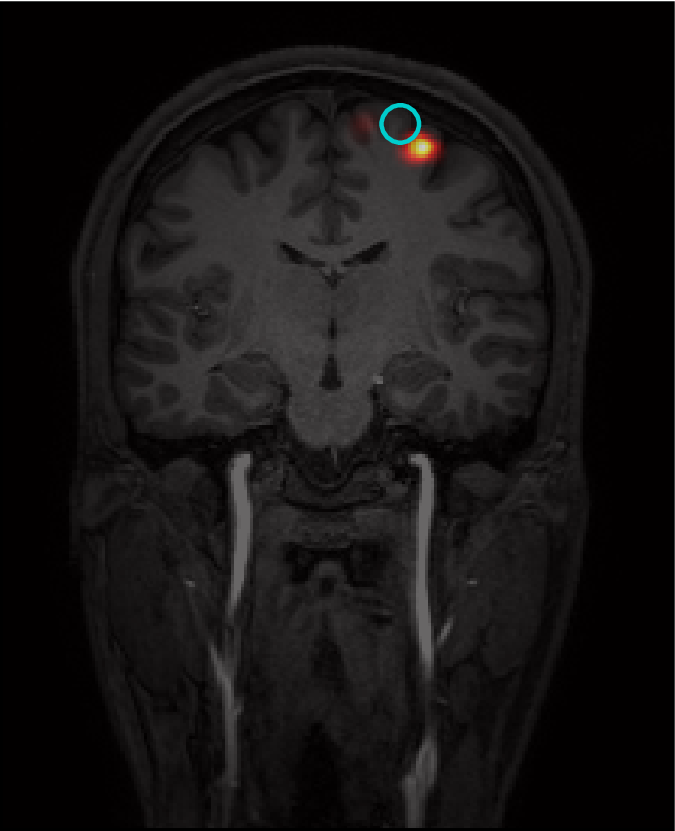}
     \end{center}
      \end{minipage}\vspace{0.5cm}\begin{minipage}{0.18\textwidth}
          \begin{center}
     \includegraphics[trim={3.5cm 5.8cm 1.8cm 2cm},clip,height=0.75\mywidth\linewidth]{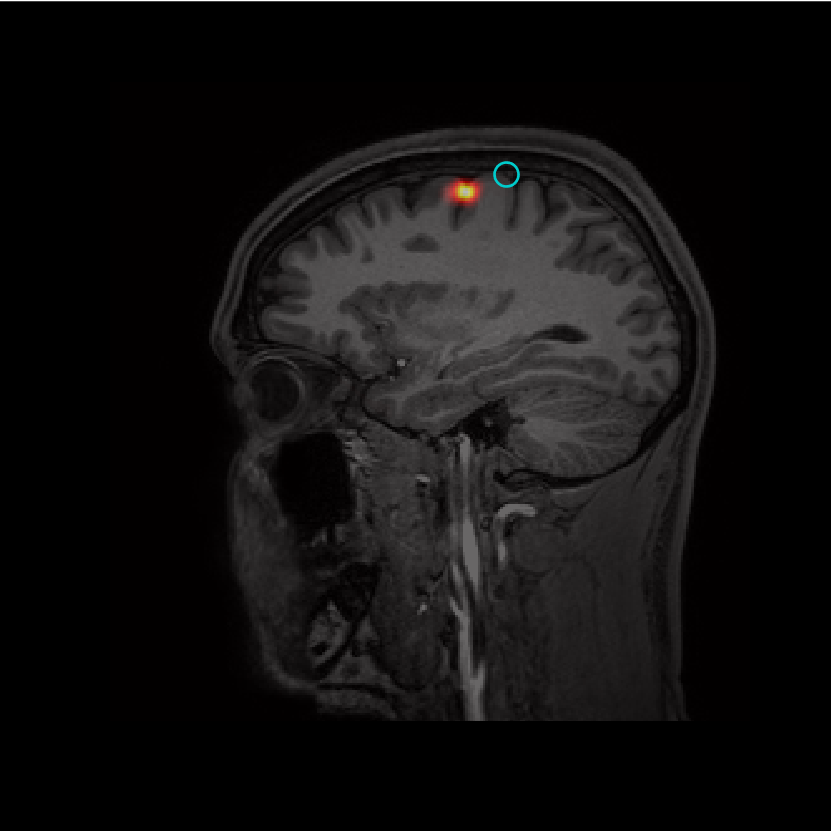}
     \end{center}
      \end{minipage}
      \begin{minipage}{0.02\textwidth}
    \rotatebox{90}{\small{\bf CG-Ga (IAS)}}
          \end{minipage}\begin{minipage}{\TransversalSz\textwidth}
          \begin{center}
     \includegraphics[trim={2cm 1.6cm 2.1cm 3cm},clip,height=\mywidth\linewidth]{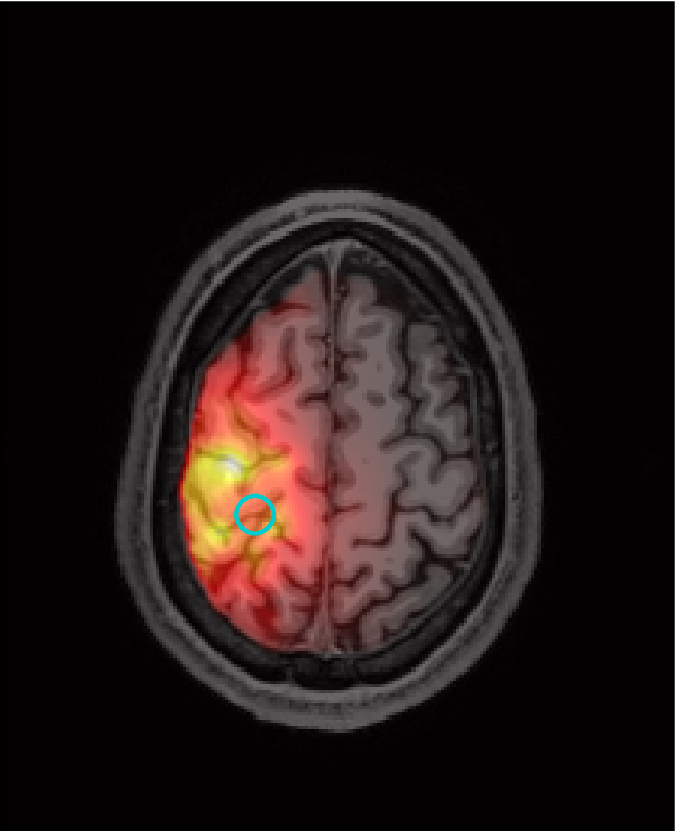}
     \end{center}
      \end{minipage}\begin{minipage}{0.18\textwidth}
          \begin{center}
     \includegraphics[trim={1cm 6cm 1.5cm 1cm},clip,height=0.75\mywidth\linewidth]{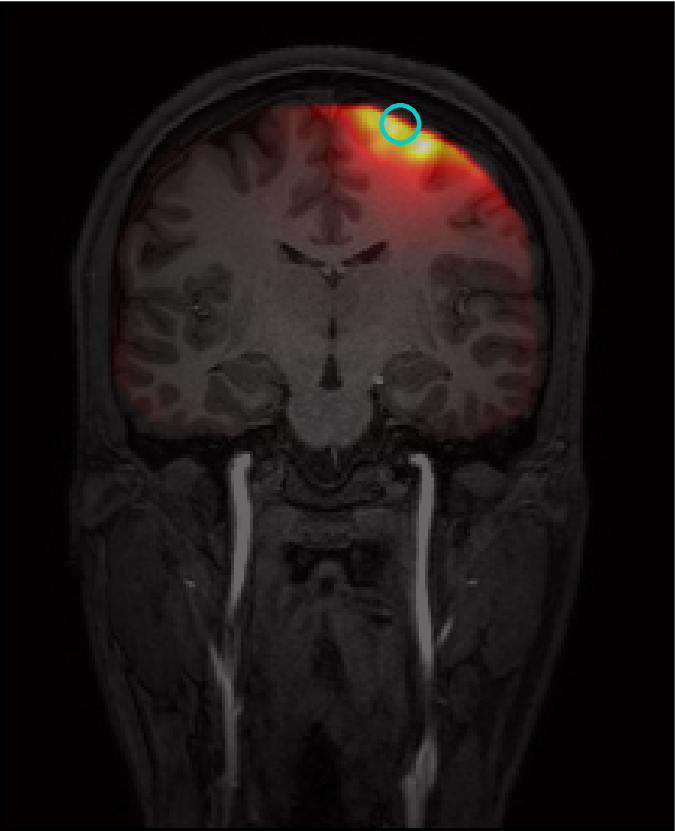}
     \end{center}
      \end{minipage}\vspace{0.5cm}\begin{minipage}{0.18\textwidth}
          \begin{center}
     \includegraphics[trim={3.5cm 5.8cm 1.8cm 2cm},clip,height=0.75\mywidth\linewidth]{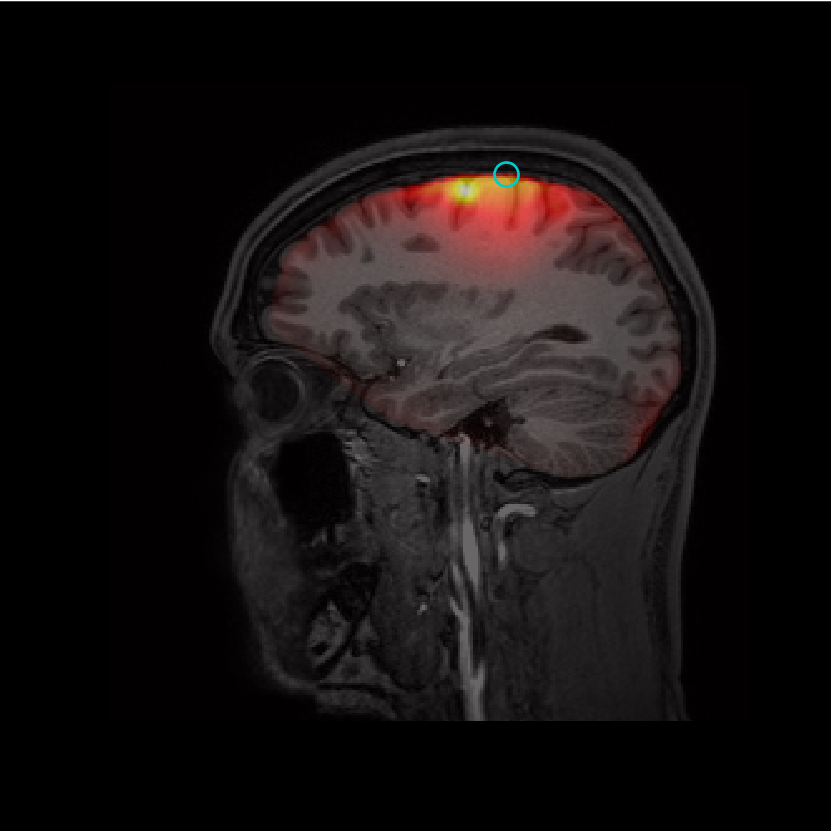}
     \end{center}
      \end{minipage}

      \begin{minipage}{0.02\textwidth}
    \rotatebox{90}{\small{\bf CG-IG (EM)}}
          \end{minipage}\begin{minipage}{\TransversalSz\textwidth}
          \begin{center}
     \includegraphics[trim={2cm 1.6cm 2.1cm 3cm},clip,height=\mywidth\linewidth]{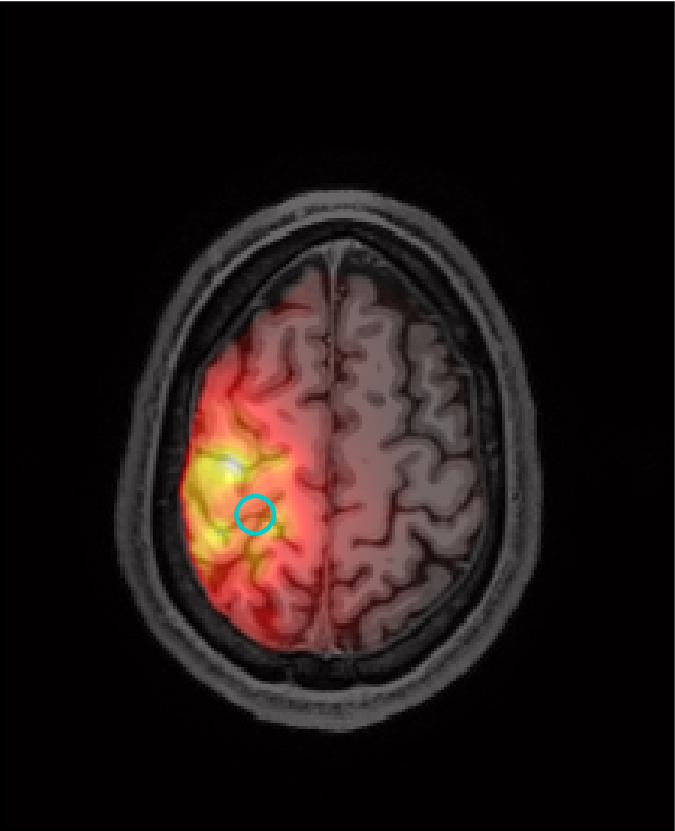}
     \end{center}
      \end{minipage}\begin{minipage}{0.18\textwidth}
          \begin{center}
     \includegraphics[trim={1cm 6cm 1.5cm 1cm},clip,height=0.75\mywidth\linewidth]{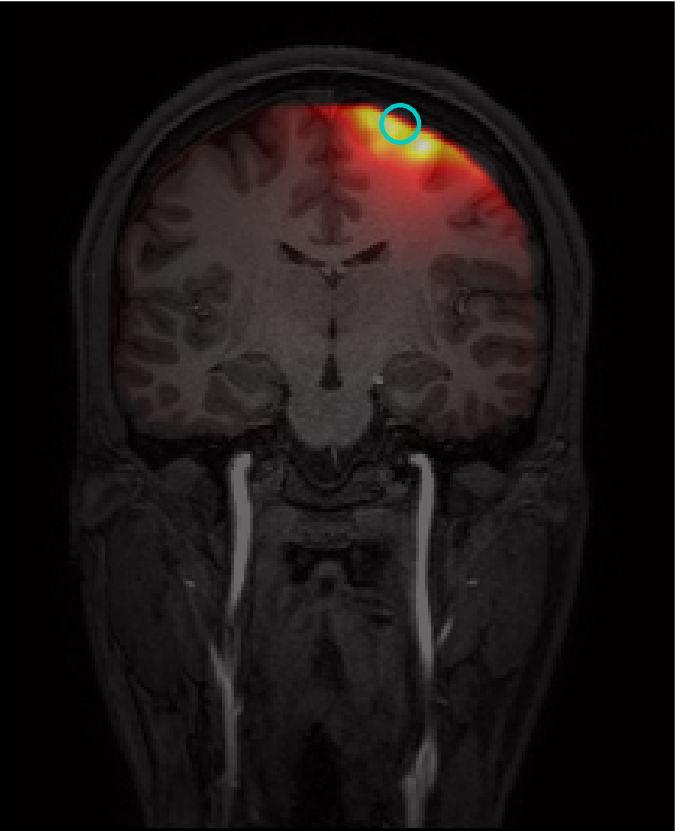}
     \end{center}
      \end{minipage}\vspace{0.5cm}\begin{minipage}{0.18\textwidth}
          \begin{center}
     \includegraphics[trim={3.5cm 5.8cm 1.8cm 2cm},clip,height=0.75\mywidth\linewidth]{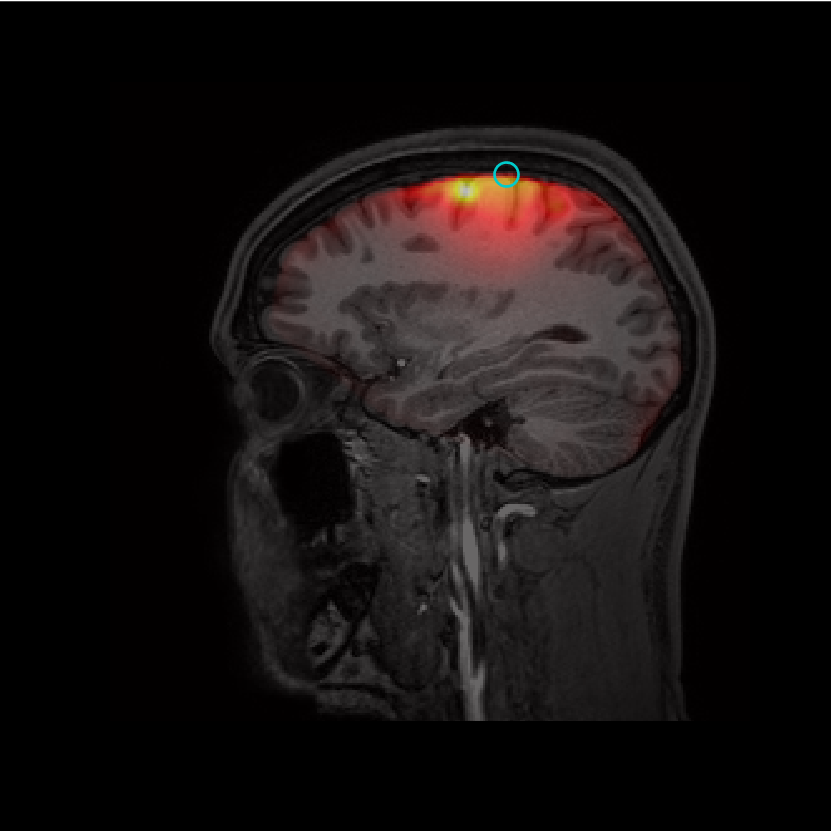}
     \end{center}
      \end{minipage}

      \begin{minipage}{0.02\textwidth}
    \rotatebox{90}{\small{\bf CG-IG (IAS)}}
          \end{minipage}\begin{minipage}{\TransversalSz\textwidth}
          \begin{center}
     \includegraphics[trim={2cm 1.6cm 2.1cm 3cm},clip,height=\mywidth\linewidth]{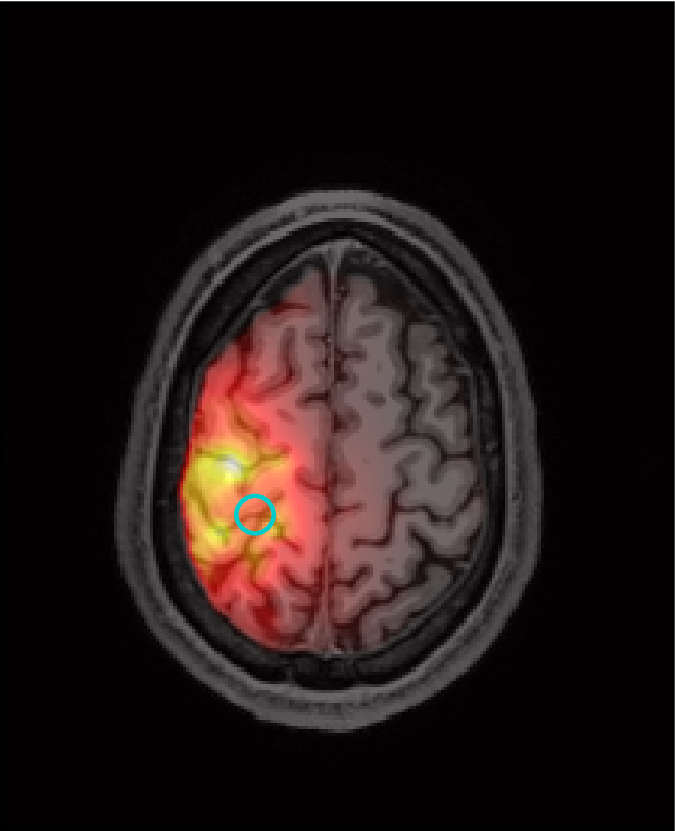}
     \end{center}
      \end{minipage}\begin{minipage}{0.18\textwidth}
          \begin{center}
     \includegraphics[trim={1cm 6cm 1.5cm 1cm},clip,height=0.75\mywidth\linewidth]{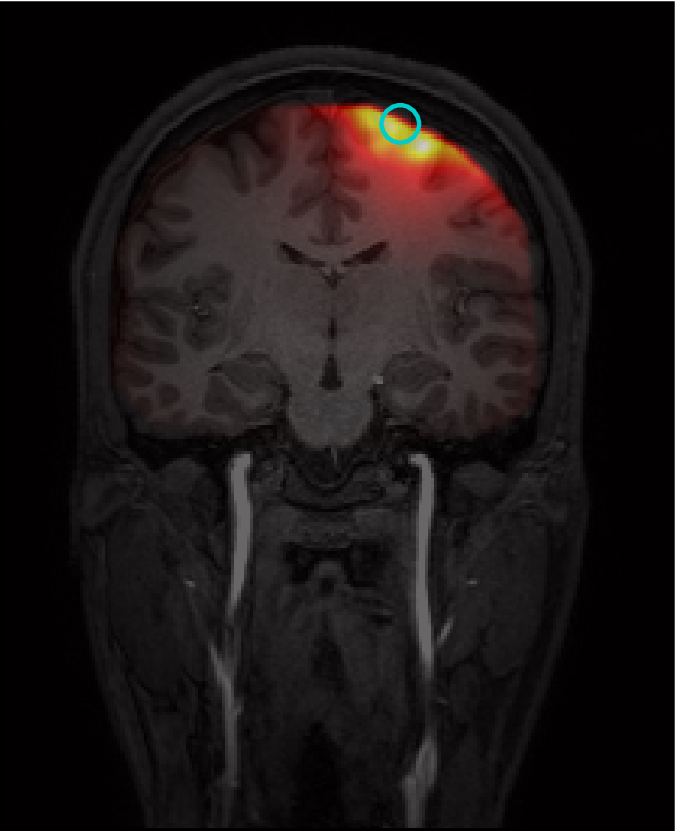}
     \end{center}
      \end{minipage}\vspace{0.5cm}\begin{minipage}{0.18\textwidth}
          \begin{center}
     \includegraphics[trim={3.5cm 5.8cm 1.8cm 2cm},clip,height=0.75\mywidth\linewidth]{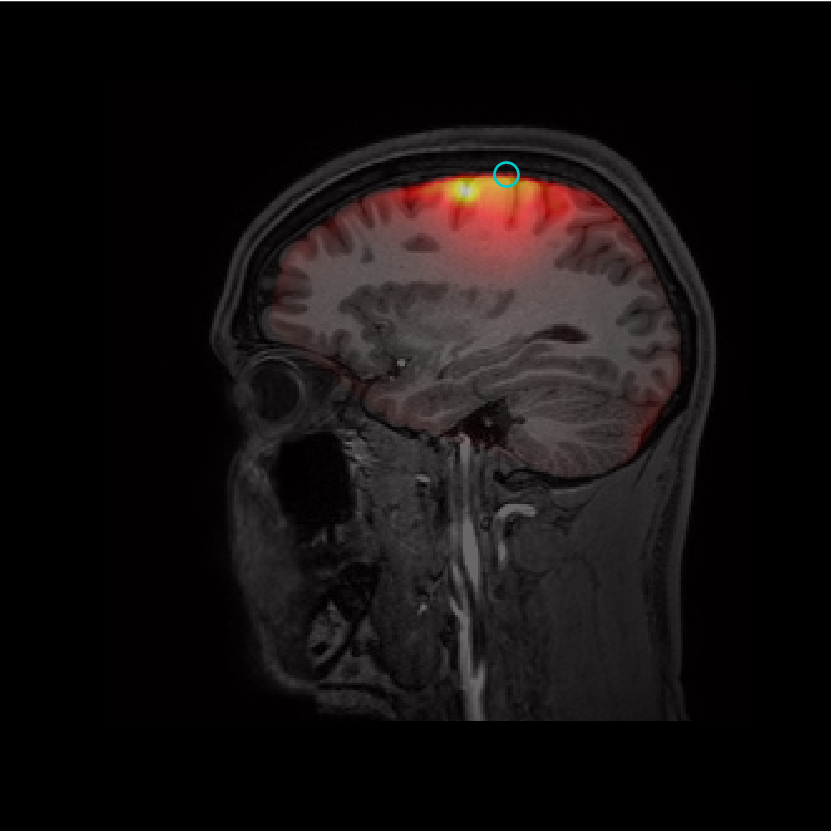}
     \end{center}
      \end{minipage}
    \caption{Estimated distributions of the simulated superficial brain activity computed using methods with Gaussian priors. The distributions are presented in three plain cuts of magnetic resonance images. The turquoise ring shows the location of the actual source to be estimated. The colored region, ranging from dark red to yellow, represents the distribution and its local magnitude. Slices have been taken at the location of the maximum estimated magnitude.}
    \label{fig:MRI2Gauss_superf}
\end{figure}

\begin{figure}
\newcommand{\mywidth}{1}
\newcommand{\TransversalSz}{0.16}
    \centering
    \begin{minipage}{0.02\textwidth}
    \rotatebox{90}{\small{\bf wL}}
          \end{minipage}\begin{minipage}{\TransversalSz\textwidth}
          \begin{center}
     \includegraphics[trim={2cm 1.6cm 2.1cm 3cm},clip,height=\mywidth\linewidth]{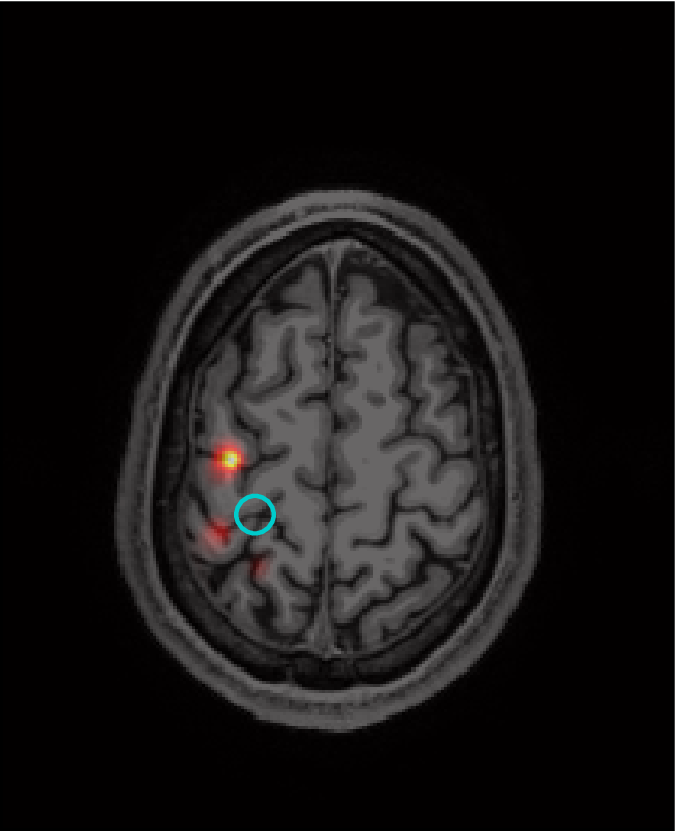}
     \end{center}
      \end{minipage}\begin{minipage}{0.18\textwidth}
          \begin{center}
     \includegraphics[trim={1cm 6cm 1.5cm 1cm},clip,height=0.75\mywidth\linewidth]{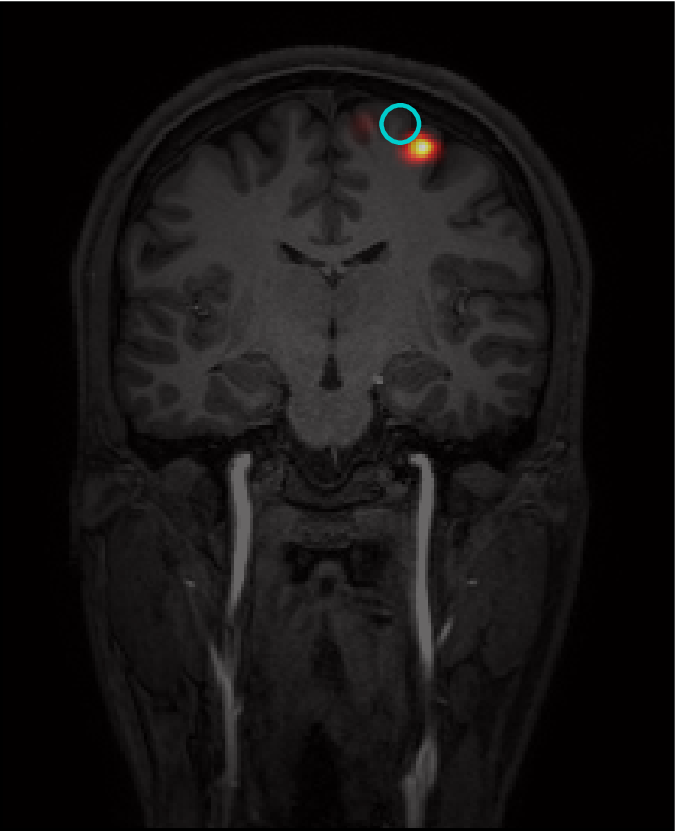}
     \end{center}
      \end{minipage}\vspace{0.5cm}\begin{minipage}{0.18\textwidth}
          \begin{center}
     \includegraphics[trim={3.5cm 5.8cm 1.8cm 2cm},clip,height=0.75\mywidth\linewidth]{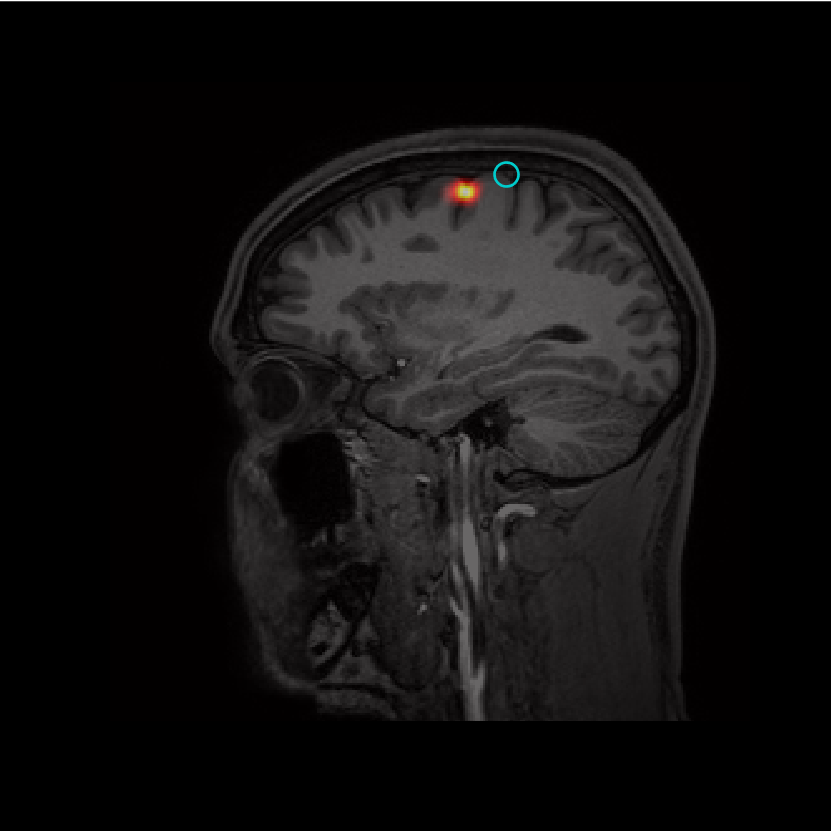}
     \end{center}
      \end{minipage}
      \begin{minipage}{0.02\textwidth}
    \rotatebox{90}{\small{\bf wGL}}
          \end{minipage}\begin{minipage}{\TransversalSz\textwidth}
          \begin{center}
     \includegraphics[trim={2cm 1.6cm 2.1cm 3cm},clip,height=\mywidth\linewidth]{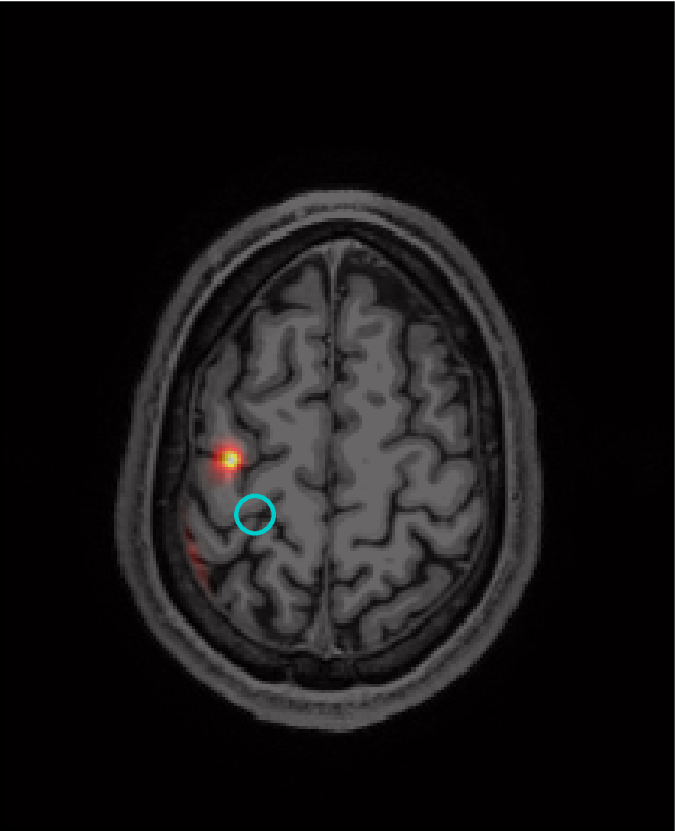}
     \end{center}
      \end{minipage}\begin{minipage}{0.18\textwidth}
          \begin{center}
     \includegraphics[trim={1cm 6cm 1.5cm 1cm},clip,height=0.75\mywidth\linewidth]{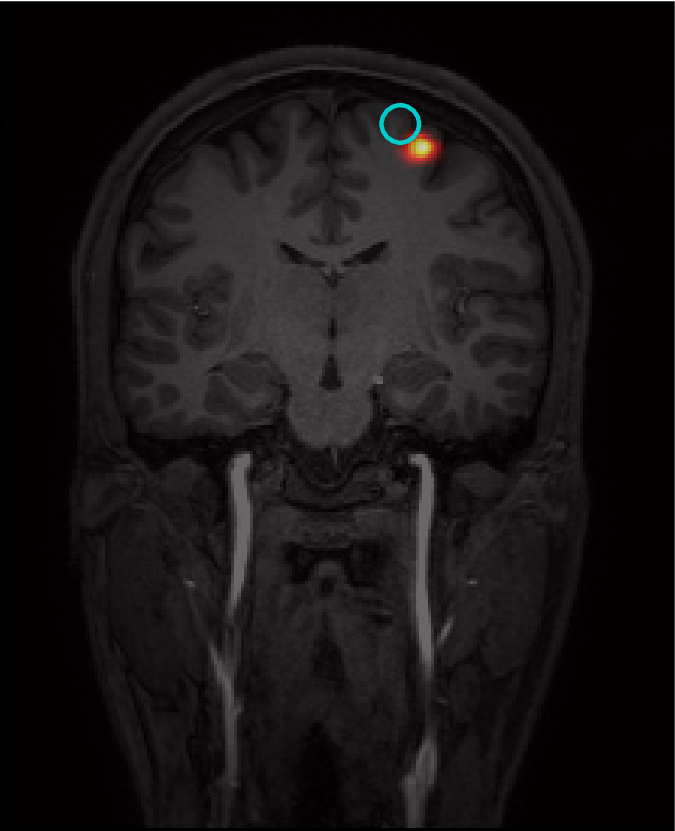}
     \end{center}
      \end{minipage}\vspace{0.5cm}\begin{minipage}{0.18\textwidth}
          \begin{center}
     \includegraphics[trim={3.5cm 5.8cm 1.8cm 2cm},clip,height=0.75\mywidth\linewidth]{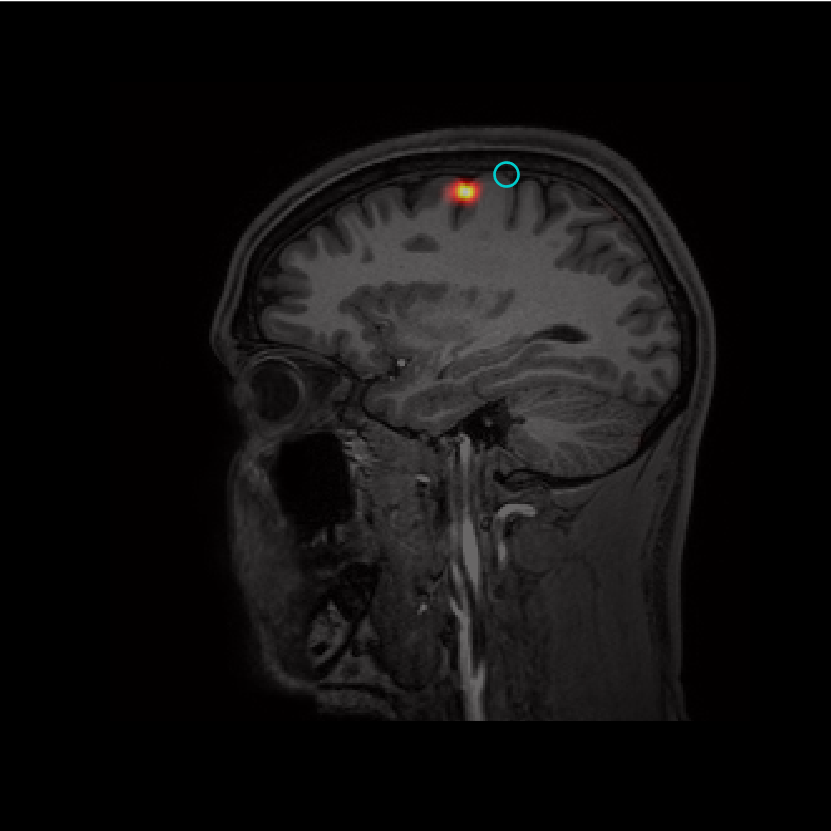}
     \end{center}
      \end{minipage}
      \begin{minipage}{0.02\textwidth}
    \rotatebox{90}{\small{\bf wCL (EM)}}
          \end{minipage}\begin{minipage}{\TransversalSz\textwidth}
          \begin{center}
     \includegraphics[trim={2cm 1.6cm 2.1cm 3cm},clip,height=\mywidth\linewidth]{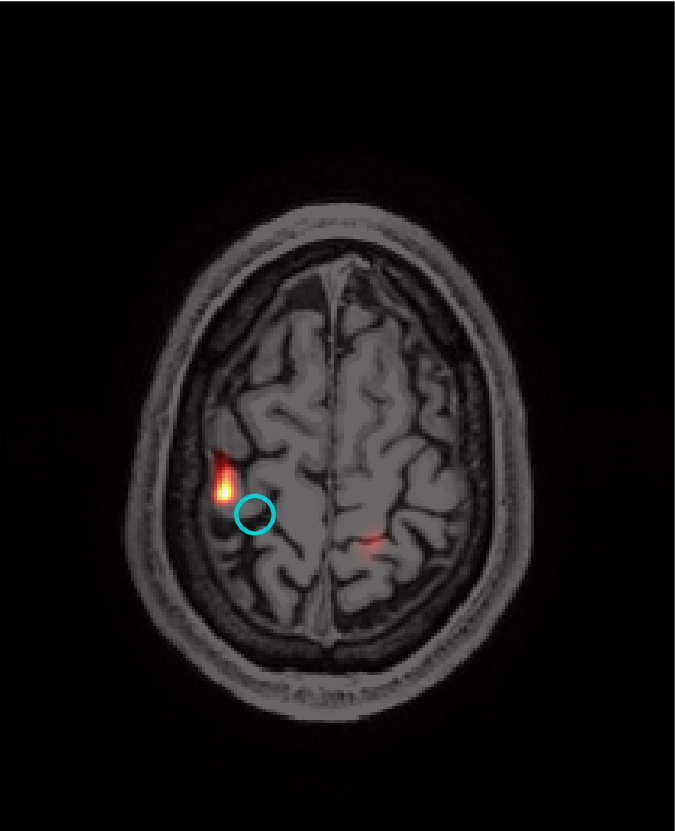}
     \end{center}
      \end{minipage}\begin{minipage}{0.18\textwidth}
          \begin{center}
     \includegraphics[trim={1cm 6cm 1.5cm 1cm},clip,height=0.75\mywidth\linewidth]{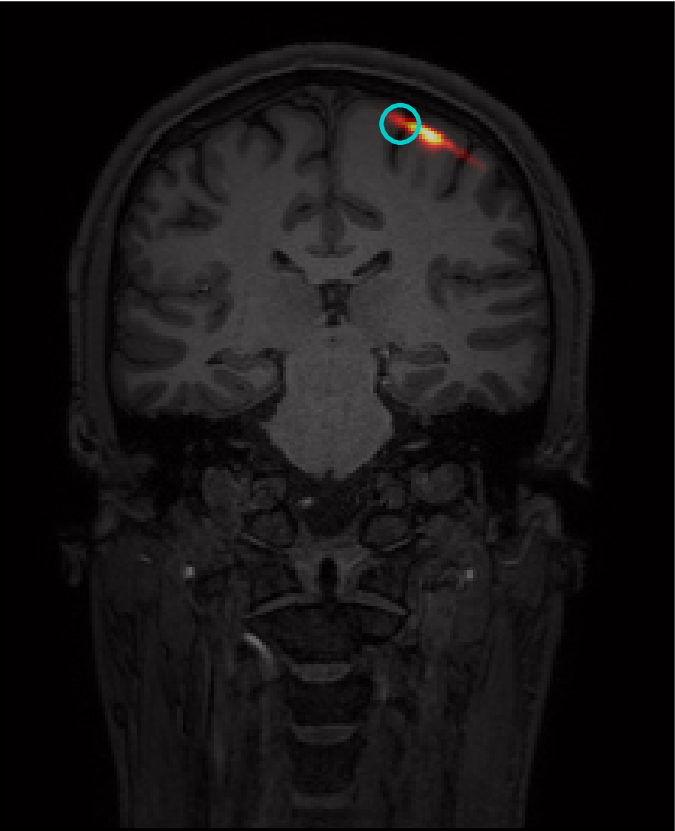}
     \end{center}
      \end{minipage}\vspace{0.5cm}\begin{minipage}{0.18\textwidth}
          \begin{center}
     \includegraphics[trim={3.5cm 5.8cm 1.8cm 2cm},clip,height=0.75\mywidth\linewidth]{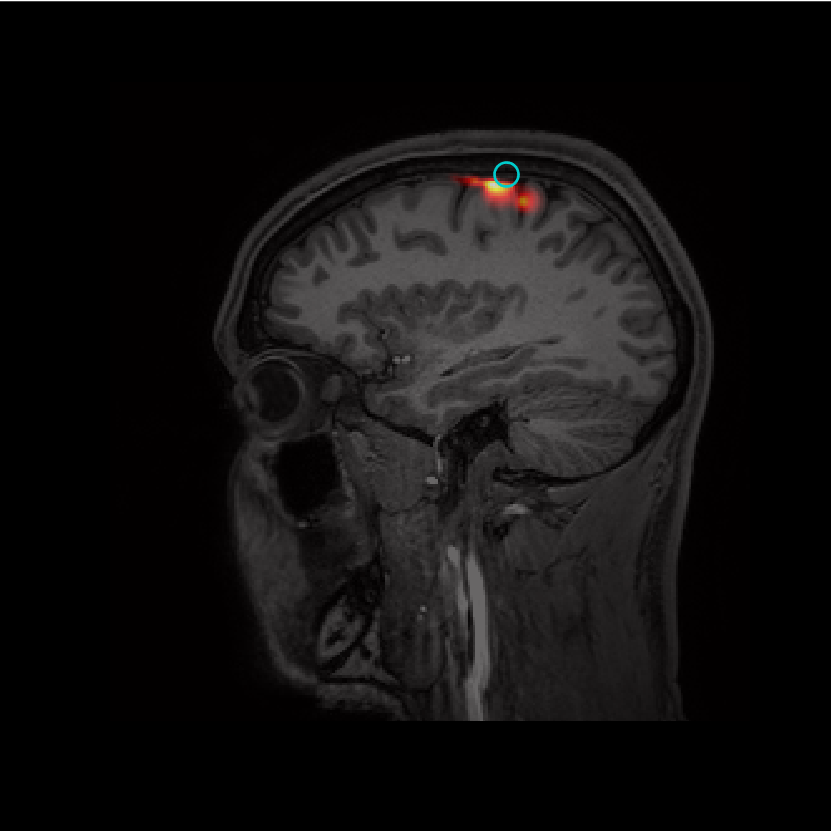}
     \end{center}
      \end{minipage}
      \begin{minipage}{0.02\textwidth}
    \rotatebox{90}{\small{\bf wCL (IAS)}}
          \end{minipage}\begin{minipage}{\TransversalSz\textwidth}
          \begin{center}
     \includegraphics[trim={2cm 1.6cm 2.1cm 3cm},clip,height=\mywidth\linewidth]{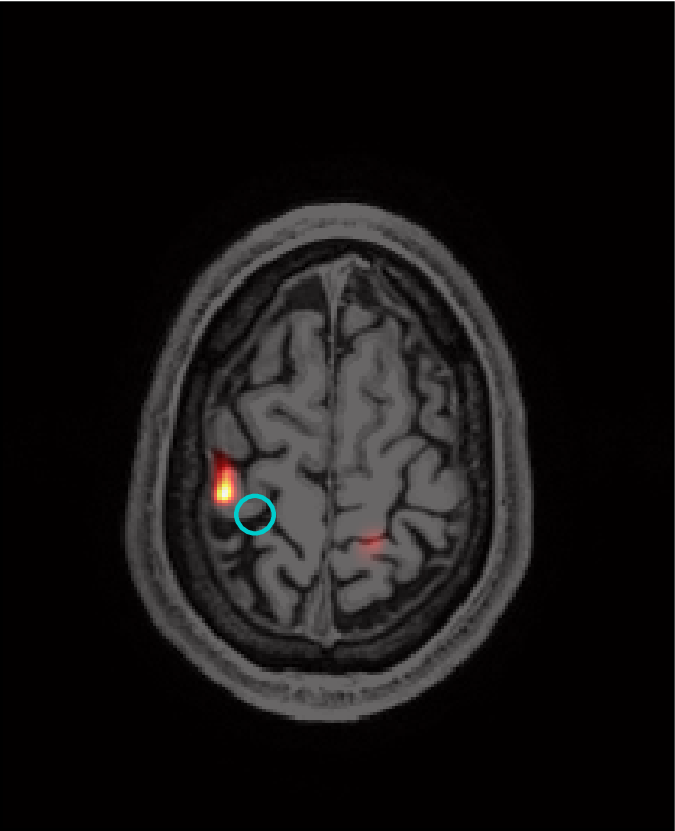}
     \end{center}
      \end{minipage}\begin{minipage}{0.18\textwidth}
          \begin{center}
     \includegraphics[trim={1cm 6cm 1.5cm 1cm},clip,height=0.75\mywidth\linewidth]{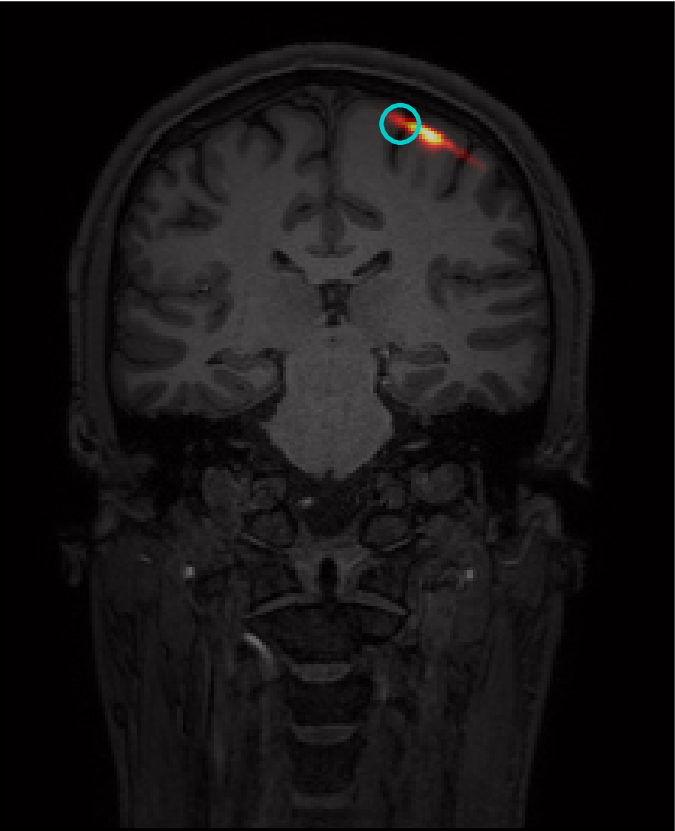}
     \end{center}
      \end{minipage}\vspace{0.5cm}\begin{minipage}{0.18\textwidth}
          \begin{center}
     \includegraphics[trim={3.5cm 5.8cm 1.8cm 2cm},clip,height=0.75\mywidth\linewidth]{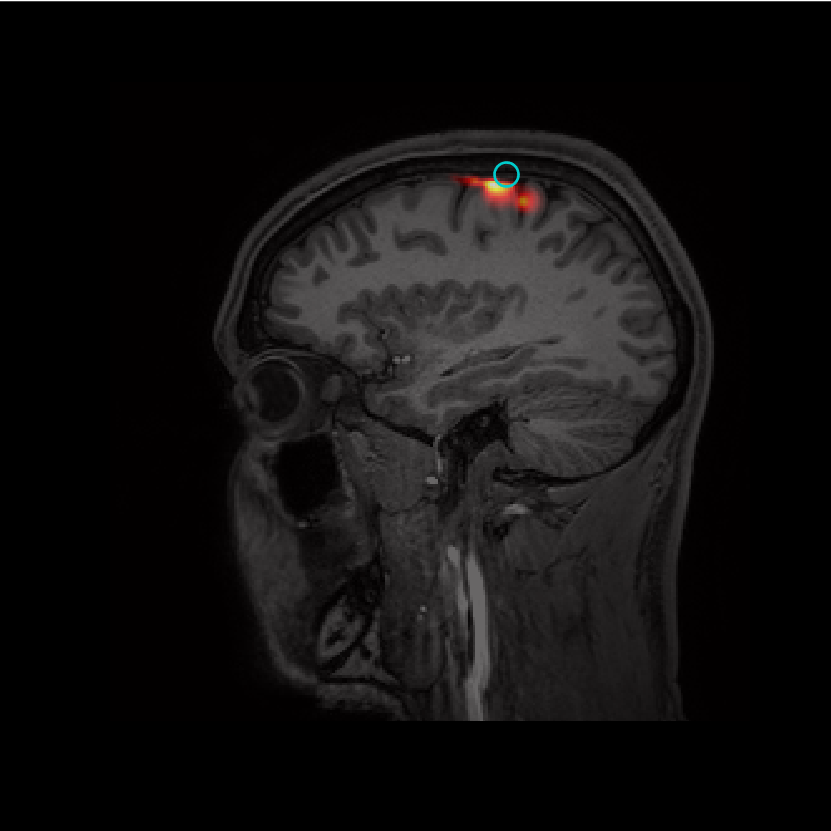}
     \end{center}
      \end{minipage}
      \begin{minipage}{0.02\textwidth}
    \rotatebox{90}{\small{\bf wCGL (EM)}}
          \end{minipage}\begin{minipage}{\TransversalSz\textwidth}
          \begin{center}
     \includegraphics[trim={2cm 1.6cm 2.1cm 3cm},clip,height=\mywidth\linewidth]{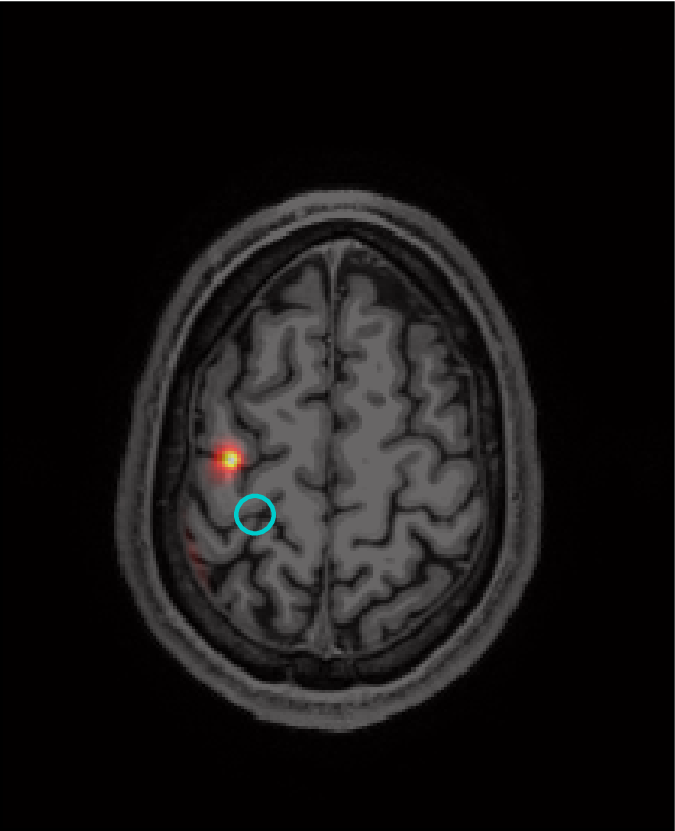}
     \end{center}
      \end{minipage}\begin{minipage}{0.18\textwidth}
          \begin{center}
     \includegraphics[trim={1cm 6cm 1.5cm 1cm},clip,height=0.75\mywidth\linewidth]{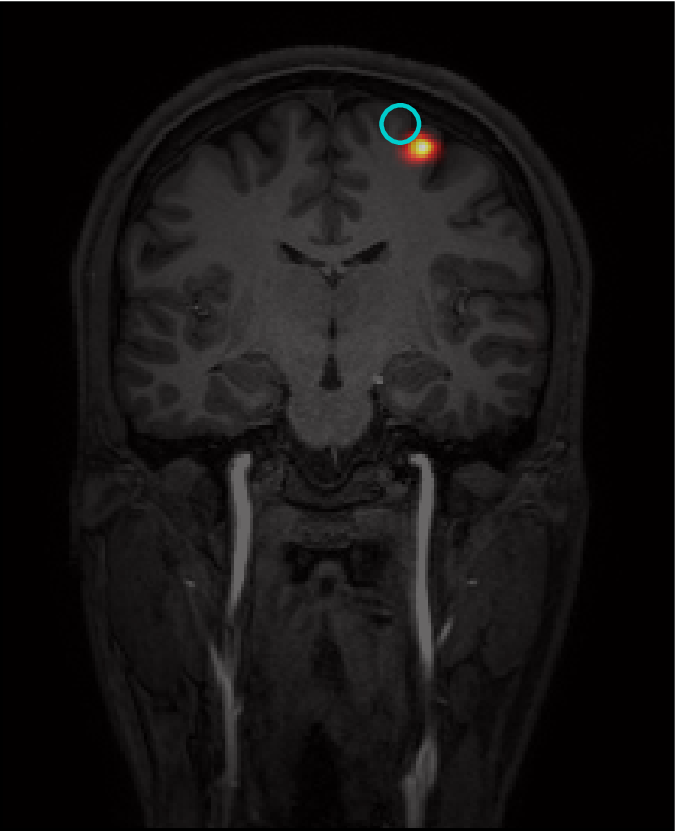}
     \end{center}
      \end{minipage}\vspace{0.5cm}\begin{minipage}{0.18\textwidth}
          \begin{center}
     \includegraphics[trim={3.5cm 5.8cm 1.8cm 2cm},clip,height=0.75\mywidth\linewidth]{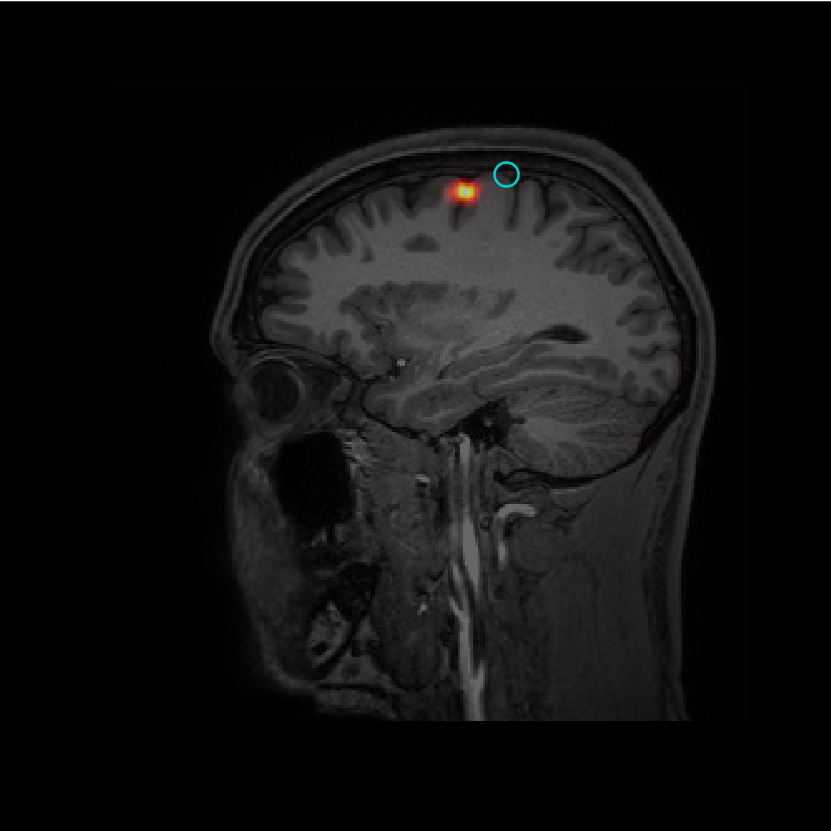}
     \end{center}
      \end{minipage}
      \begin{minipage}{0.02\textwidth}
    \rotatebox{90}{\small{\bf wCGL (IAS)}}
          \end{minipage}\begin{minipage}{\TransversalSz\textwidth}
          \begin{center}
     \includegraphics[trim={2cm 1.6cm 2.1cm 3cm},clip,height=\mywidth\linewidth]{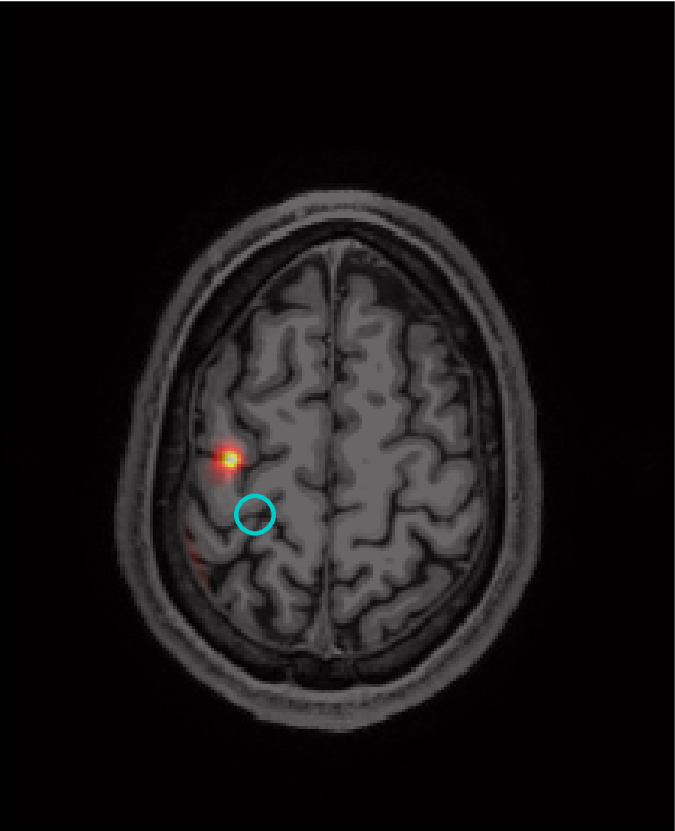}
     \end{center}
      \end{minipage}\begin{minipage}{0.18\textwidth}
          \begin{center}
     \includegraphics[trim={1cm 6cm 1.5cm 1cm},clip,height=0.75\mywidth\linewidth]{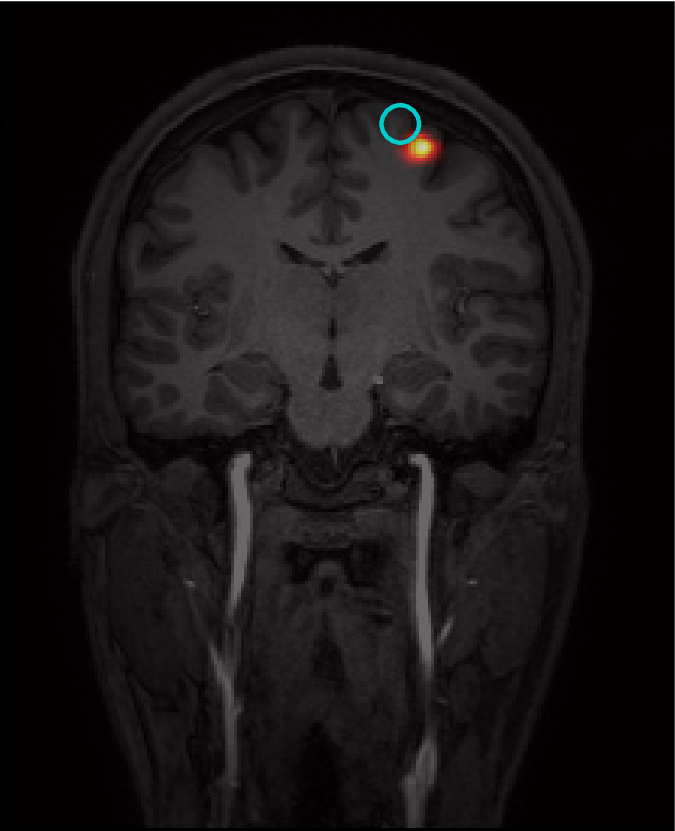}
     \end{center}
      \end{minipage}\vspace{0.5cm}\begin{minipage}{0.18\textwidth}
          \begin{center}
     \includegraphics[trim={3.5cm 5.8cm 1.8cm 2cm},clip,height=0.75\mywidth\linewidth]{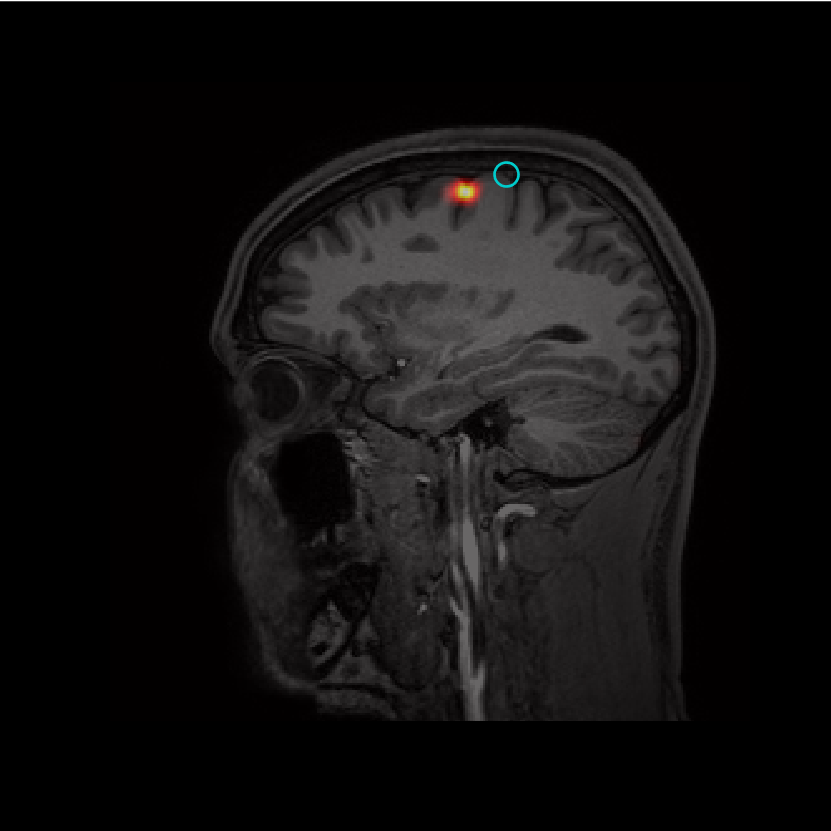}
     \end{center}
      \end{minipage}
    \caption{Estimated distributions of the simulated superficial brain activity computed using methods with Laplace priors. The distributions are presented in three plain cuts of magnetic resonance images. The turquoise ring shows the location of the actual source to be estimated. The colored region, ranging from dark red to yellow, represents the distribution and its local magnitude. Slices have been taken at the location of the maximum estimated magnitude.}
    \label{fig:MRI2Lap_superf}
\end{figure}

\begin{figure}
\newcommand{\mywidth}{1}
\newcommand{\TransversalSz}{0.16}
    \centering\begin{minipage}{0.2\textwidth}
        \hspace{1.5cm}\footnotesize{Transversal}
              \vspace{0.2cm}
    \end{minipage}\begin{minipage}{0.2\textwidth}
    \centering
        \footnotesize{Coronal}
              \vspace{0.2cm}
    \end{minipage}\begin{minipage}{0.2\textwidth}
    \hspace{1cm}\footnotesize{Sagittal}
              \vspace{0.2cm}
    \end{minipage}
    \begin{minipage}{0.02\textwidth}
    \rotatebox{90}{\small{\bf wMNE}}
          \end{minipage}\begin{minipage}{\TransversalSz\textwidth}
          \begin{center}
     \includegraphics[trim={2cm 1.6cm 2.1cm 3cm},clip,height=\mywidth\linewidth]{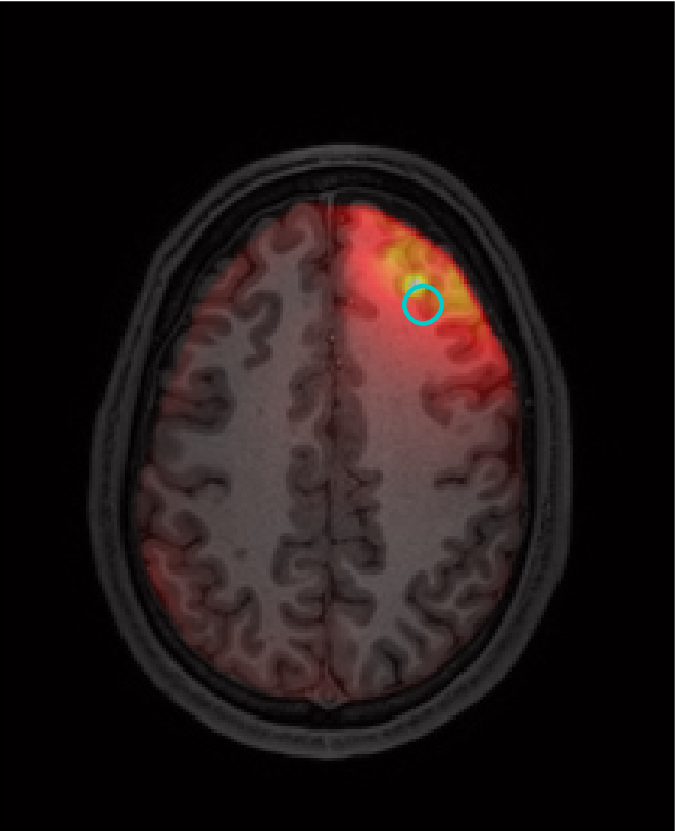}
     \end{center}
      \end{minipage}\begin{minipage}{0.18\textwidth}
          \begin{center}
     \includegraphics[trim={1cm 6cm 1.5cm 1cm},clip,height=0.75\mywidth\linewidth]{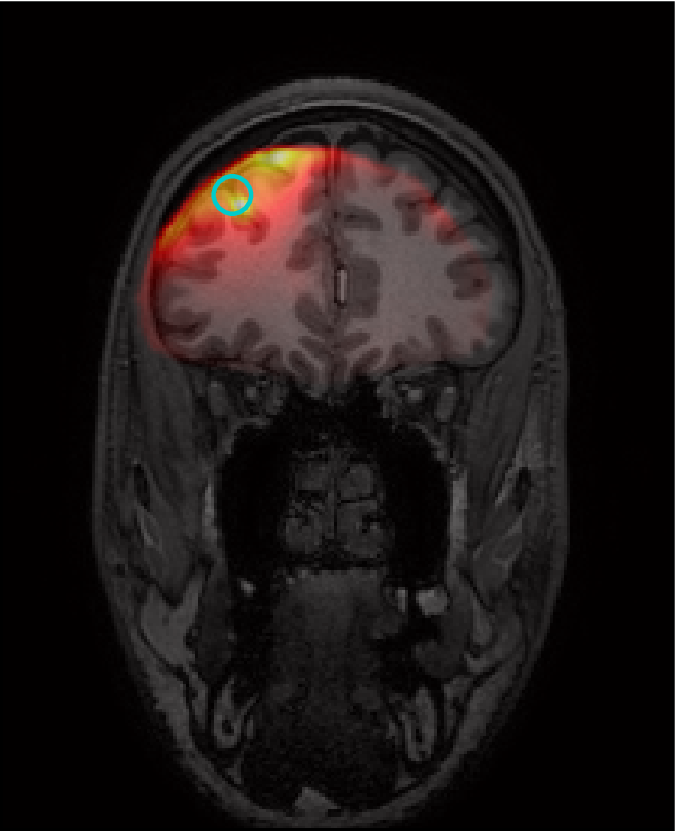}
     \end{center}
      \end{minipage}\vspace{0.5cm}\begin{minipage}{0.18\textwidth}
          \begin{center}
     \includegraphics[trim={10.5cm 7.5cm 8cm 2.8cm},clip,height=0.75\mywidth\linewidth]{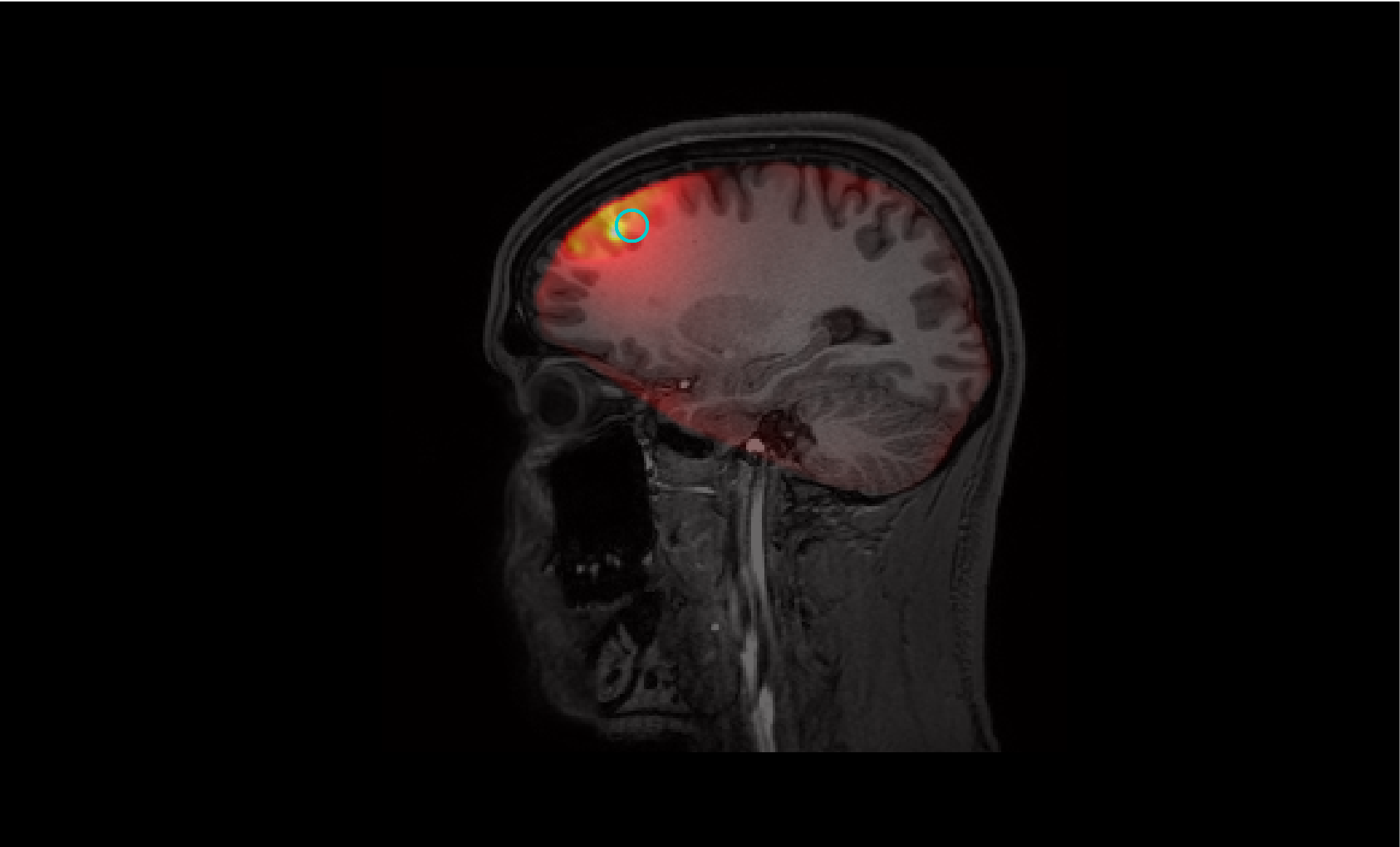}
     \end{center}
      \end{minipage}
      \begin{minipage}{0.02\textwidth}
    \rotatebox{90}{\small{\bf CG-Ga (EM)}}
          \end{minipage}\begin{minipage}{\TransversalSz\textwidth}
          \begin{center}
     \includegraphics[trim={2cm 1.6cm 2.1cm 3cm},clip,height=\mywidth\linewidth]{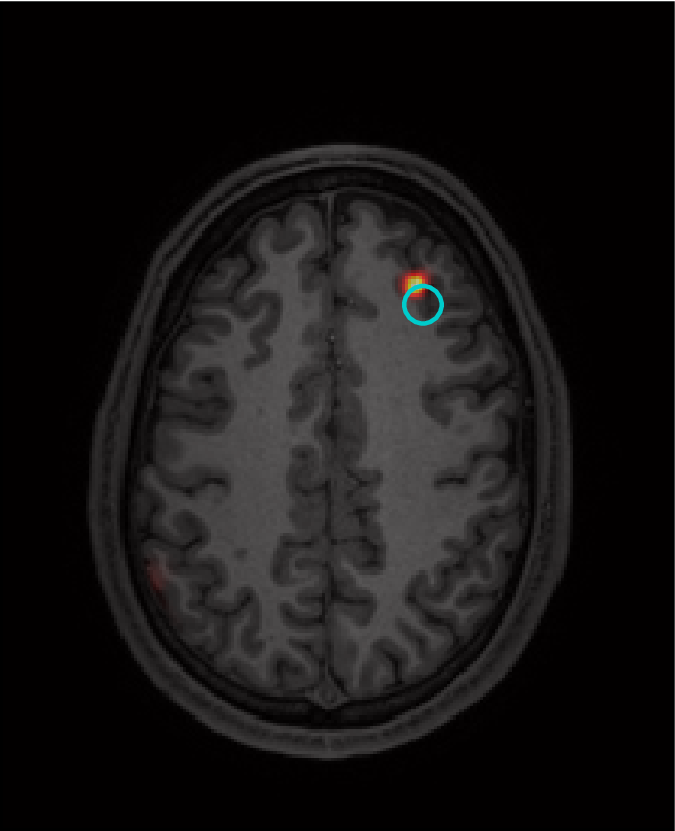}
     \end{center}
      \end{minipage}\begin{minipage}{0.18\textwidth}
          \begin{center}
     \includegraphics[trim={1cm 6cm 1.5cm 1cm},clip,height=0.75\mywidth\linewidth]{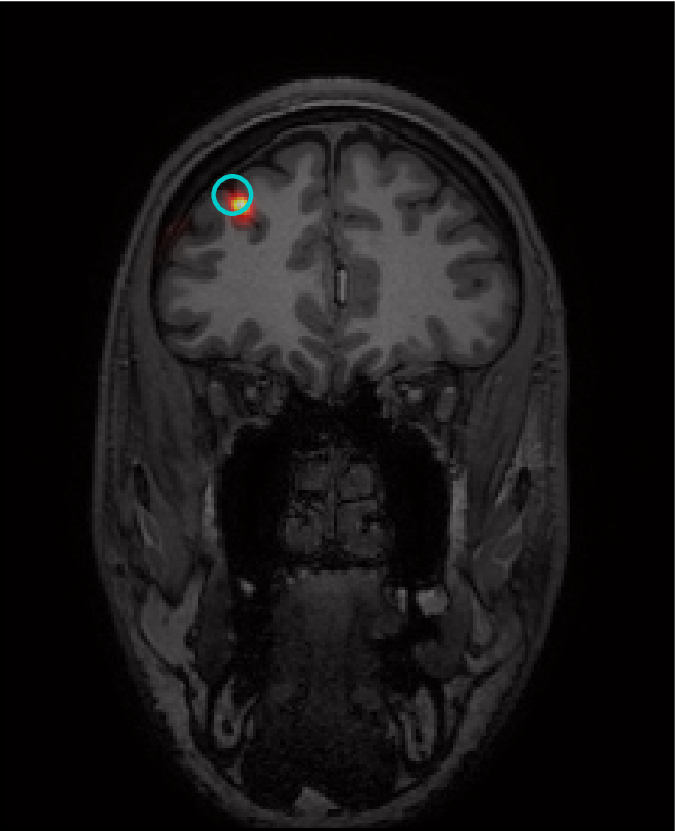}
     \end{center}
      \end{minipage}\vspace{0.5cm}\begin{minipage}{0.18\textwidth}
          \begin{center}
     \includegraphics[trim={10.5cm 7.5cm 8cm 2.8cm},clip,height=0.75\mywidth\linewidth]{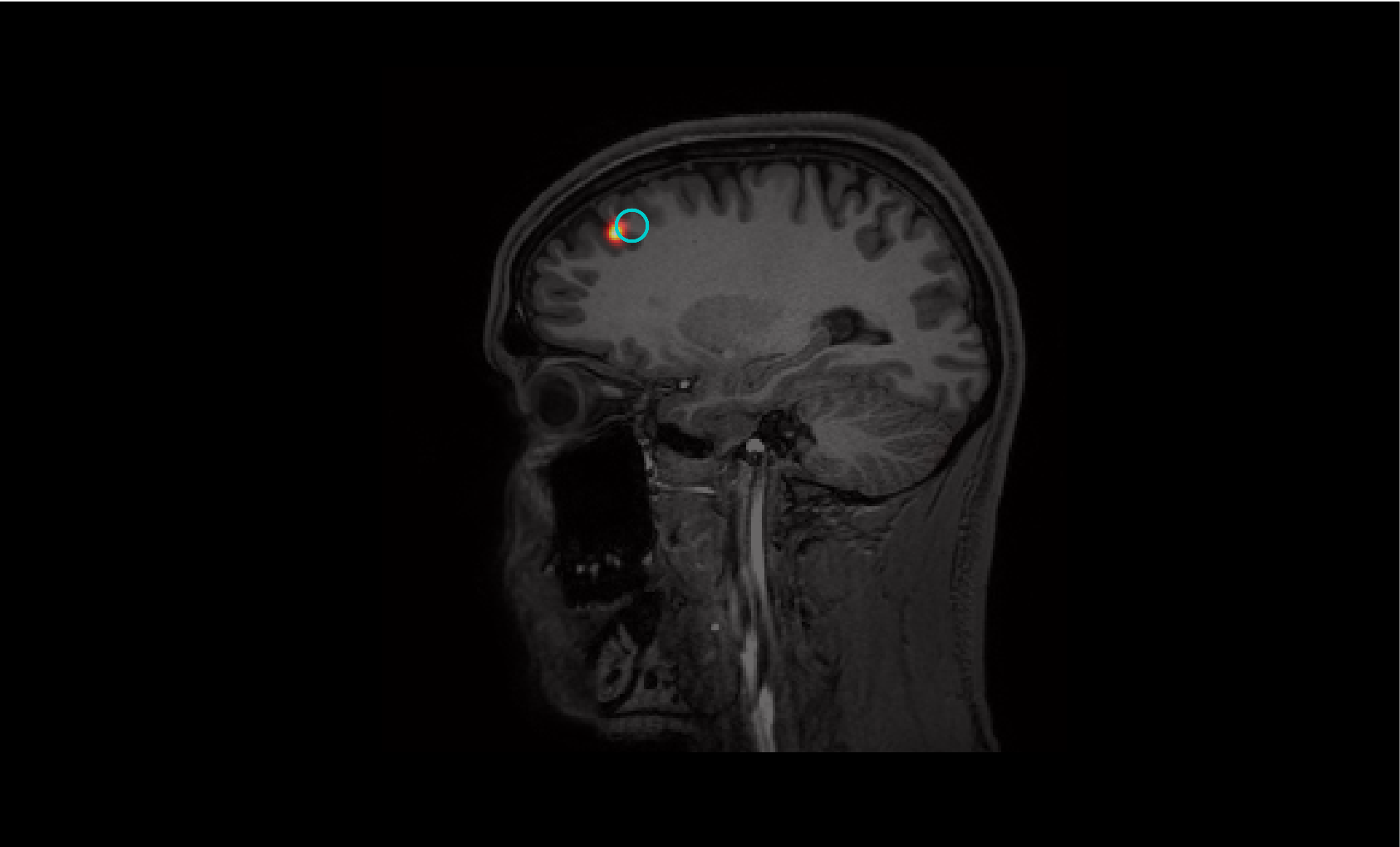}
     \end{center}
      \end{minipage}
      \begin{minipage}{0.02\textwidth}
    \rotatebox{90}{\small{\bf CG-Ga (IAS)}}
          \end{minipage}\begin{minipage}{\TransversalSz\textwidth}
          \begin{center}
     \includegraphics[trim={2cm 1.6cm 2.1cm 3cm},clip,height=\mywidth\linewidth]{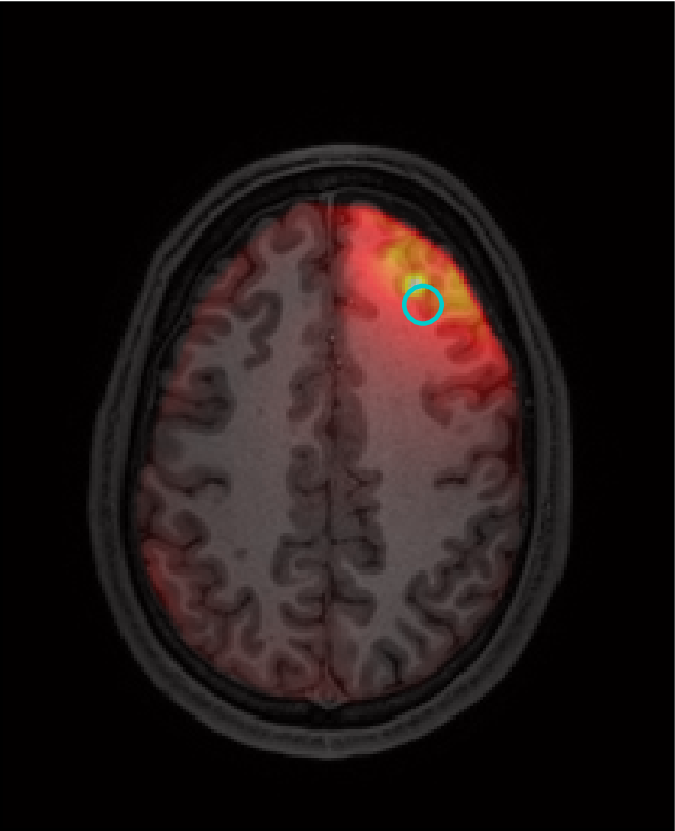}
     \end{center}
      \end{minipage}\begin{minipage}{0.18\textwidth}
          \begin{center}
     \includegraphics[trim={1cm 6cm 1.5cm 1cm},clip,height=0.75\mywidth\linewidth]{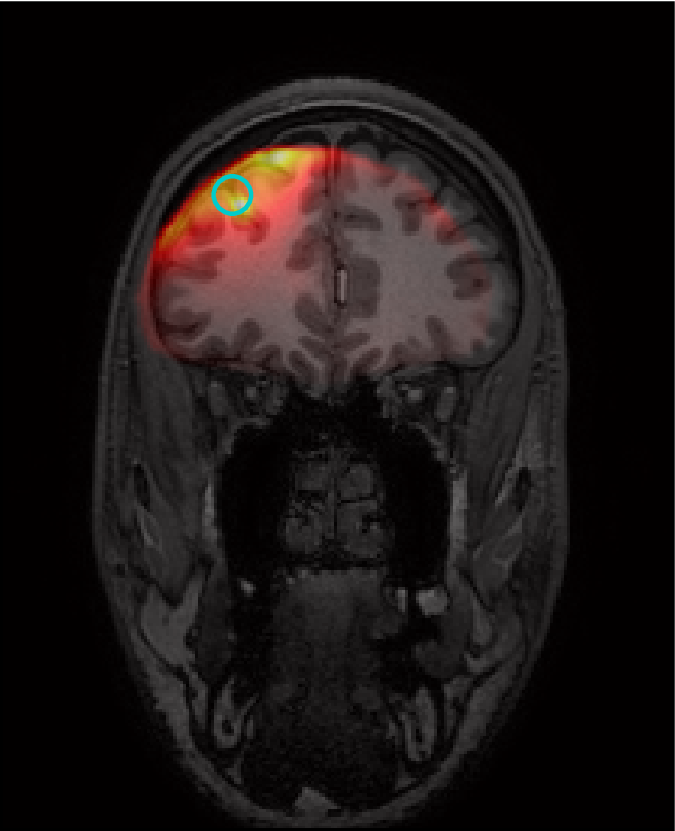}
     \end{center}
      \end{minipage}\vspace{0.5cm}\begin{minipage}{0.18\textwidth}
          \begin{center}
     \includegraphics[trim={10.5cm 7.5cm 8cm 2.8cm},clip,height=0.75\mywidth\linewidth]{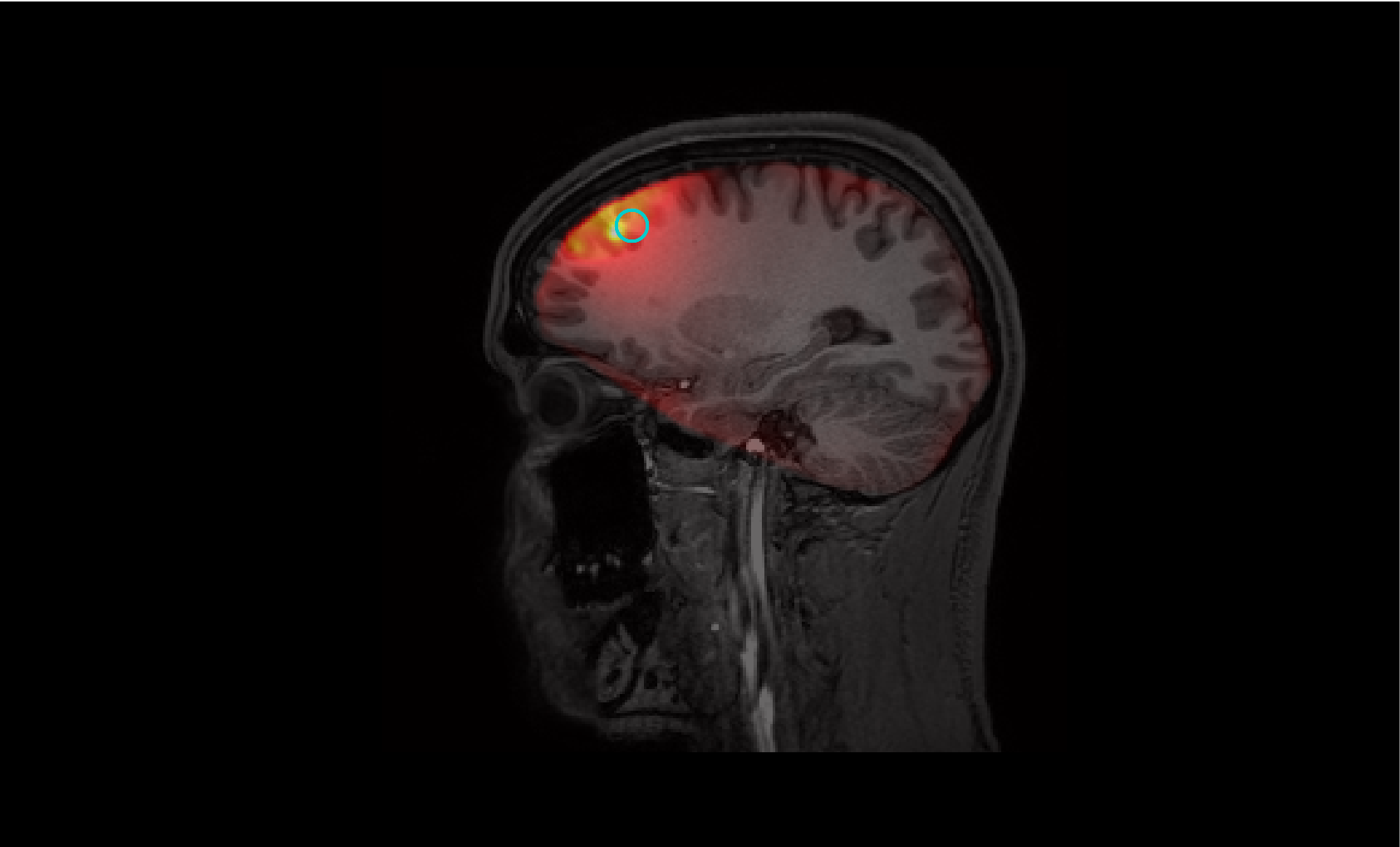}
     \end{center}
      \end{minipage}

      \begin{minipage}{0.02\textwidth}
    \rotatebox{90}{\small{\bf CG-IG (EM)}}
          \end{minipage}\begin{minipage}{\TransversalSz\textwidth}
          \begin{center}
     \includegraphics[trim={2cm 1.6cm 2.1cm 3cm},clip,height=\mywidth\linewidth]{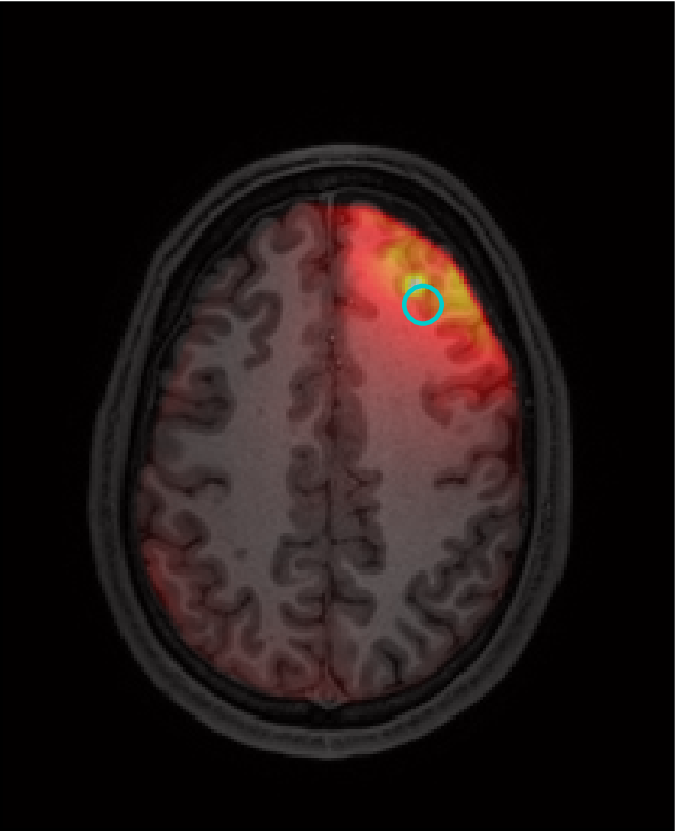}
     \end{center}
      \end{minipage}\begin{minipage}{0.18\textwidth}
          \begin{center}
     \includegraphics[trim={1cm 6cm 1.5cm 1cm},clip,height=0.75\mywidth\linewidth]{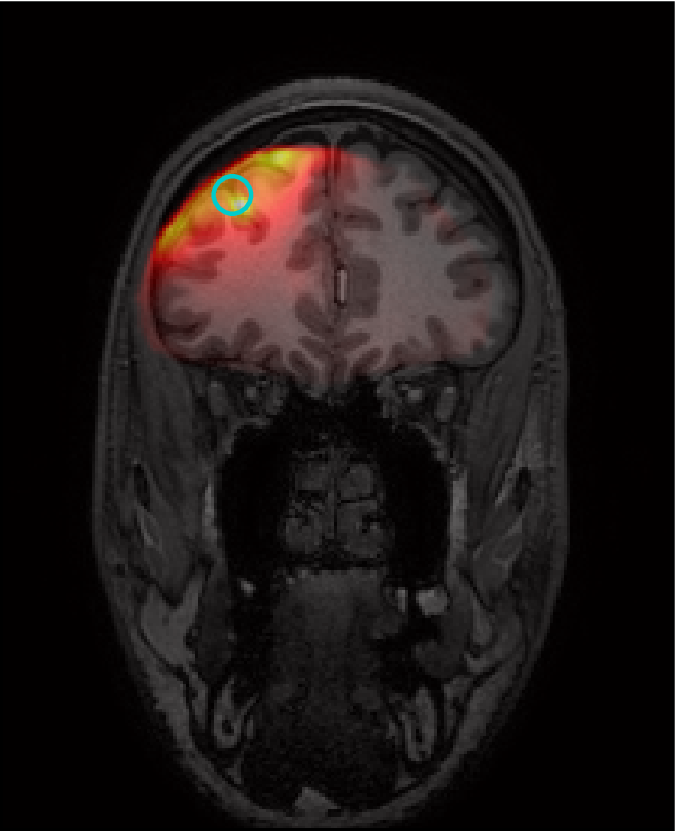}
     \end{center}
      \end{minipage}\vspace{0.5cm}\begin{minipage}{0.18\textwidth}
          \begin{center}
     \includegraphics[trim={10.5cm 7.5cm 8cm 2.8cm},clip,height=0.75\mywidth\linewidth]{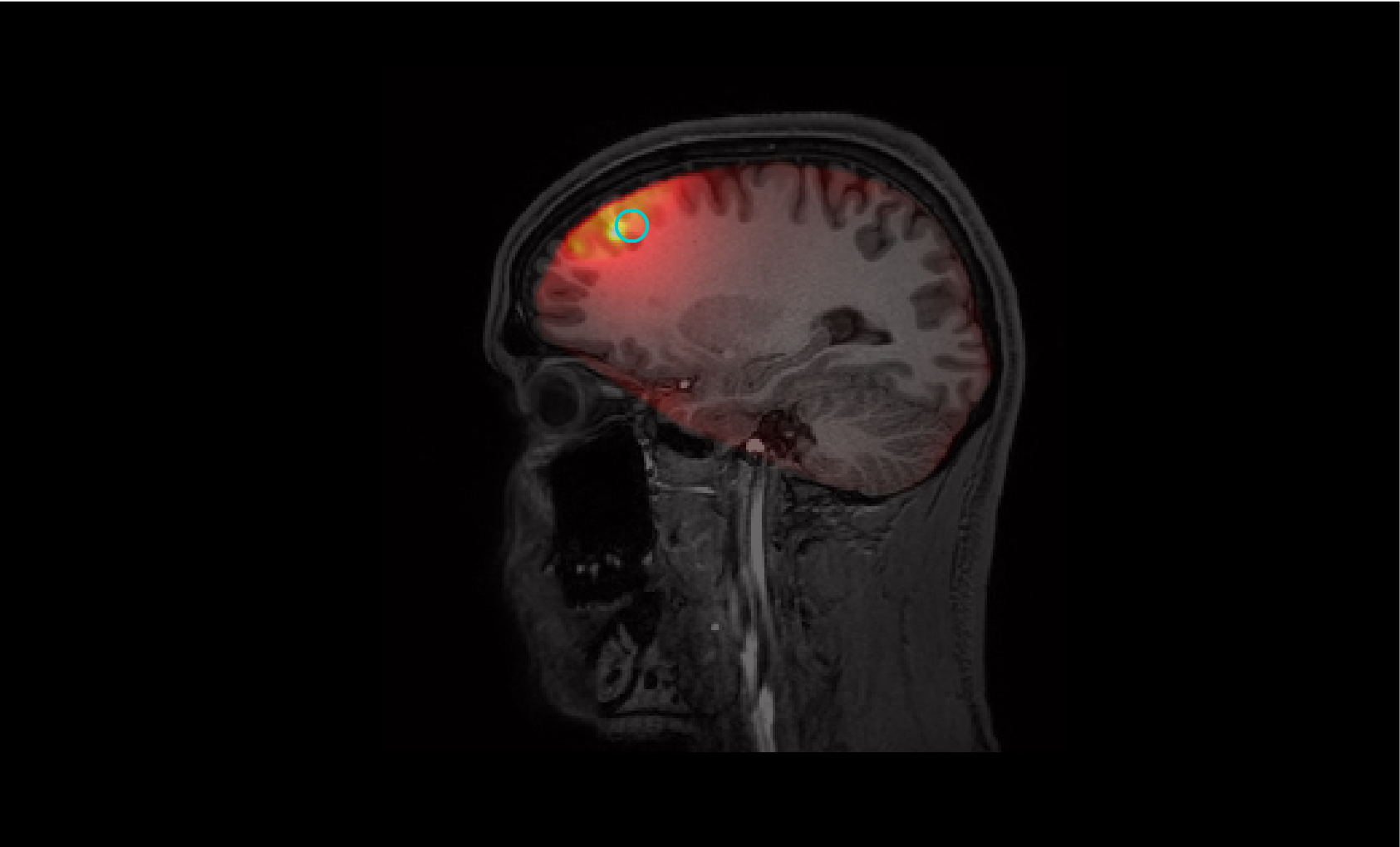}
     \end{center}
      \end{minipage}

      \begin{minipage}{0.02\textwidth}
    \rotatebox{90}{\small{\bf CG-IG (IAS)}}
          \end{minipage}\begin{minipage}{\TransversalSz\textwidth}
          \begin{center}
     \includegraphics[trim={2cm 1.6cm 2.1cm 3cm},clip,height=\mywidth\linewidth]{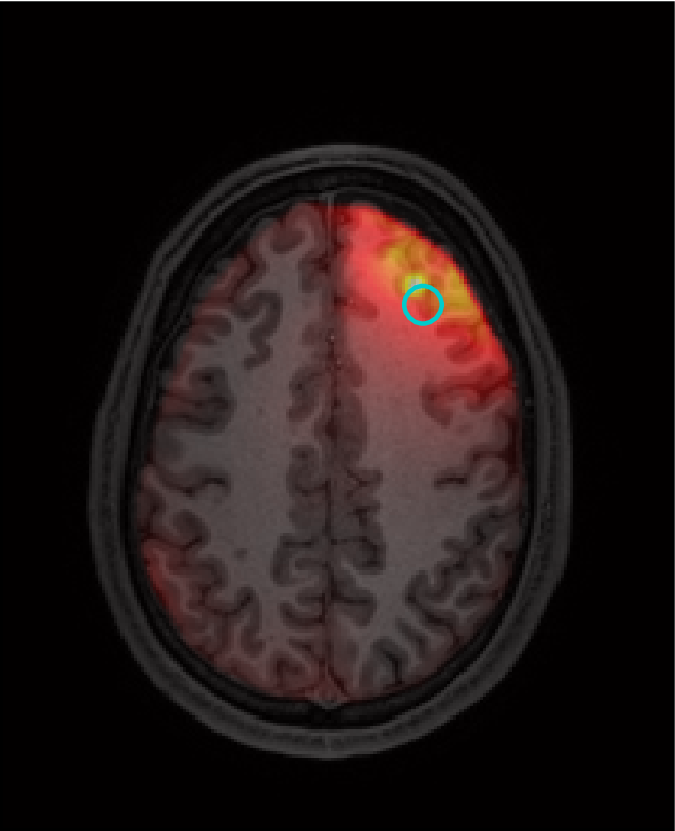}
     \end{center}
      \end{minipage}\begin{minipage}{0.18\textwidth}
          \begin{center}
     \includegraphics[trim={1cm 6cm 1.5cm 1cm},clip,height=0.75\mywidth\linewidth]{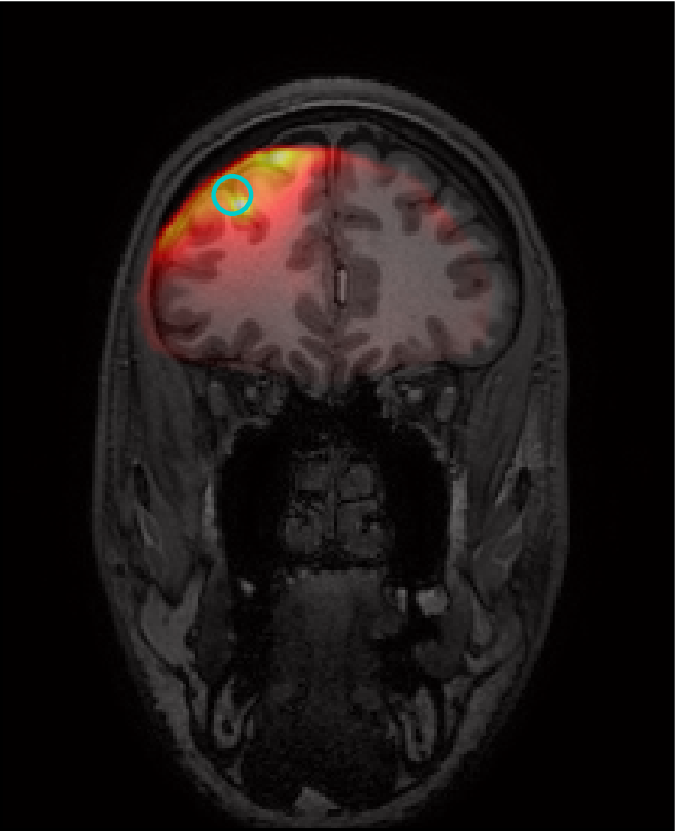}
     \end{center}
      \end{minipage}\vspace{0.5cm}\begin{minipage}{0.18\textwidth}
          \begin{center}
     \includegraphics[trim={10.5cm 7.5cm 8cm 2.8cm},clip,height=0.75\mywidth\linewidth]{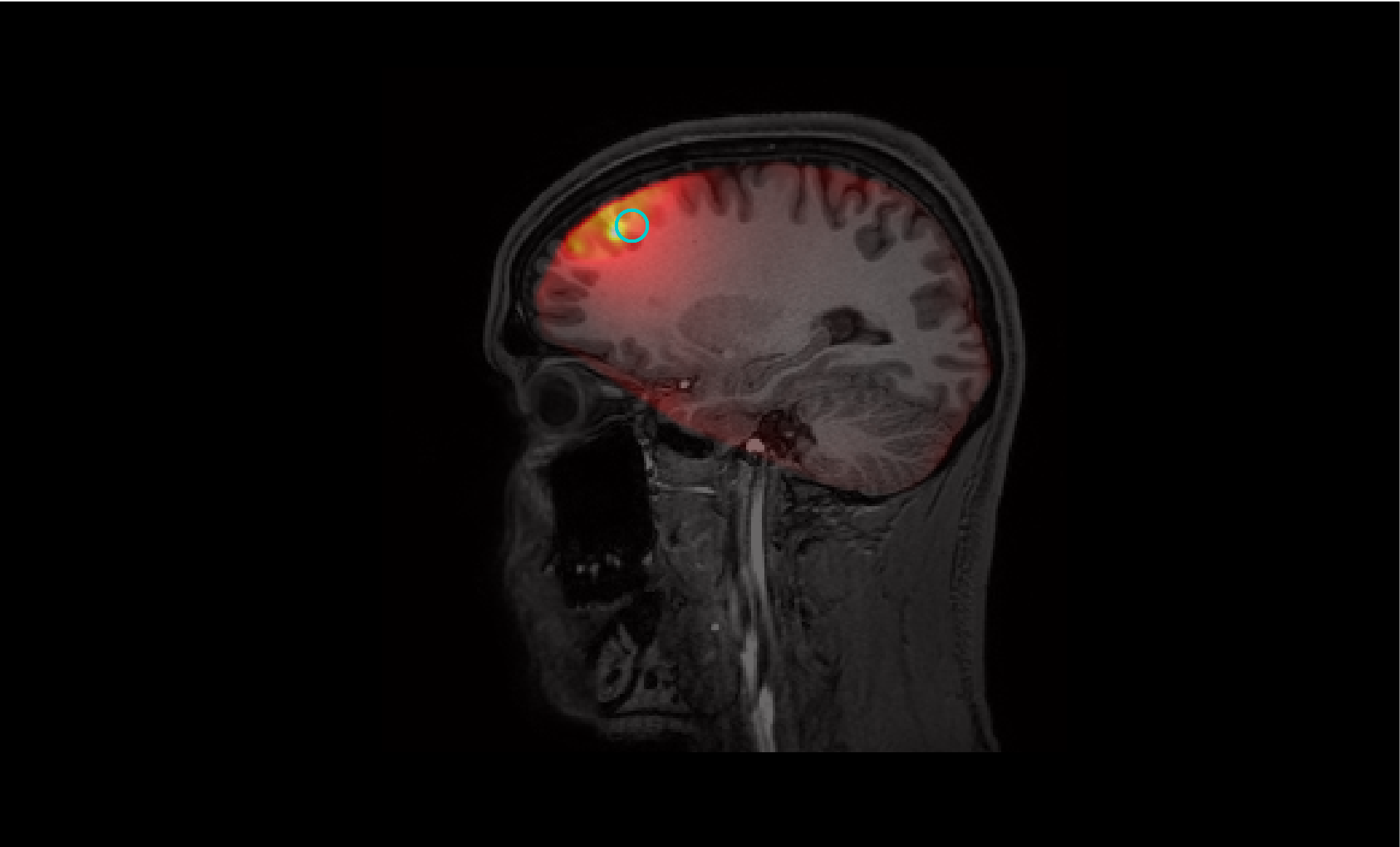}
     \end{center}
      \end{minipage}
    \caption{Estimated distributions of the simulated brain activity at 12 mm depth computed using methods with Gaussian priors. The distributions are presented in three plain cuts of magnetic resonance images. The turquoise ring shows the location of the actual source to be estimated. The colored region, ranging from dark red to yellow, represents the distribution and its local magnitude. Slices have been taken at the location of the maximum estimated magnitude.}
    \label{fig:MRI2Gauss_deep}
\end{figure}

\begin{figure}
\newcommand{\mywidth}{1}
\newcommand{\TransversalSz}{0.16}
    \centering
    \begin{minipage}{0.02\textwidth}
    \rotatebox{90}{\small{\bf wL}}
          \end{minipage}\begin{minipage}{\TransversalSz\textwidth}
          \begin{center}
     \includegraphics[trim={2cm 1.6cm 2.1cm 3cm},clip,height=\mywidth\linewidth]{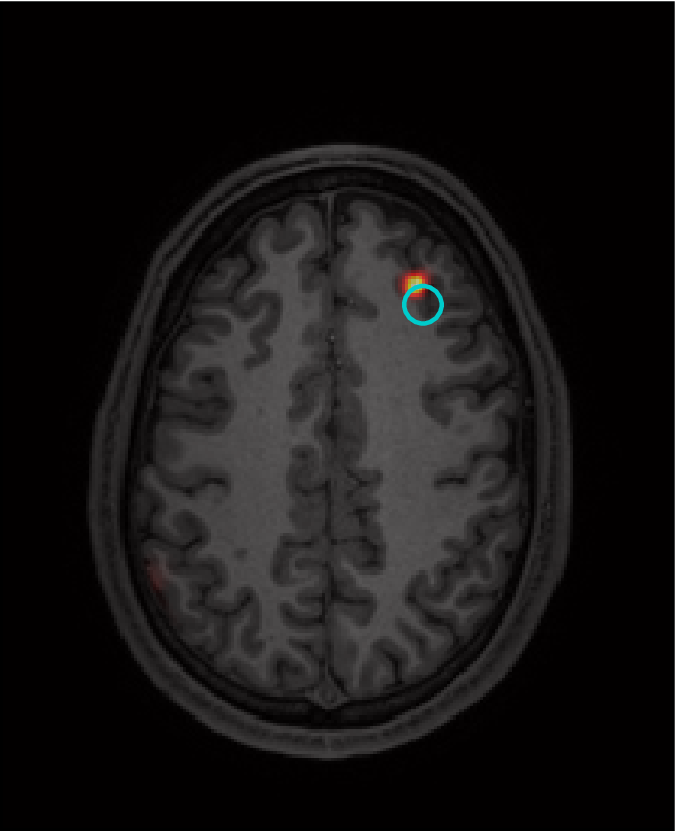}
     \end{center}
      \end{minipage}\begin{minipage}{0.18\textwidth}
          \begin{center}
     \includegraphics[trim={1cm 6cm 1.5cm 1cm},clip,height=0.75\mywidth\linewidth]{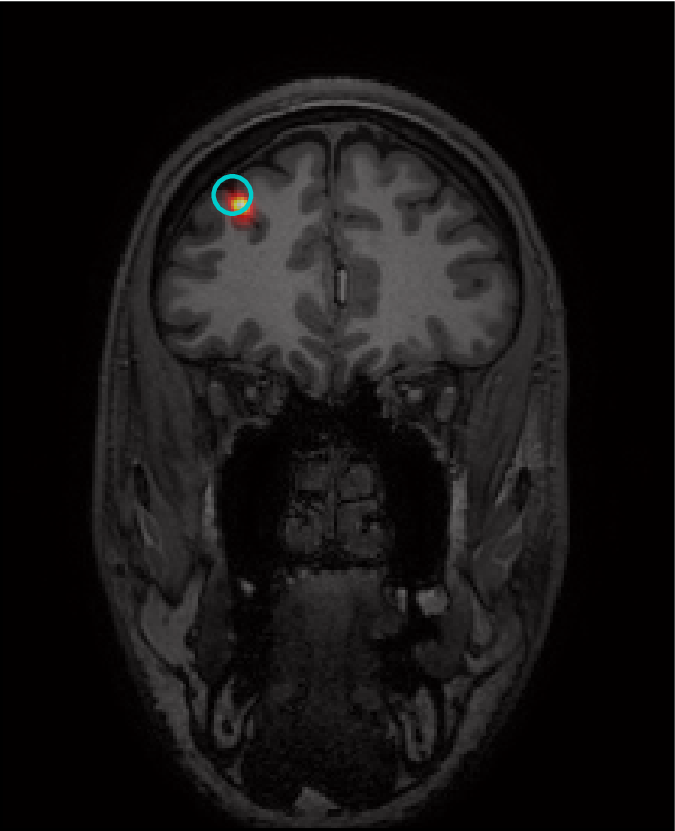}
     \end{center}
      \end{minipage}\vspace{0.5cm}\begin{minipage}{0.18\textwidth}
          \begin{center}
     \includegraphics[trim={10.5cm 7.5cm 8cm 2.8cm},clip,height=0.75\mywidth\linewidth]{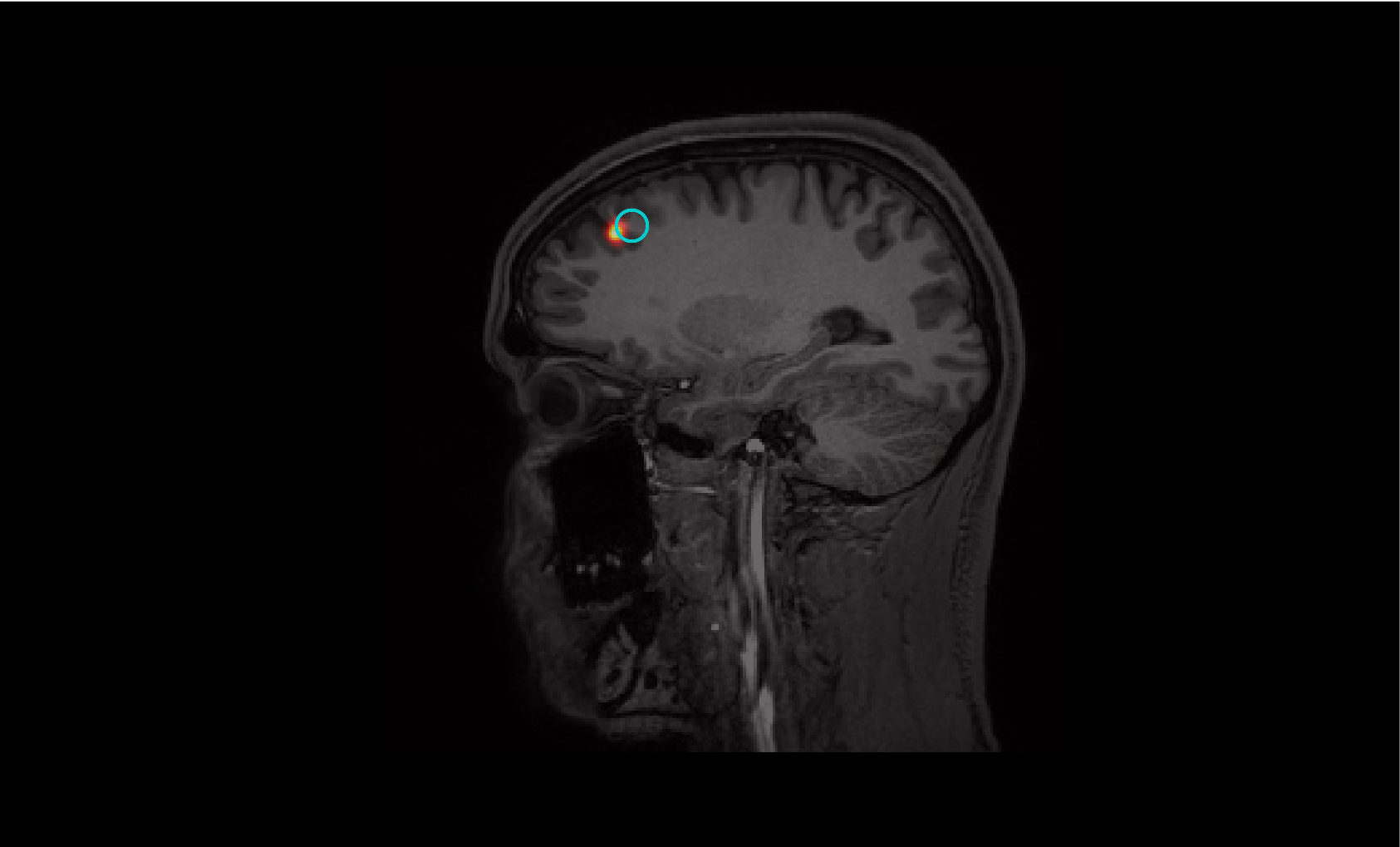}
     \end{center}
      \end{minipage}
      \begin{minipage}{0.02\textwidth}
    \rotatebox{90}{\small{\bf wGL}}
          \end{minipage}\begin{minipage}{\TransversalSz\textwidth}
          \begin{center}
     \includegraphics[trim={2cm 1.6cm 2.1cm 3cm},clip,height=\mywidth\linewidth]{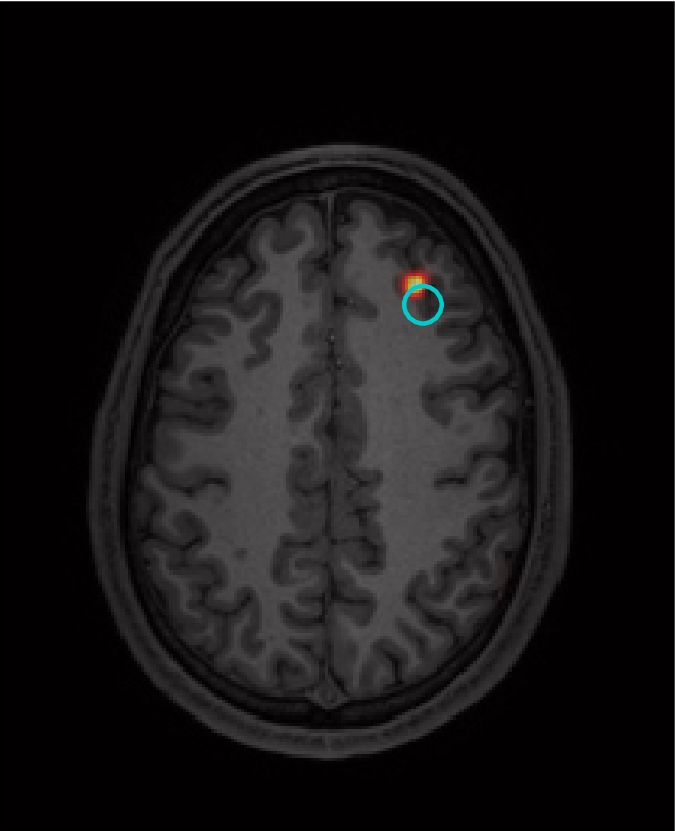}
     \end{center}
      \end{minipage}\begin{minipage}{0.18\textwidth}
          \begin{center}
     \includegraphics[trim={1cm 6cm 1.5cm 1cm},clip,height=0.75\mywidth\linewidth]{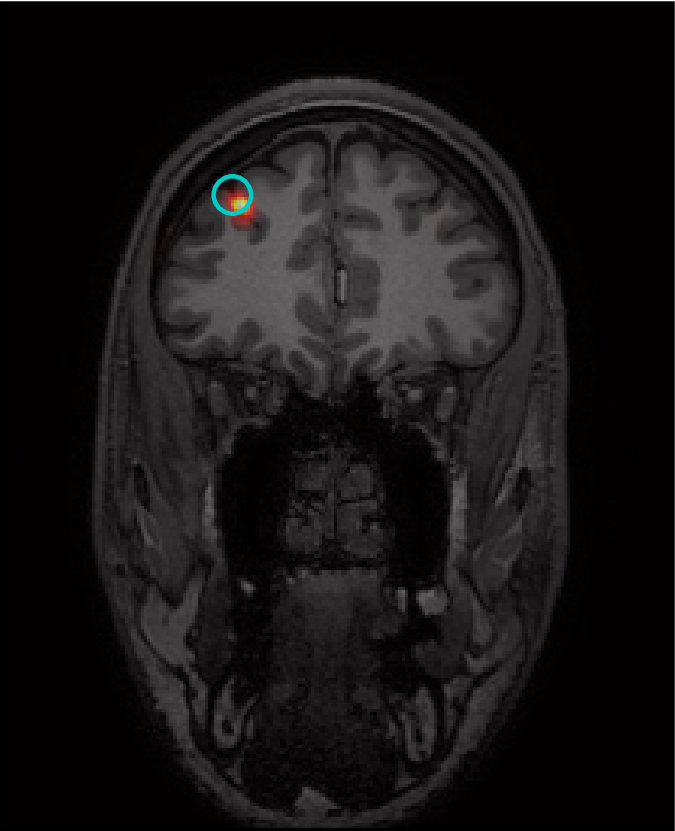}
     \end{center}
      \end{minipage}\vspace{0.5cm}\begin{minipage}{0.18\textwidth}
          \begin{center}
     \includegraphics[trim={10.5cm 7.5cm 8cm 2.8cm},clip,height=0.75\mywidth\linewidth]{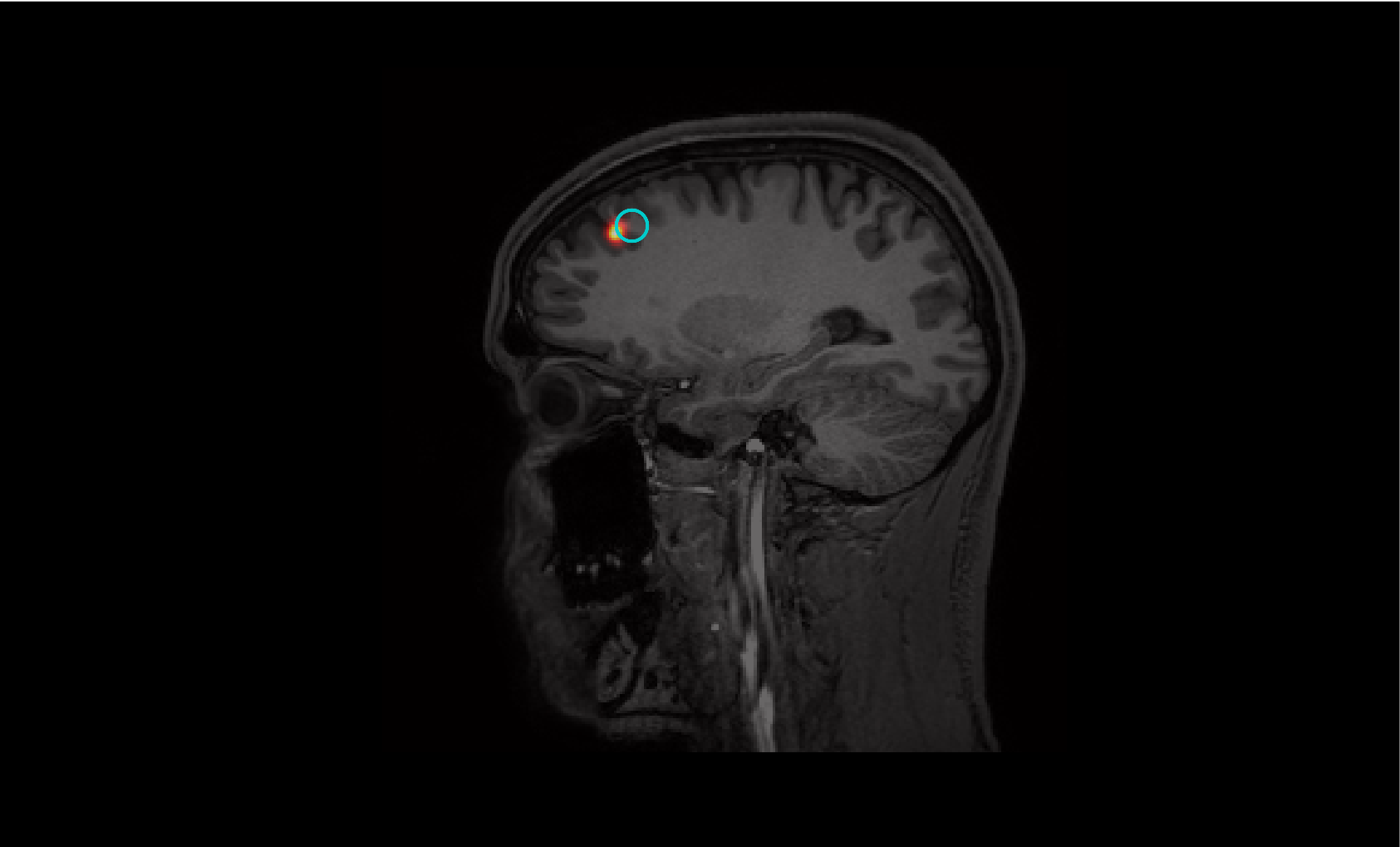}
     \end{center}
      \end{minipage}
      \begin{minipage}{0.02\textwidth}
    \rotatebox{90}{\small{\bf wCL (EM)}}
          \end{minipage}\begin{minipage}{\TransversalSz\textwidth}
          \begin{center}
     \includegraphics[trim={2cm 1.6cm 2.1cm 3cm},clip,height=\mywidth\linewidth]{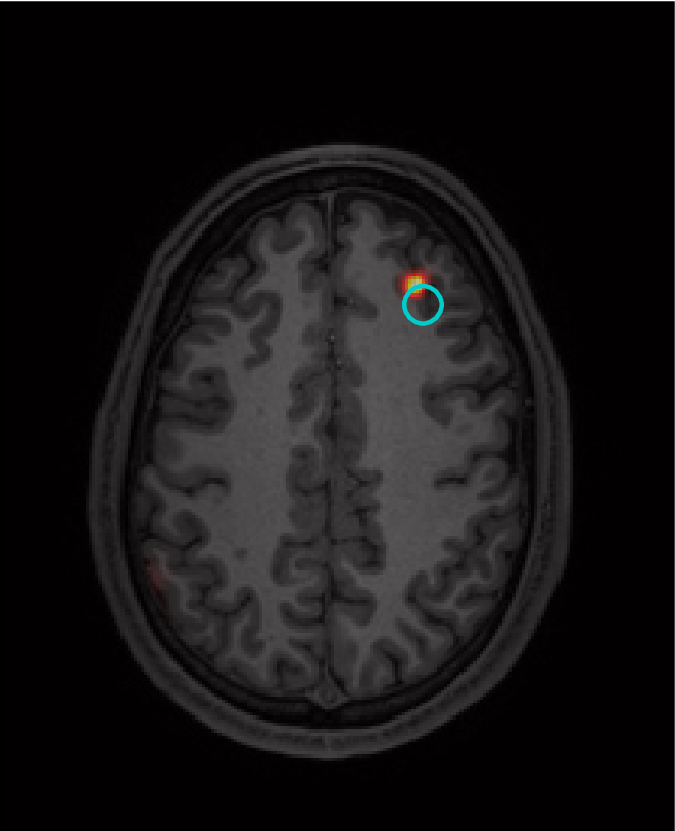}
     \end{center}
      \end{minipage}\begin{minipage}{0.18\textwidth}
          \begin{center}
     \includegraphics[trim={1cm 6cm 1.5cm 1cm},clip,height=0.75\mywidth\linewidth]{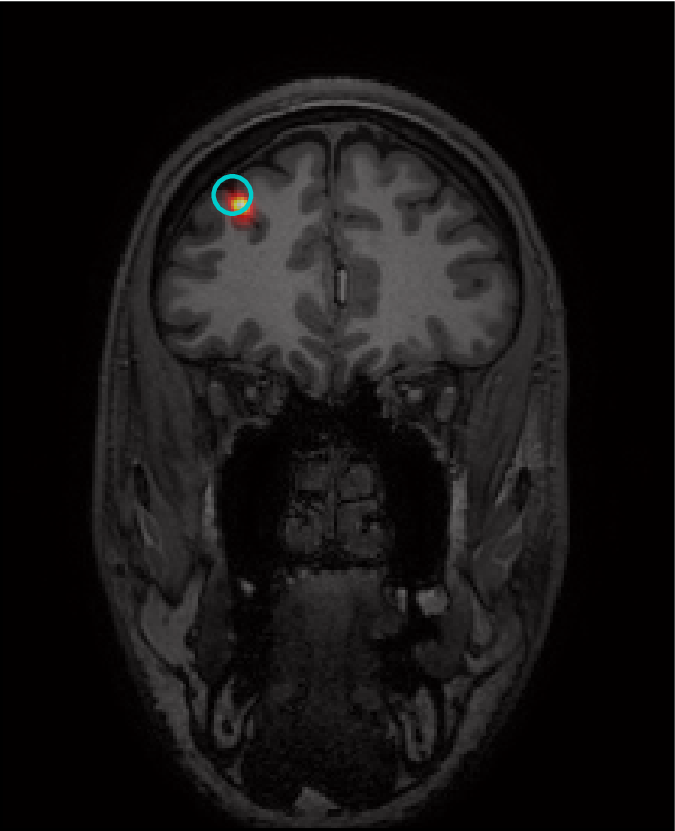}
     \end{center}
      \end{minipage}\vspace{0.5cm}\begin{minipage}{0.18\textwidth}
          \begin{center}
     \includegraphics[trim={10.5cm 7.5cm 8cm 2.8cm},clip,height=0.75\mywidth\linewidth]{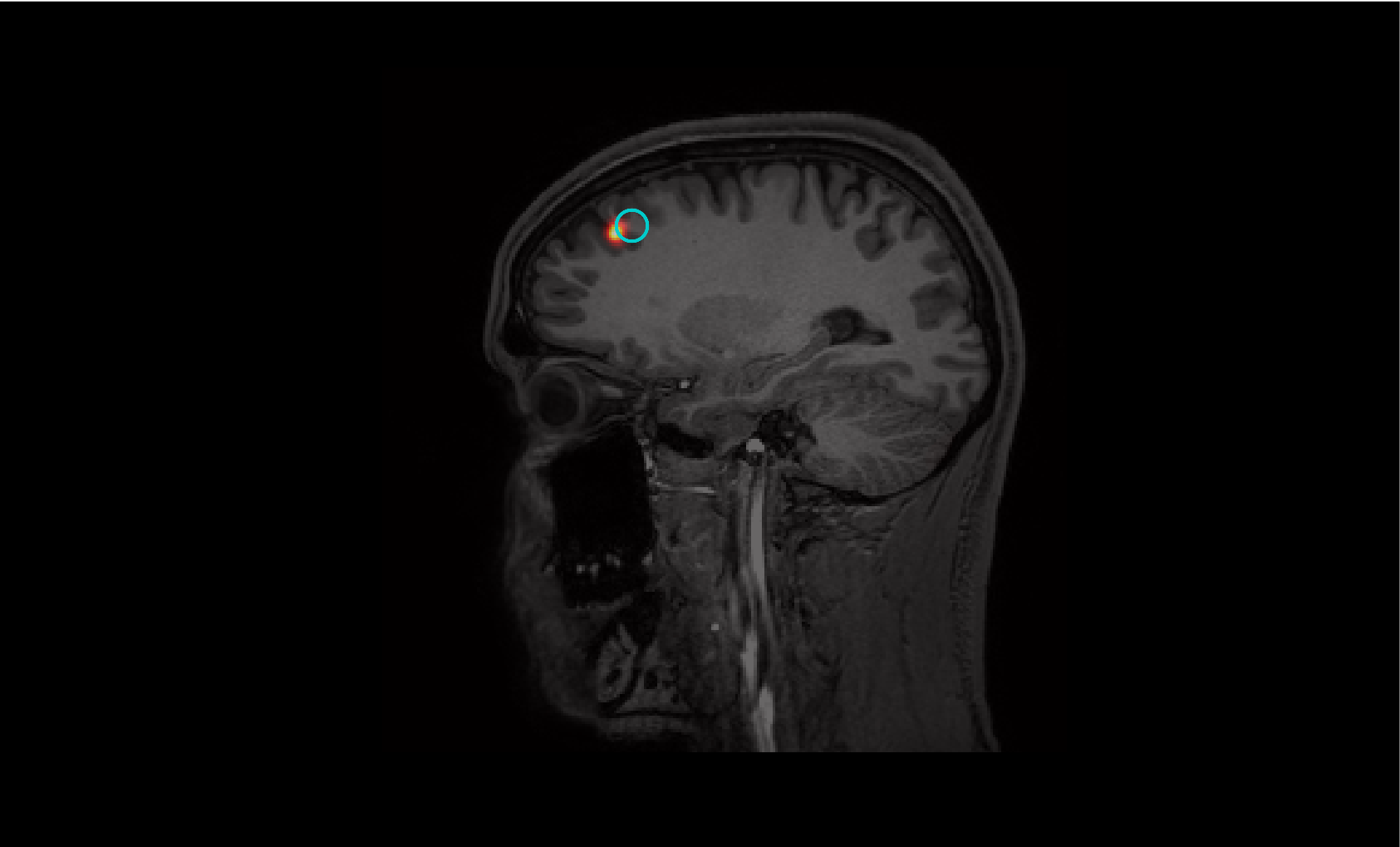}
     \end{center}
      \end{minipage}
      \begin{minipage}{0.02\textwidth}
    \rotatebox{90}{\small{\bf wCL (IAS)}}
          \end{minipage}\begin{minipage}{\TransversalSz\textwidth}
          \begin{center}
     \includegraphics[trim={2cm 1.6cm 2.1cm 3cm},clip,height=\mywidth\linewidth]{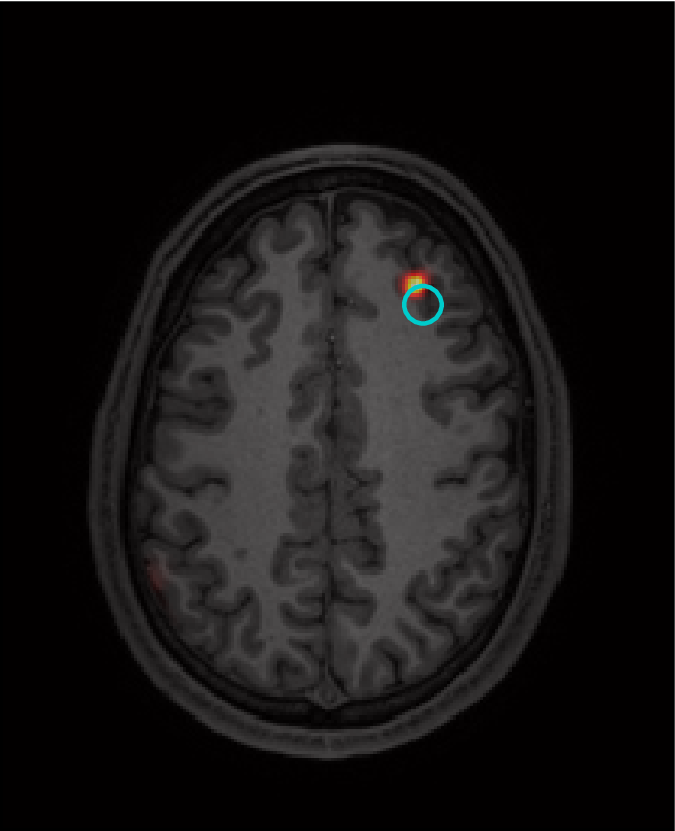}
     \end{center}
      \end{minipage}\begin{minipage}{0.18\textwidth}
          \begin{center}
     \includegraphics[trim={1cm 6cm 1.5cm 1cm},clip,height=0.75\mywidth\linewidth]{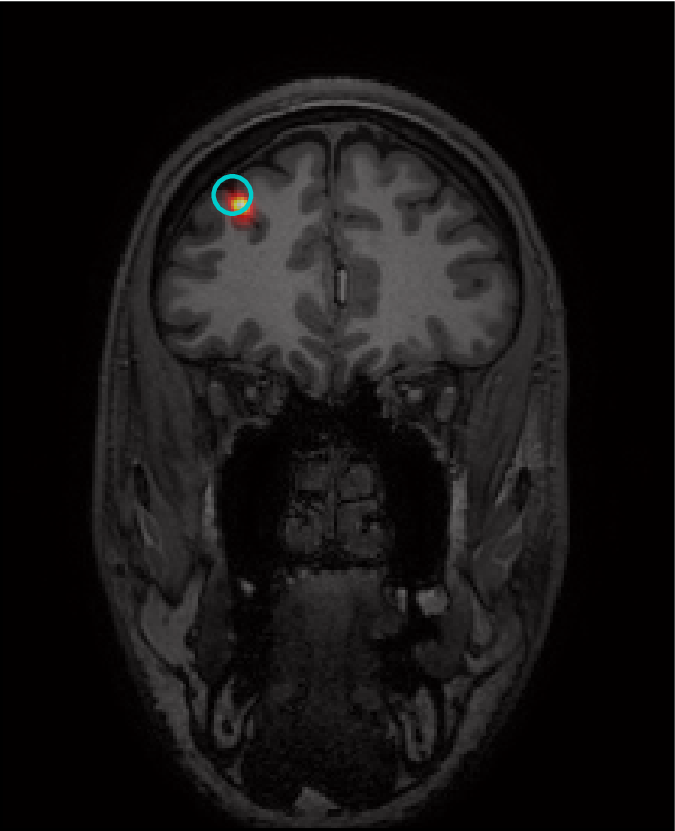}
     \end{center}
      \end{minipage}\vspace{0.5cm}\begin{minipage}{0.18\textwidth}
          \begin{center}
     \includegraphics[trim={10.5cm 7.5cm 8cm 2.8cm},clip,height=0.75\mywidth\linewidth]{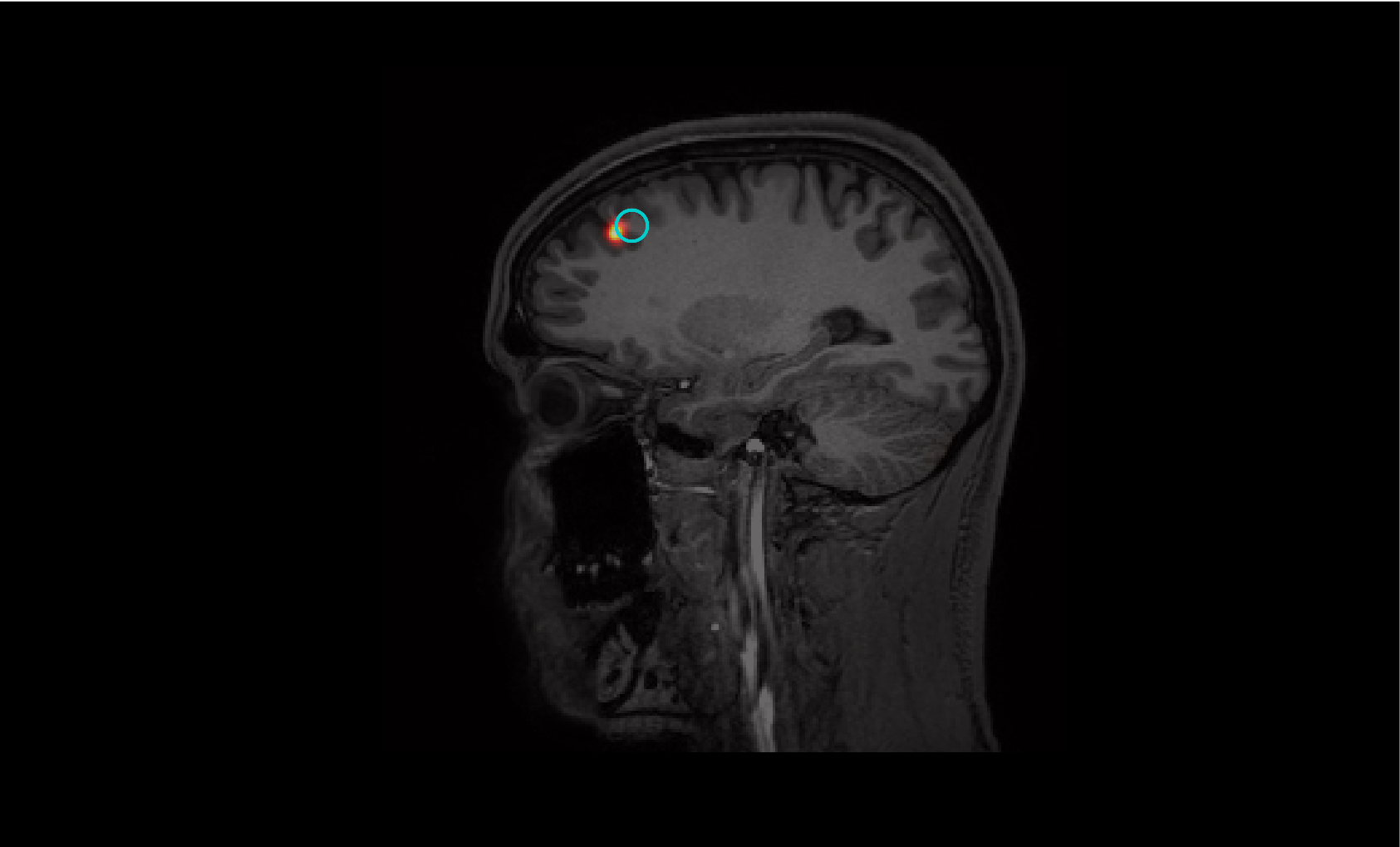}
     \end{center}
      \end{minipage}
      \begin{minipage}{0.02\textwidth}
    \rotatebox{90}{\small{\bf wCGL (EM)}}
          \end{minipage}\begin{minipage}{\TransversalSz\textwidth}
          \begin{center}
     \includegraphics[trim={2cm 1.6cm 2.1cm 3cm},clip,height=\mywidth\linewidth]{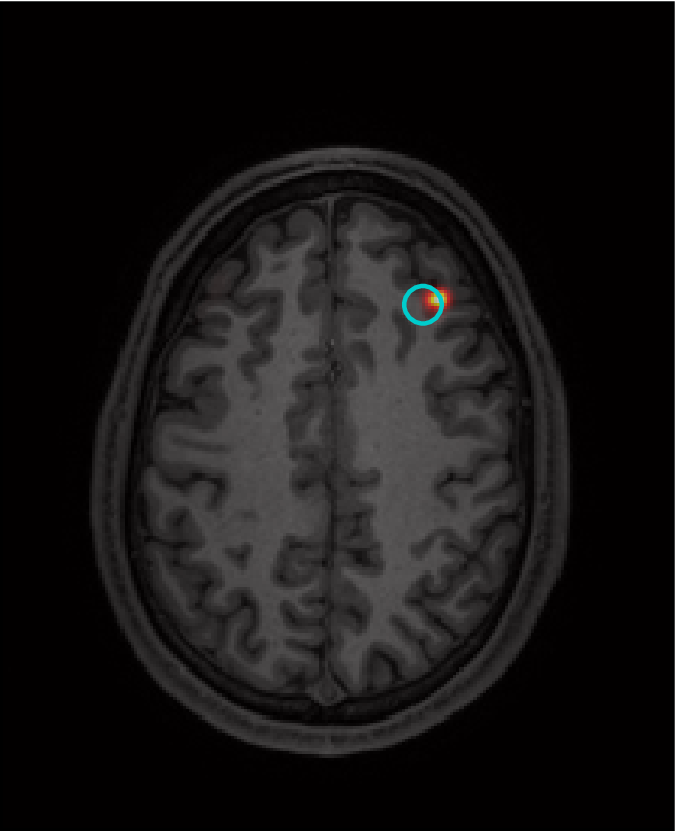}
     \end{center}
      \end{minipage}\begin{minipage}{0.18\textwidth}
          \begin{center}
     \includegraphics[trim={1cm 6cm 1.5cm 1cm},clip,height=0.75\mywidth\linewidth]{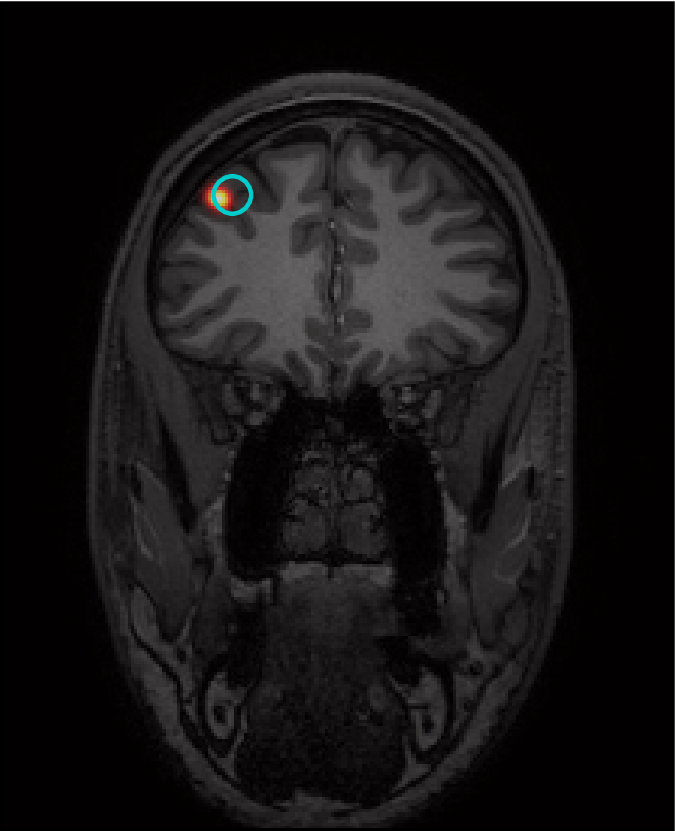}
     \end{center}
      \end{minipage}\vspace{0.5cm}\begin{minipage}{0.18\textwidth}
          \begin{center}
     \includegraphics[trim={10.5cm 7.5cm 8cm 2.8cm},clip,height=0.75\mywidth\linewidth]{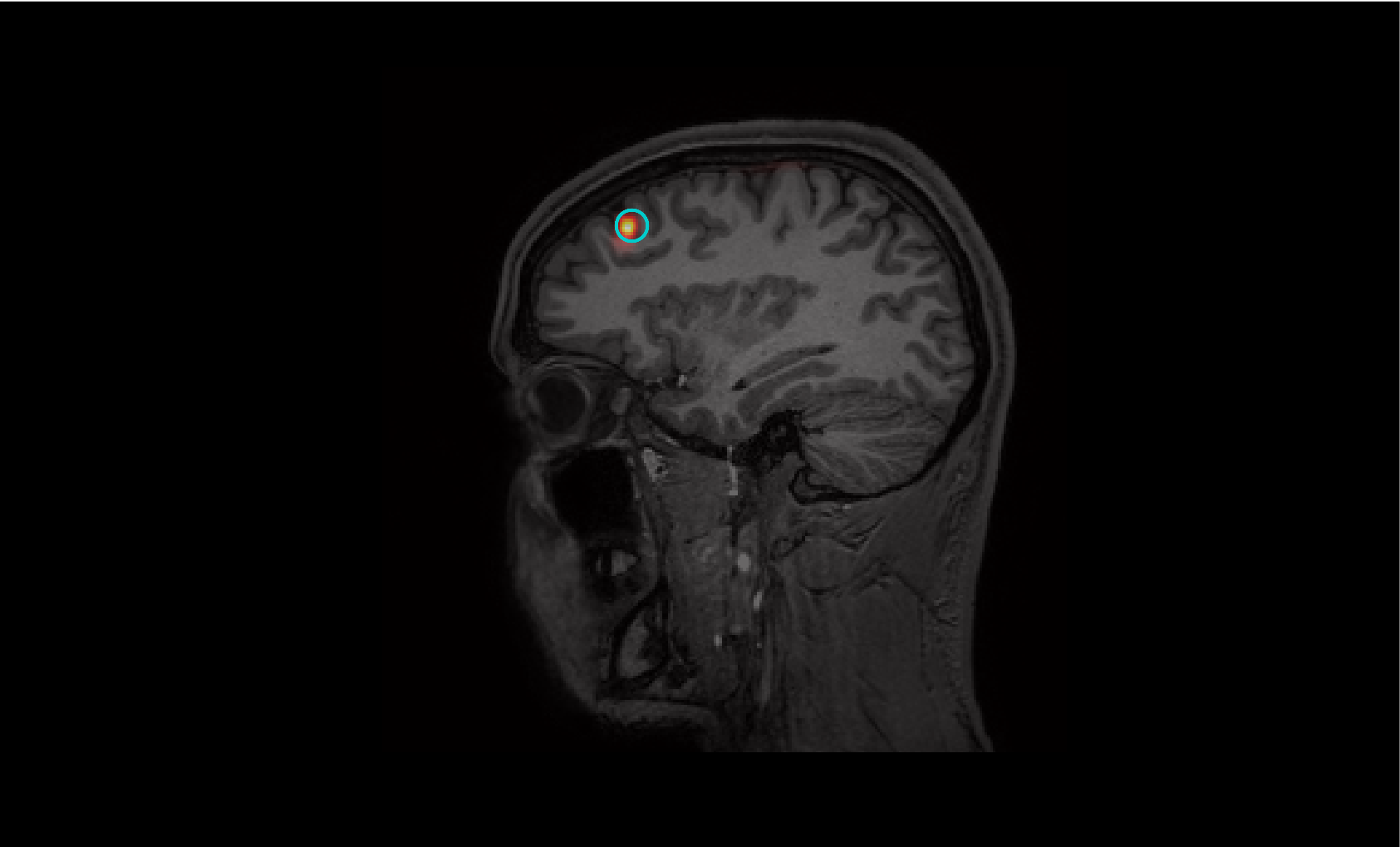}
     \end{center}
      \end{minipage}
      \begin{minipage}{0.02\textwidth}
    \rotatebox{90}{\small{\bf wCGL (IAS)}}
          \end{minipage}\begin{minipage}{\TransversalSz\textwidth}
          \begin{center}
     \includegraphics[trim={2cm 1.6cm 2.1cm 3cm},clip,height=\mywidth\linewidth]{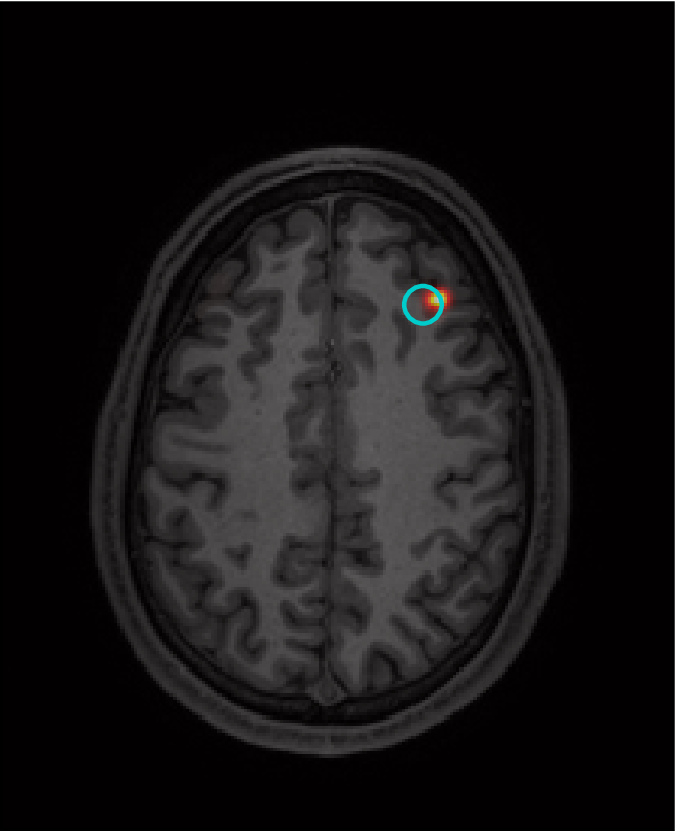}
     \end{center}
      \end{minipage}\begin{minipage}{0.18\textwidth}
          \begin{center}
     \includegraphics[trim={1cm 6cm 1.5cm 1cm},clip,height=0.75\mywidth\linewidth]{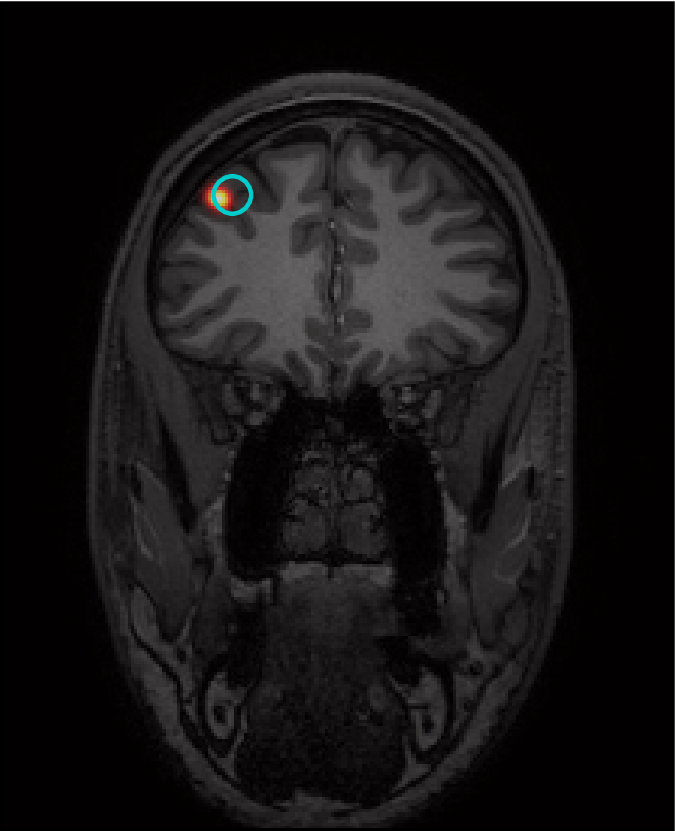}
     \end{center}
      \end{minipage}\vspace{0.5cm}\begin{minipage}{0.18\textwidth}
          \begin{center}
     \includegraphics[trim={10.5cm 7.5cm 8cm 2.8cm},clip,height=0.75\mywidth\linewidth]{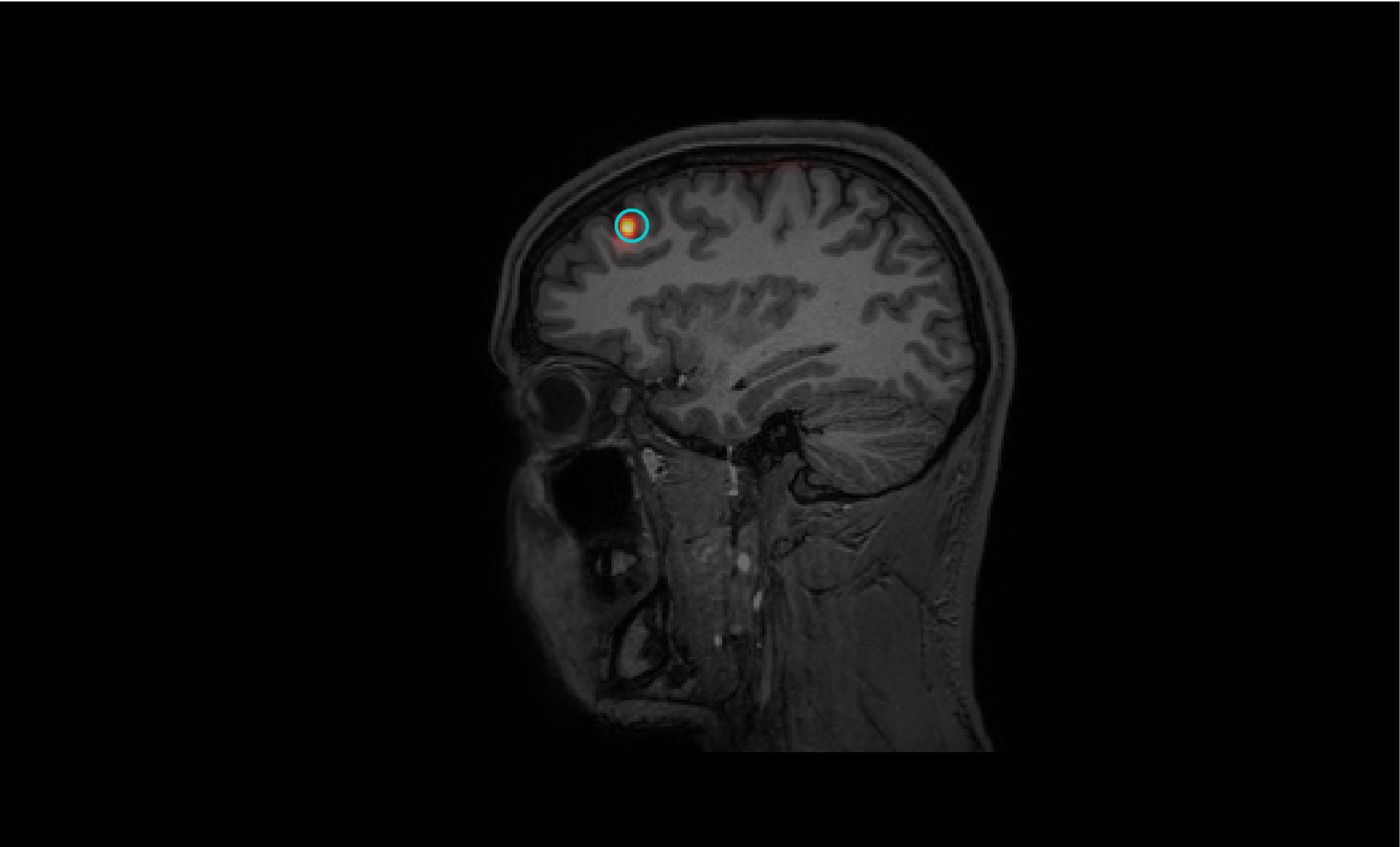}
     \end{center}
      \end{minipage}
    \caption{Estimated distributions of the simulated brain activity at 12 mm depth computed using methods with Laplace priors. The distributions are presented in three plain cuts of magnetic resonance images. The turquoise ring shows the location of the actual source to be estimated. The colored region, ranging from dark red to yellow, represents the distribution and its local magnitude. Slices have been taken at the location of the maximum estimated magnitude.}
    \label{fig:MRI2Lap_deep}
\end{figure}

In this subsection, we visualize the reconstruction results on the MRI images. The algorithms are grouped according to their ability to recover either spread or focal source activity. In particular, Figures \ref{fig:MRI2Gauss_superf} and \ref{fig:MRI2Lap_superf} illustrate the corresponding reconstruction results for methods using Gaussian priors and Laplace or group Laplace priors, respectively, with data from a cortical source. In another setup, we picked a source at 12 mm depth with results shown for spread Gaussian methods in Figure \ref{fig:MRI2Gauss_deep} and focal Laplace methods in Figure \ref{fig:MRI2Lap_deep}. The measurements are contaminated by 5 \% of additive Gaussian noise. Moreover, we computed the EMDs of each reconstruction and presented those in Table \ref{tab:MRIEMD}.

By observing the estimations interpolated to the MRI slices when conditionally Gaussian models are used in Figures \ref{fig:MRI2Gauss_superf} and \ref{fig:MRI2Gauss_deep}, we see very similar and widely spread estimates for wMNE and the weighted Conditionally Gaussian model with Inverse-Gamma hyperprior, and in the case of the Gamma prior when the IAS algorithm is used. Maximum estimate values in yellow are too frontal and lateral with wMNE and two Conditionally Gaussian models when estimating the superficial source. However, the spread of the mentioned high-value region does reach the true source location. With the deeper source, the estimation spreads are wider, and the estimated maxima are slightly too frontal and significantly too lateral. There is no visible difference in estimate distribution among most of the compared Gaussian methods. However, when the EM algorithm is used with the Conditional Gaussian model and the gamma hyperprior, the estimation becomes highly focal.

The estimates from focal-by-design methods, presented in Figures \ref{fig:MRI2Lap_superf} and \ref{fig:MRI2Lap_deep}, are highly similar and estimate both sources close to the true locations. Methods using Group Laplace estimated the superficial source more frontal than methods with Laplace priors. The deeper source is estimated to be nearly at the correct depth by all these methods. Although Laplace methods can, Weighted Group Laplace (wGL) estimates the source slightly too frontal, while Weighted Conditional Group Laplace (wCGL) estimates the source a bit too laterally.

The best estimation among all of the compared methods is obtained with the Weighted Conditional Laplace with EM algorithm for both sources, based on the EMD results provided in Table \ref{tab:MRIEMD}. Comparing hyperparameter updating algorithms, we find that EM yields smaller EMDs than IAS across all methods in both source configurations. Overall, focal methods, including CG-Ga-EM, provide more accurate estimates of the deeper source than methods with wider spreads, as they tend to yield more superficial estimates.

\clearpage

\begin{table}[t]
    \centering
    \small{
    \begin{tabular}{|p{30pt}|p{34pt}|p{32pt}|p{30pt}|p{30pt}|p{30pt}|p{30pt}|p{30pt}|p{30pt}|p{30pt}|p{30pt}|p{30pt}|}
    \hline
    \multicolumn{12}{|c|}{Superficial source}\\ \hline
    Method &  CG-Ga EM & CG-Ga IAS & CG-IG EM & CG-IG IAS & wCGL EM & wCGL IAS & wCL EM & wCL IAS & wGL & wL & wMNE\\ \hline    
    EMD & 47.7 & 56.1 & 50.1 & 54.0 &  45.7 & 46.3 & 44.3 & 44.8 & 48.2 & 44.7 & 56.5\\ \hline
    \multicolumn{12}{|c|}{ Source at 12 mm depth}\\ \hline
    EMD & 58.1 & 69.0 & 64.2 & 67.8 & 52.3 & 52.9 & 53.3 & 53.7  & 55.4 & 55.4 & 69.2 \\ \hline
    \end{tabular}
    }
    \caption{Earth Mover's Distances in millimeters for the compared method to estimate the superficial source and a source at 12 mm depth.}
    \label{tab:MRIEMD}
\end{table}

\newpage 
\subsection{Numerical Evaluation of Focality, Depth Bias, and Noise Robustness}
\subsubsection{Overall algorithmic performance in the presence of noise}

To better understand the overall performance of the Bayesian algorithms and their robustness to noise, we evaluate the reconstructions and the corresponding EMD values using observations generated from radially oriented sources (with respect to the cortical surface). The simulated sources are uniformly distributed across different depths, with approximately 28–30 sources per depth level. The depth is defined with respect to the inner skull surface (i.e., the closest distance from the source to the inner skull boundary). The EMD quantifies the spatial discrepancy between the true and reconstructed source distributions; lower values indicate that the estimated activity is closer to the true source location (i.e., more accurate and focal reconstruction), whereas higher values reflect increased spatial spread, mislocalization, or depth bias.

\begin{figure}[t]
      \centering
      \includegraphics[width=1\linewidth]{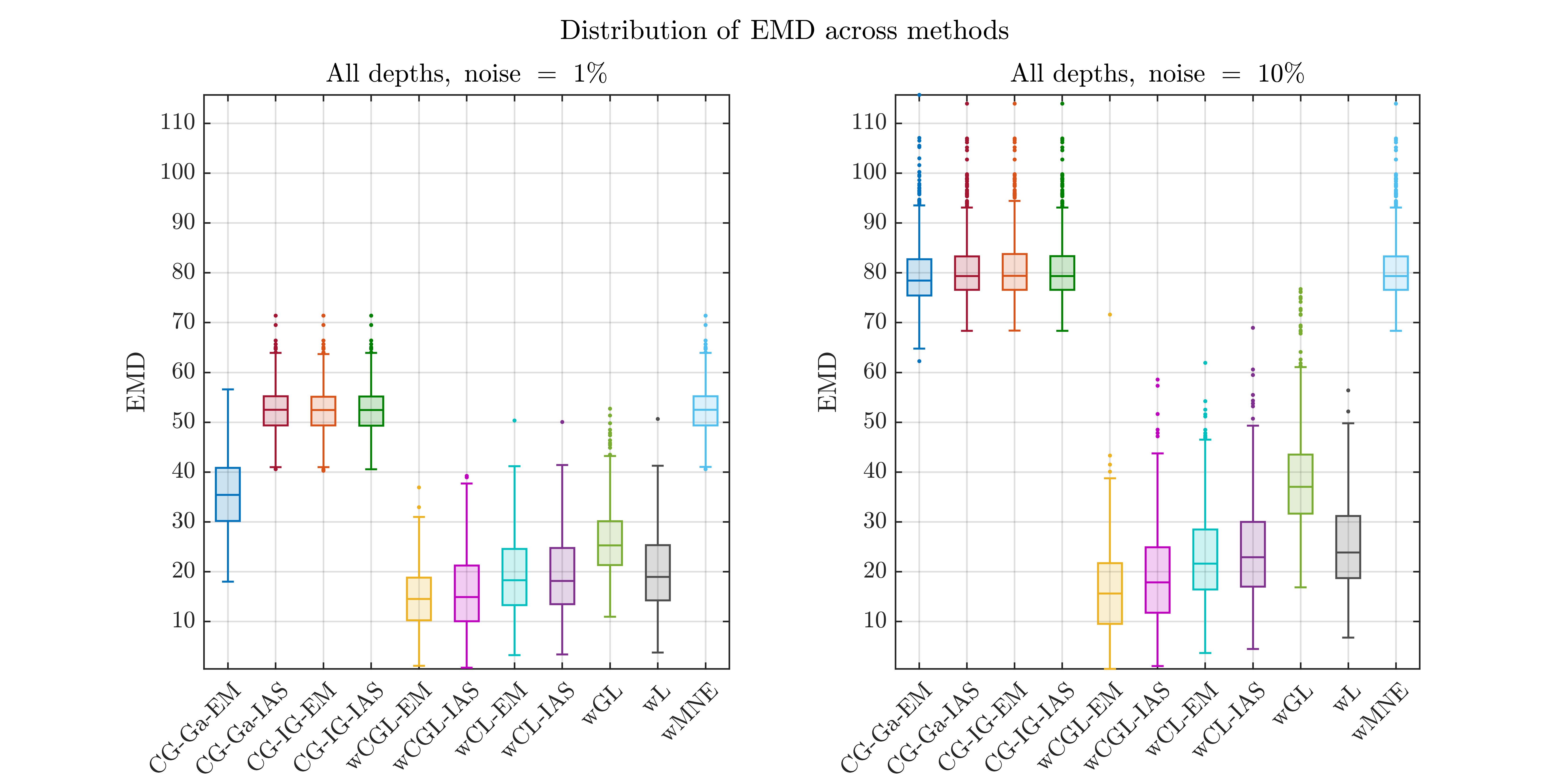}
      \caption{Distribution of EMD values across reconstruction methods for all source depths at two noise levels (1\% and 10\%). Lower EMD indicates better localization accuracy and }\label{fig:EMD_all_depths}
    \end{figure}

Figure~\ref{fig:EMD_all_depths} shows the distribution of EMD values across all tested methods for two noise levels (1\% and 10\%), evaluated over all source depths. At the lower noise level (1\%), the weighted focal models (wCGL and wCL) exhibit the lowest median EMD values, particularly when hyperparameters are estimated using the EM algorithm, indicating more accurate and spatially precise reconstructions. In contrast, the classical conditionally Gaussian approaches (CG variants) and wMNE yield significantly higher EMD values, reflecting more spread and less accurate localization. However, a notable difference can be obtained with CG-Ga-EM. which produces a lower EMD, especially with 1\% of noise. 

When the noise level increases to 10\%, all methods show performance degradation, as evidenced by higher EMD values and greater variability. Nevertheless, the wCGL and wCL with EM updates remain the most robust, maintaining relatively low median EMD and tighter distributions compared to the other approaches, as we have shown in our previous study \citep{Lahtinen2022}. The IAS-based variants generally exhibit higher variability and slightly worse median performance than their EM counterparts. Among the non-Bayesian approaches, wMNE continues to produce relatively high EMD values, indicating limited robustness to noise. 
    Overall, the present results demonstrate that incorporating sensitivity weighting and EM-based hyperparameter estimation leads to improved focality, reduced depth bias, and greater robustness to noise.

\subsubsection{Effect of Source Depth on Reconstruction Performance (Average EMD per depth)}
In this subsection, we investigate how the performance of different algorithms is affected by the depth of the simulated source. To enhance clarity and reduce visual complexity in the subsequent figures, we limit the number of compared methods. As the EM-based approaches consistently show slightly better overall performance, the IAS-based variants are excluded from the  current analysis.

\begin{figure}[t]
    \centering
           \centering
        \includegraphics[width=1\textwidth]{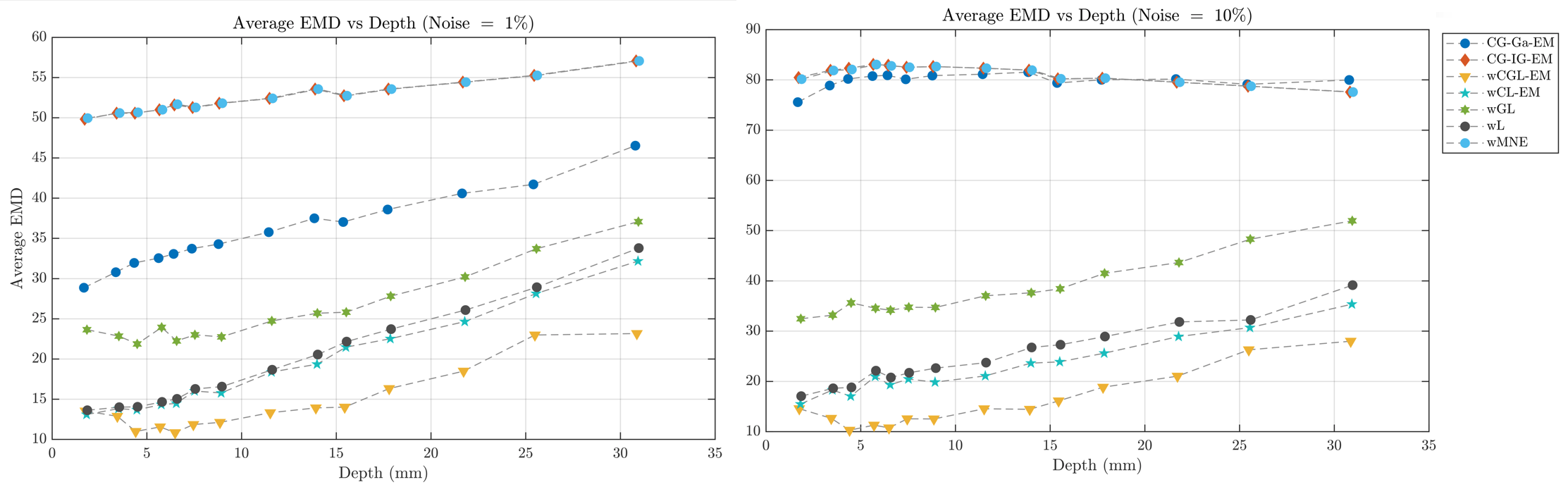}
       \caption{Average EMD as a function of source depth for different reconstruction algorithms. Lower EMD values indicate more accurate source localization. The results show that localization accuracy degrades with increasing depth, highlighting depth bias across methods, while higher noise further amplifies this effect.}
    \label{fig:EMD_vs_depth}
\end{figure}

Figures~\ref{fig:EMD_vs_depth} shows the average EMD as a function of source depth for the considered algorithms under low (1\%) and moderate (10\%) noise levels. For all methods, the EMD generally increases with depth (especially after 15mm depth), indicating a degradation in localization accuracy for deeper sources and highlighting the presence of depth bias.

Among the compared approaches, the weighted conditionally Laplace models (wCGL-EM and wCL-EM) consistently achieve the lowest EMD values across depths, demonstrating better robustness to depth-related effects. In contrast, classical conditionally Gaussian methods (CG variants) and wMNE exhibit significantly higher EMD values, with relatively weak sensitivity to depth but overall poorer localization accuracy.

Increasing the noise level to 10\% amplifies these effects, leading to higher EMD values across all depths and methods. Nevertheless, the relative performance ranking remains similar, with weighted models maintaining superior performance, particularly for deeper sources.

\subsubsection{Statistical analysis for sources in two different depths}

Based on the previous results (Figure~\ref{fig:EMD_vs_depth}), we observed that for shallow sources (0–15 mm), the algorithms exhibit similar performance on average, whereas for deeper sources (greater than 15 mm), a clear degradation in performance is evident. To better understand this depth-dependent behavior, we further analyze the estimated EMD values at two representative depth ranges.

\begin{figure}[t]
    \centering
    \includegraphics[width=\linewidth]{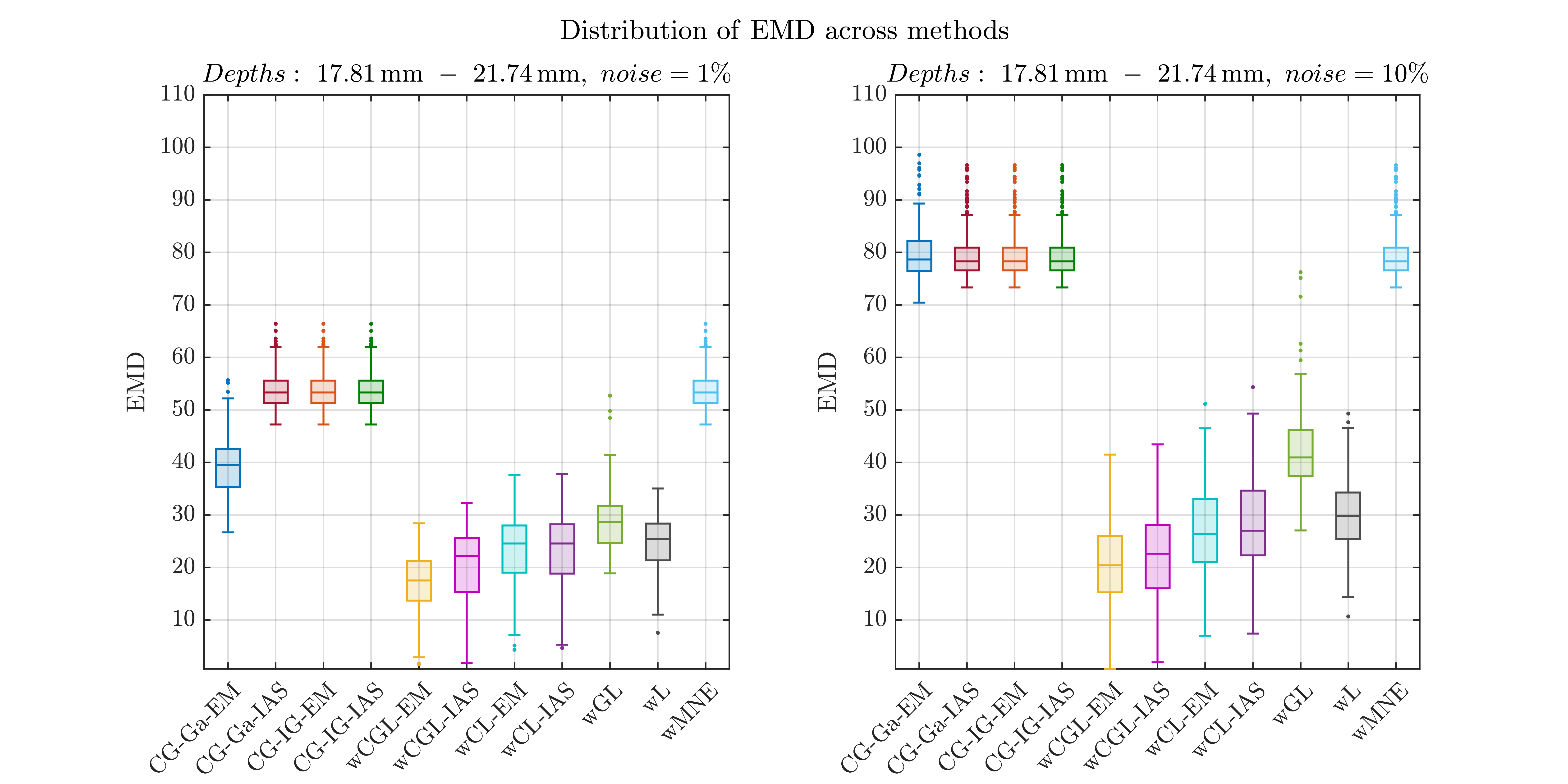}
    \caption{Deep sources (17.81--21.74 mm). Each panel shows results under 1\% (left) and 10\% (right) noise levels. Performance degrades significantly for deeper sources, especially under higher noise.}
    \label{fig:emd_deep}
\end{figure}

In the following, figure~\ref{fig:emd_deep} illustrates the distribution of EMD across methods for deep sources (17.81--21.74 mm) under low (1\%) and high (10\%) noise levels, showing higher error values and slightly increased dispersion across most methods.

\begin{figure}[t]
    \centering
\includegraphics[width=\linewidth]{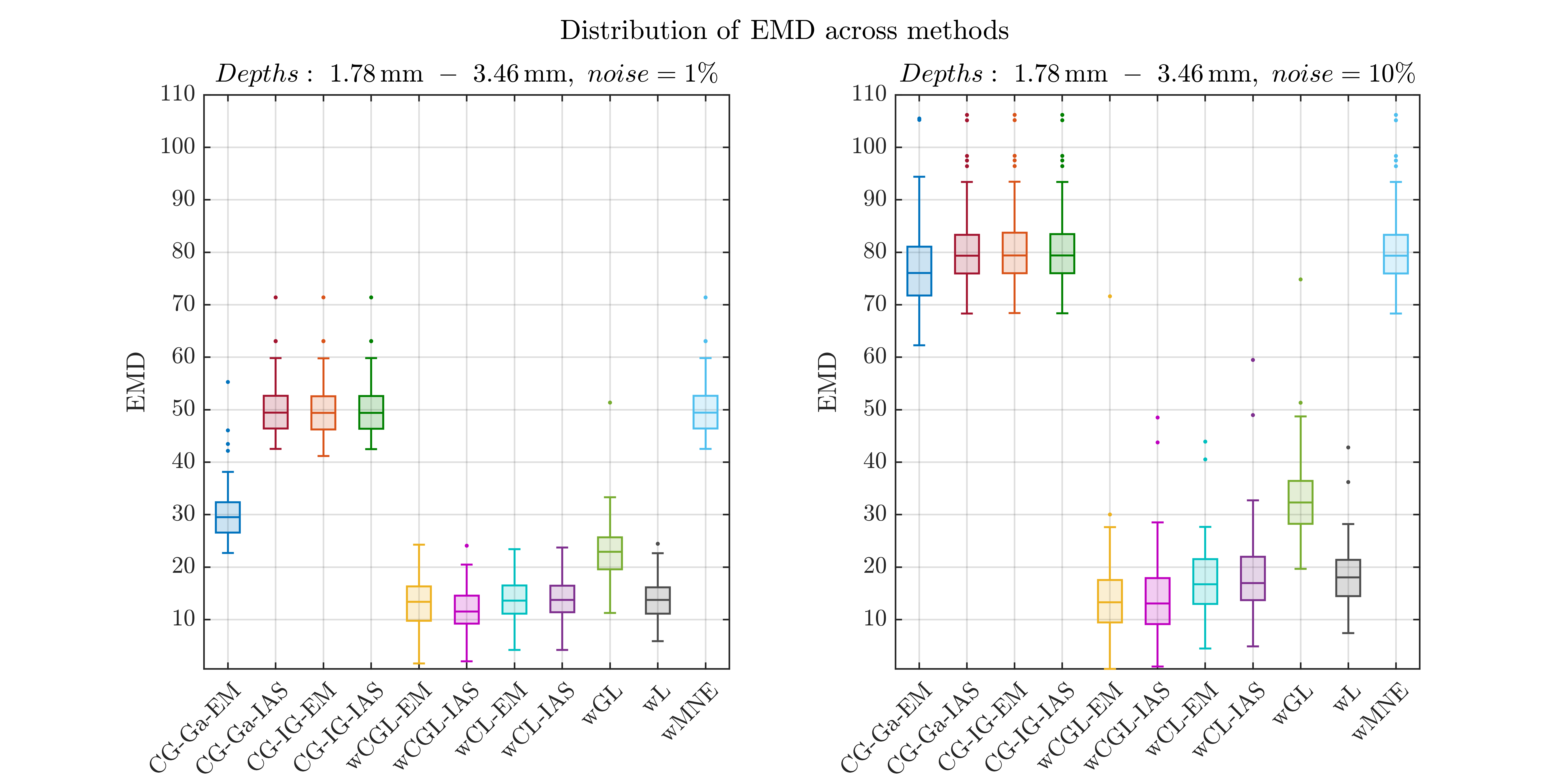}
    \caption{Superfical sources (1.78--3.5 mm). Each panel shows results under 1\% (left) and 10\% (right) noise levels. Shallow sources exhibit lower and more stable EMD values.}
    \label{fig:emd_shallow}
\end{figure}
Figure~\ref{fig:emd_shallow} illustrates the EMD distributions for shallow sources (1.78--3.46 mm), where several methods exhibit lower errors and tighter distributions, though this improvement is not consistent across all approaches.

Overall, for deeper sources (17.8–21.7 mm), a clear degradation in performance is observed across all methods. The EMD values are significantly higher and more dispersed, particularly under higher noise levels (10\%). While methods such as wCGL-EM (IAS) and wCL-EM (IAS) remain relatively stable, others exhibit increased variability, indicating reduced robustness to depth and noise.

 For shallow sources, several focal methods demonstrate comparable performance with generally lower EMD values than for deep simulated source. For most of the methods, the distributions are slightly tighter and less sensitive to noise, indicating improved reliability of source localization in superficial regions. 
Table~\ref{tab:EMDtbl} provides a quantitative summary of the EMD statistics across all methods for both depth ranges and noise levels. The reported median values confirm the trends observed in Figures~\ref{fig:emd_deep} and~\ref{fig:emd_shallow}, with consistently higher errors for deeper sources and under increased noise. For shallow sources (1.78--3.5 mm) at low noise (1\%), several focal methods, such as wCGL and wCL variants, achieve substantially lower median EMD values compared to CG-based approaches, indicating improved localization accuracy. This behavior persists under higher noise (10\%), where these methods remain comparatively robust, while CG-Ga and CG-IG methods exhibit significantly higher median errors.

For deeper sources (17.81--21.74 mm), all methods show increased median EMD values, reflecting the greater difficulty of reconstructing deeper activity. Although wCGL and wCL approaches still outperform others in terms of median error, their variability (as indicated by higher Std and IQR values) increases under higher noise levels but not significantly. In contrast, methods such as wGL and wL exhibit moderate performance. Overall, the table highlights the combined impact of depth and noise on reconstruction accuracy, and confirms that the relative performance of the methods remains strongly method-dependent.

\begin{table}[t]
    \centering
    \small{
    \begin{tabular}{|p{30pt}|p{34pt}|p{32pt}|p{30pt}|p{30pt}|p{30pt}|p{30pt}|p{30pt}|p{30pt}|p{30pt}|p{30pt}|p{30pt}|}
    \hline
    \multicolumn{12}{|c|}{\bf 1 \% of measurement noise}\\ \hline\hline
    \multicolumn{12}{|c|}{ 1.78--3.5 mm depth}\\ \hline
    Method &  CG-Ga EM & CG-Ga IAS & CG-IG EM & CG-IG IAS & wCGL EM & wCGL IAS & wCL EM & wCL IAS & wGL & wL & wMNE\\ \hline    
    Median & 29.52 & 49.43 & 49.42 & 49.40 & 13.40 & 11.56 & 13.61 & 13.75 & 22.92 & 13.73 & 49.43 \\ \hline
Std & 5.24 & 4.87 & 4.87 & 4.82 & 5.15 & 4.16 & 4.40 & 4.32 & 5.41 & 4.16 & 4.87 \\ \hline
IQR & 5.80 & 6.21 & 6.32 & 6.20 & 6.56 & 5.35 & 5.35 & 5.06 & 6.12 & 5.03 & 6.21\\
    \hline
    \multicolumn{12}{|c|}{ 17.81–21.74 mm depth}\\ \hline
    Median & 39.58 & 53.35 & 53.35 & 53.35 & 17.52 & 22.19 & 24.57 & 24.58 & 28.61 & 25.37 & 53.35 \\ \hline
Std & 5.83 & 3.54 & 3.54 & 3.54 & 5.55 & 7.37 & 6.66 & 6.65 & 5.62 & 5.68 & 3.54 \\ \hline
IQR & 7.26 & 4.27 & 4.27 & 4.27 & 7.59 & 10.29 & 8.96 & 9.40 & 7.05 & 6.99 & 4.27\\
    \hline \hline
    \multicolumn{12}{|c|}{\bf 10 \% of measurement noise}\\ \hline\hline
    \multicolumn{12}{|c|}{ 1.78--3.5 mm depth}\\ \hline
    Method &  CG-Ga EM & CG-Ga IAS & CG-IG EM & CG-IG IAS & wCGL EM & wCGL IAS & wCL EM & wCL IAS & wGL & wL & wMNE \\ \hline
Median & 76.05 & 79.36 & 79.41 & 79.42 & 13.32 & 13.07 & 16.75 & 16.97 & 32.33 & 18.04 & 79.36 \\ \hline
Std & 8.94 & 7.82 & 7.89 & 7.93 & 8.21 & 7.34 & 6.43 & 7.58 & 7.51 & 5.47 & 7.82 \\ \hline
IQR & 9.32 & 7.35 & 7.75 & 7.48 & 8.10 & 8.77 & 8.55 & 8.25 & 8.20 & 6.93 & 7.35\\ \hline
    \multicolumn{12}{|c|}{17.81–21.74 mm depth}\\ \hline
   Median & 78.67 & 78.31 & 78.30 & 78.30 & 20.40 & 22.61 & 26.40 & 27.00 & 40.95 & 29.77 & 78.31 \\ \hline
Std & 5.75 & 5.11 & 5.11 & 5.11 & 8.65 & 9.28 & 8.67 & 9.11 & 8.32 & 7.14 & 5.11 \\ \hline
IQR & 5.75 & 4.34 & 4.35 & 4.34 & 10.77 & 12.07 & 12.04 & 12.33 & 8.76 & 8.83 & 4.34 \\ \hline
\\ \hline
    \end{tabular}
    }
    \caption{Summary of Earth Mover’s Distance (EMD) statistics across methods for two representative depth ranges (1.78–3.5 mm and 17.81–21.74 mm) under 1\% and 10\% measurement noise. For each method, the median, standard deviation (Std), and interquartile range (IQR) are reported. Lower EMD values indicate improved localization accuracy, with performance generally degrading for deeper sources and higher noise levels.}
    \label{tab:EMDtbl}
\end{table}

\clearpage

\subsubsection{Assessment of Depth Bias in Source Reconstruction Algorithms}
In this subsection, we examine the relationship between the true source depth and the depth of the reconstructed maximum. By comparing these quantities, we assess whether the algorithms accurately localize sources in correct depths or exhibit systematic biases, such as a tendency to shift the estimated activity toward the cortical surface. It is important to note that this analysis focuses specifically on depth estimation and does not fully reflect overall localization error. In particular, a source may be reconstructed at the correct depth while still being mislocalized in the tangential directions.

\begin{table}[b!]
    \centering
    \small{
    \begin{tabular}{|c|p{23pt}|p{32pt}|p{32pt}|p{28pt}|p{28pt}|p{20pt}|}
    \hline
        \scriptsize{\diagbox{Method}{Depth error}} & $>20$ \unit{\milli\meter} & $(15,20]$ \unit{\milli\meter} & $(10,15]$ \unit{\milli\meter} & $(5,10]$ \unit{\milli\meter} & $(1,5]$ \unit{\milli\meter} & $\leq 1$ \unit{\milli\meter} \\ \hline
        wMNE & 8.3 & 12.0 & 17.5 & 20.0 & 31.2 & 11.0\\ \hline 
        wMCE & 2.1 & 5.0 & 9.7 & 24.8 & 42.2 & 16.2\\ \hline 
        wGL & 0.3 & 2.5 & 6.7 & 24.7 & 45.0 & 20.8\\ \hline 
        CG-IG (IAS) & 8.3 & 12.0 & 17.5 & 20.1 & 30.7 & 11.4\\ \hline 
        CG-IG (EM) & 8.3 & 12.0 & 17.5 & 20.1 & 30.7 & 11.4\\ \hline 
        CG-Ga (IAS) & 8.3 & 12.0 & 17.5 & 20.0 & 31.1 & 11.1\\ \hline 
        CG-Ga (EM) & 2.5 & 4.6 & 13.4 & 24.1 & 38.9 & 16.5\\ \hline 
        wCL (IAS) & 1.8 & 5.0 & 9.6 & 23.1 & 45.0 & 15.5\\ \hline 
        wCL (EM) & 2.3 & 5.6 & 9.6 & 23.1 & 44.1 & 15.3\\ \hline 
        wCGL (IAS) & 2.0 & 3.0 & 11.1 & 21.0 & 40.8 & 22.2 \\ \hline 
        wCGL (EM) & 0.2 & 1.5 & 6.3 & 21.0 & 48.8 & 22.2\\ \hline
    \end{tabular}
    }
    \caption{The table reports the percentage of reconstructed sources whose absolute depth error, defined as $|d_\mathrm{recon}-d_\mathrm{true}|$, falls within predefined error intervals in millimeters for each method. Percentages are computed over 1000 simulated sources and provide a quantitative complement to the reconstructed-versus-true depth plots, highlighting differences in depth localization accuracy across methods. }
    \label{tab:depth_error}
\end{table}

\begin{figure}[t]
     \centering
     \begin{minipage}{0.22\linewidth}
         \centering
         \tiny{CG-Ga (EM)}
         \includegraphics[width=1\linewidth]{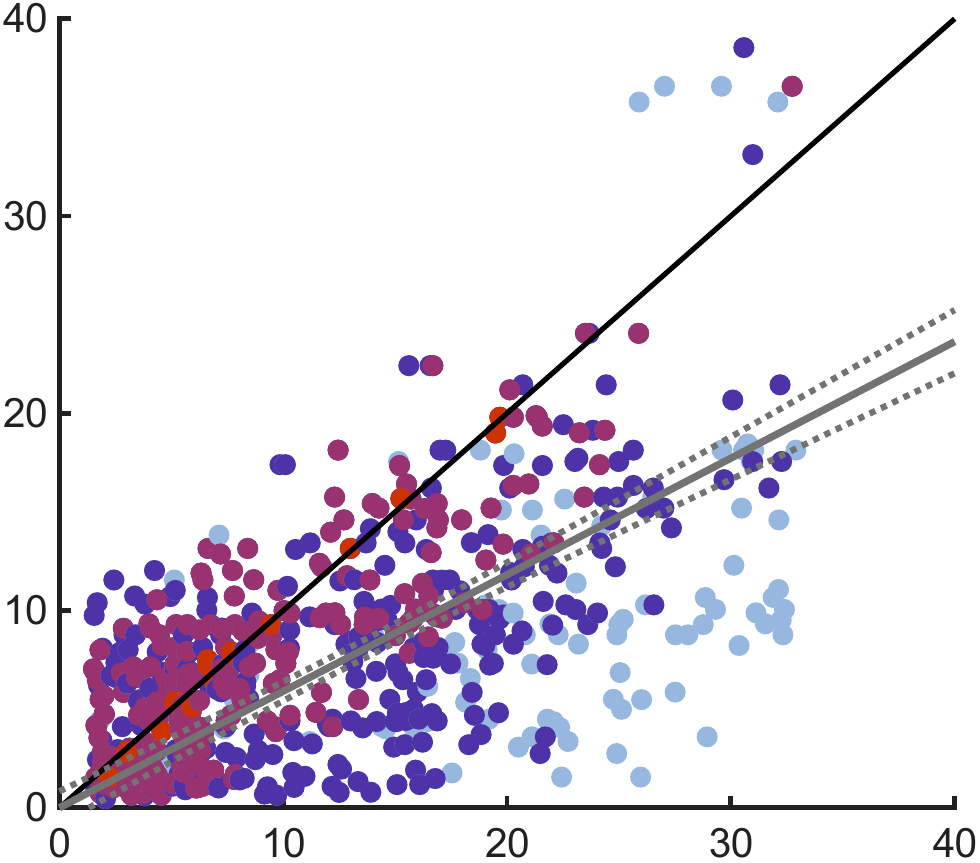}
     \end{minipage}\begin{minipage}{0.22\linewidth}
         \centering
         \tiny{CG-Ga (IAS)}
         \includegraphics[width=1\linewidth]{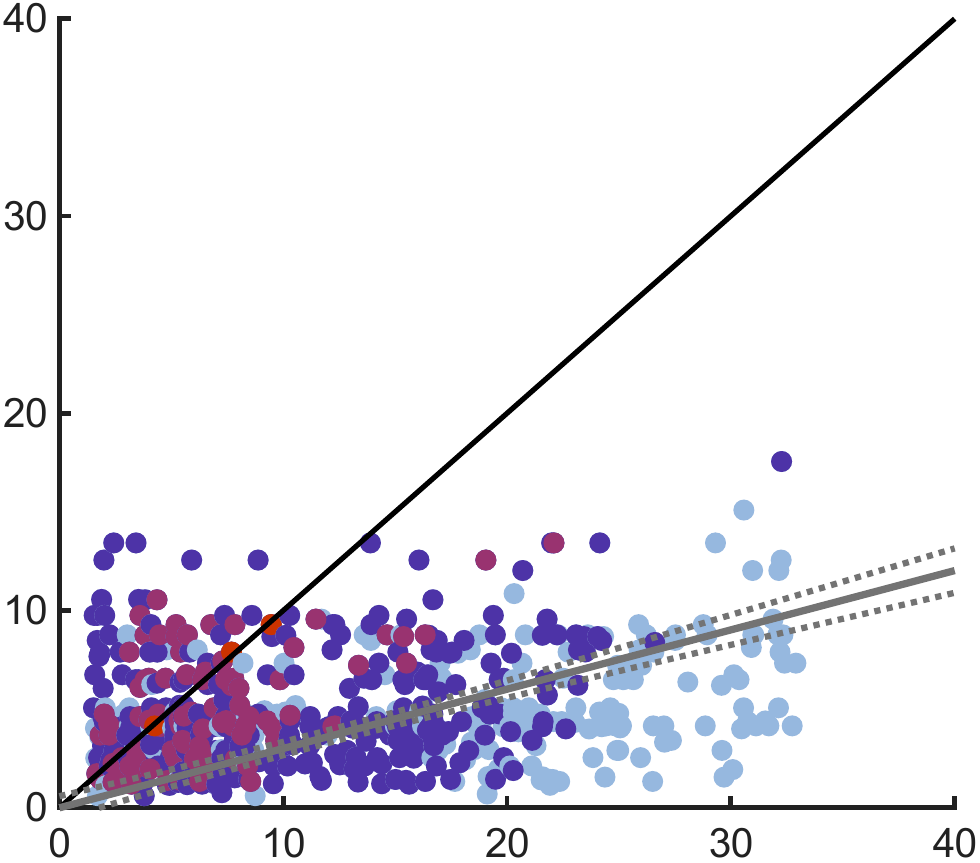}
     \end{minipage}\begin{minipage}{0.22\linewidth}
         \centering
         \tiny{CG-IG (EM)}
         \includegraphics[width=1\linewidth]{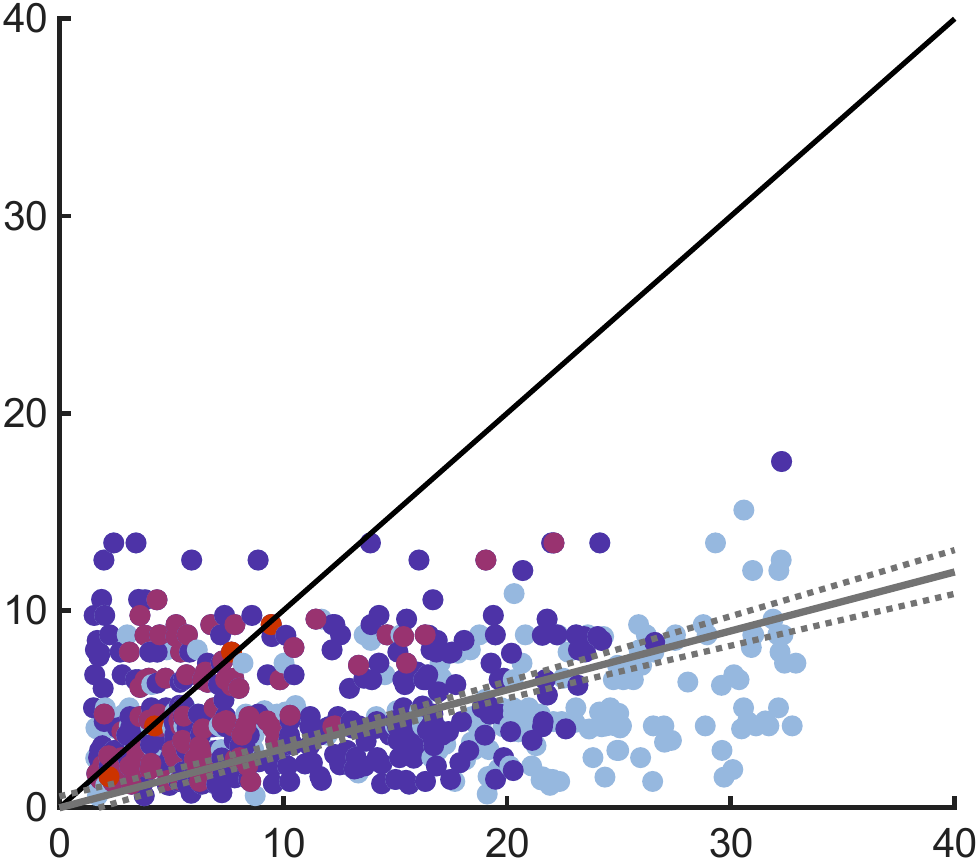}
     \end{minipage}\begin{minipage}{0.22\linewidth}
         \centering
         \tiny{CG-IG (IAS)}
         \includegraphics[width=1\linewidth]{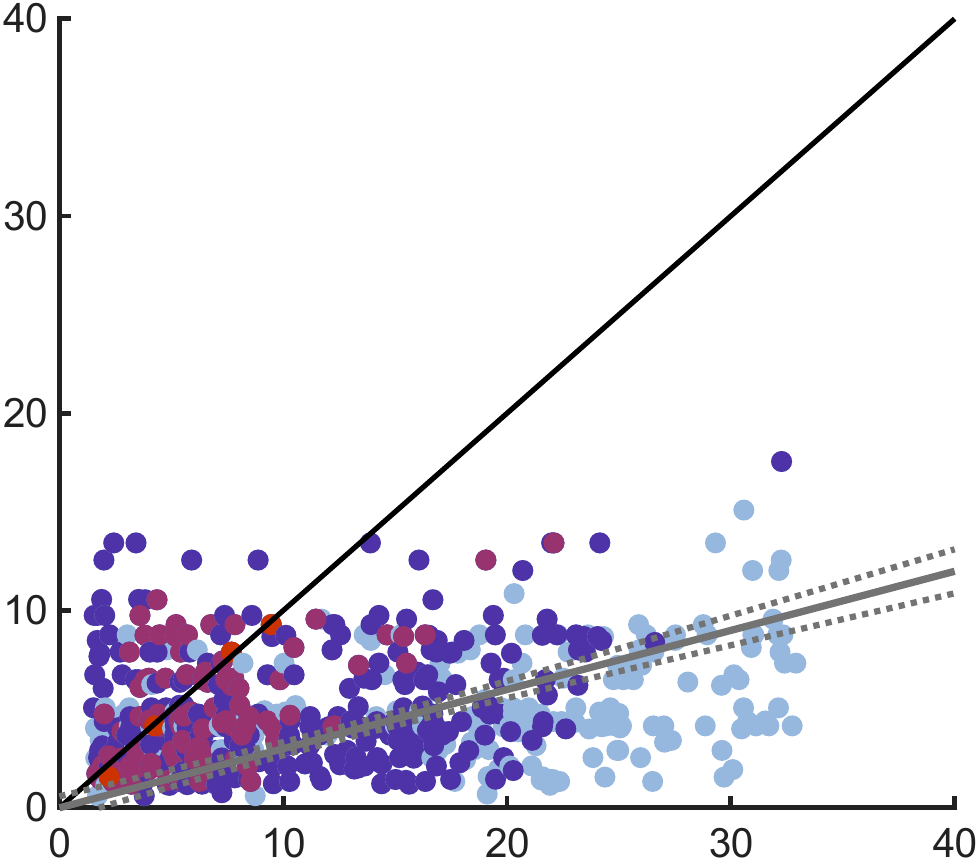}
     \end{minipage}\vspace{0.3cm}
     \begin{minipage}{0.22\linewidth}
         \centering
         \tiny{wCGL (EM)}
         \includegraphics[width=1\linewidth]{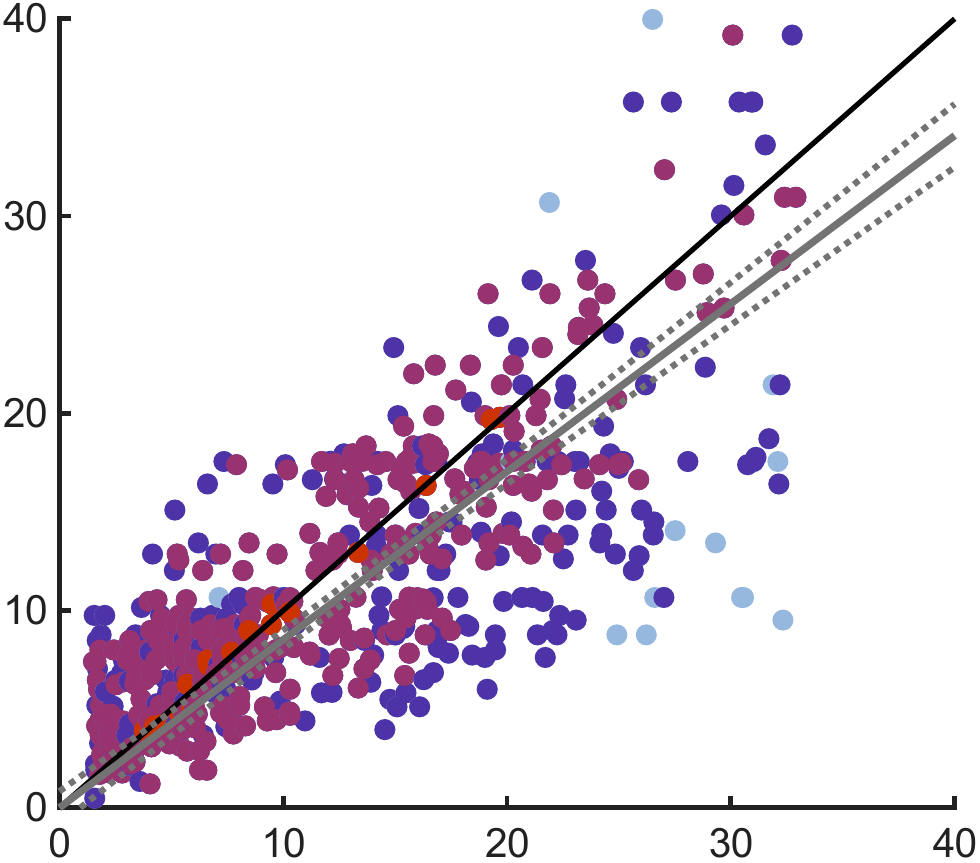}
     \end{minipage}\begin{minipage}{0.22\linewidth}
         \centering
         \tiny{wCGL (IAS)}
         \includegraphics[width=1\linewidth]{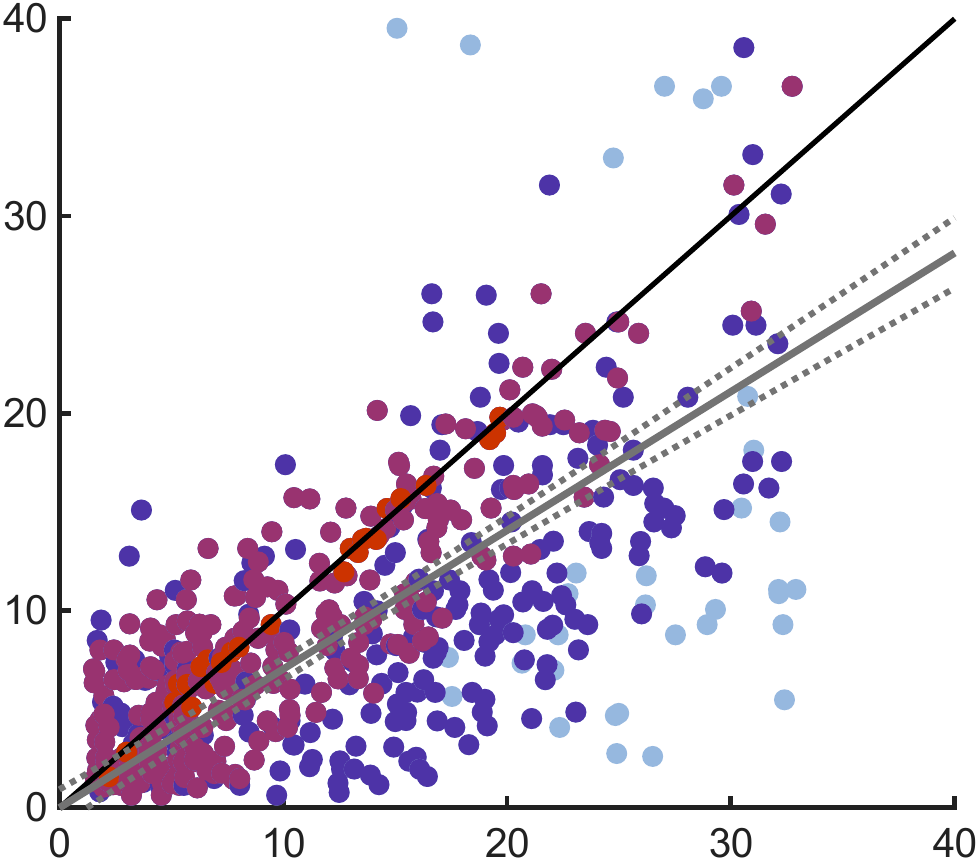}
     \end{minipage}\begin{minipage}{0.22\linewidth}
         \centering
         \tiny{wCL (EM)}
         \includegraphics[width=1\linewidth]{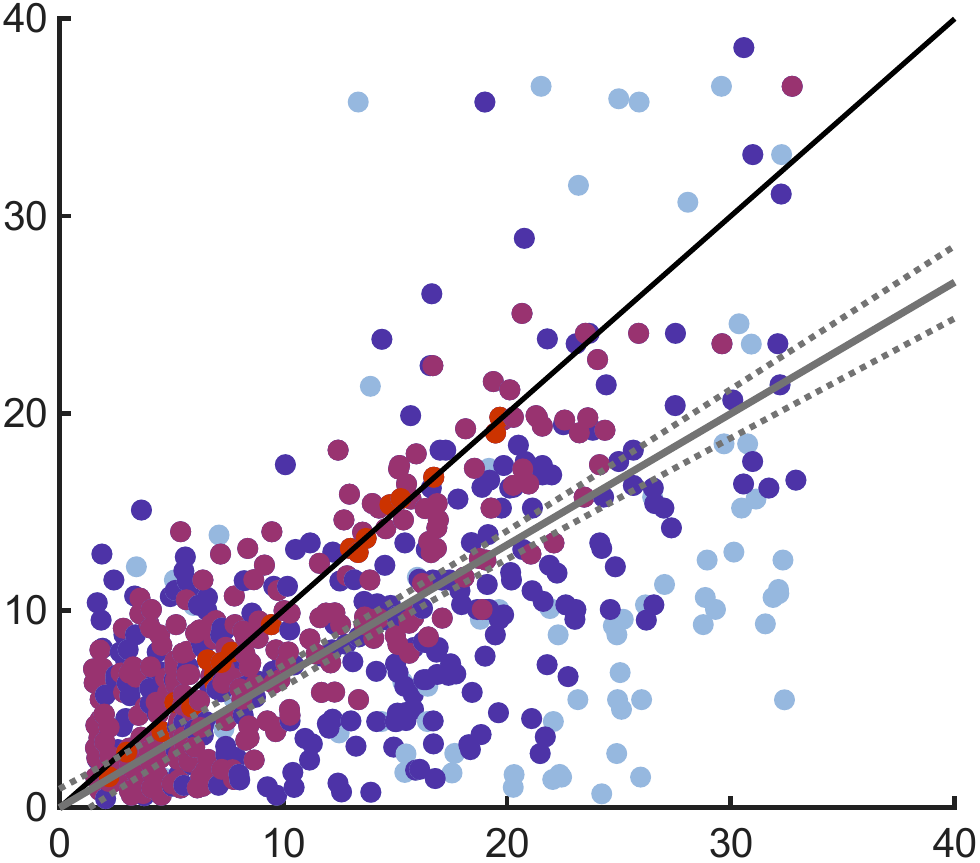}
     \end{minipage}\begin{minipage}{0.22\linewidth}
         \centering
         \tiny{wCL (IAS)}
         \includegraphics[width=1\linewidth]{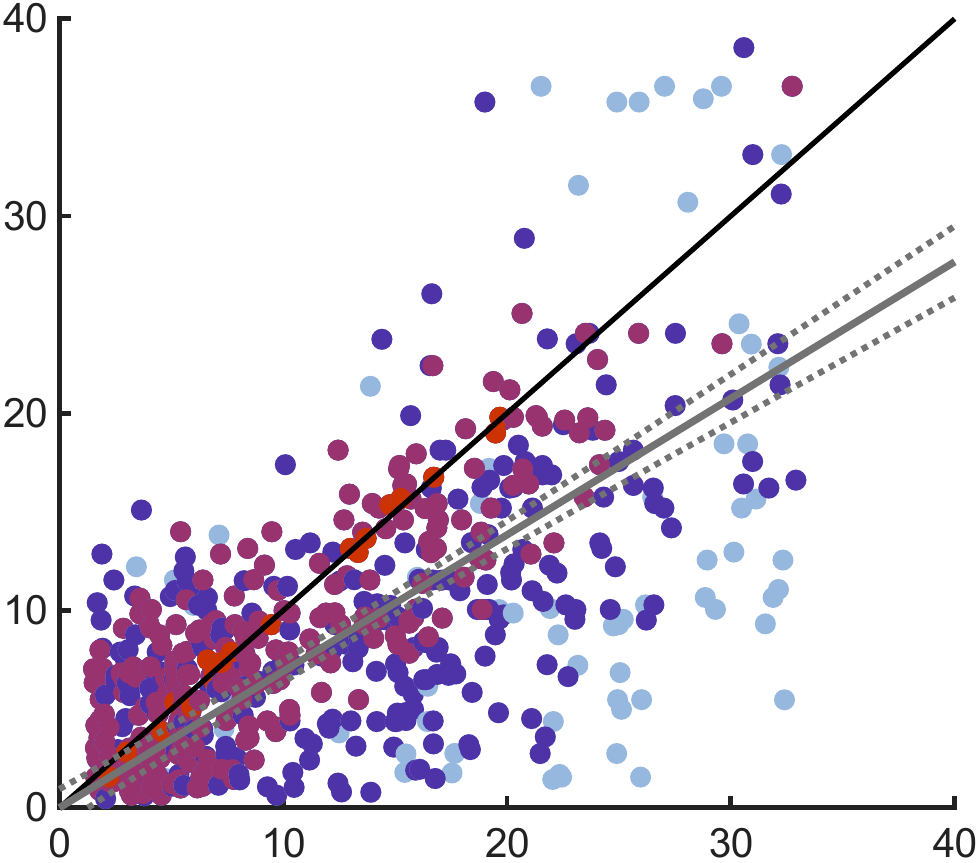}
     \end{minipage}\vspace{0.3cm}
     \begin{minipage}{0.22\linewidth}
         \centering
         \tiny{wMNE}
         \includegraphics[width=1\linewidth]{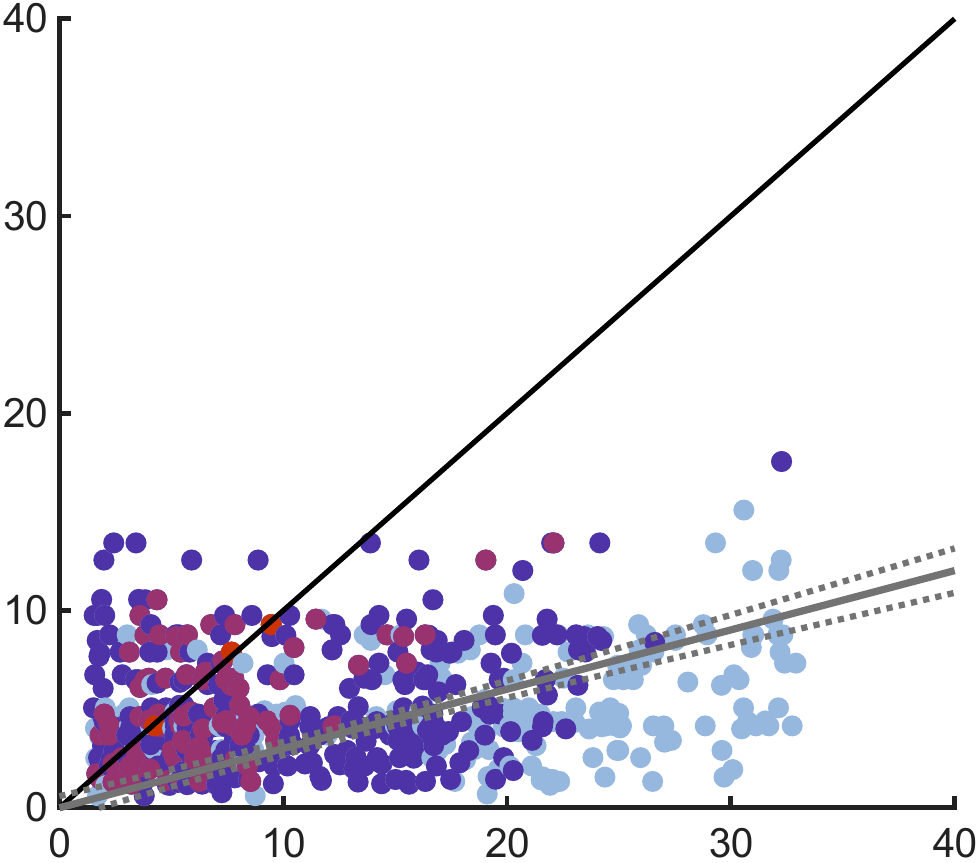}
     \end{minipage}\begin{minipage}{0.22\linewidth}
         \centering
         \tiny{wGL}
         \includegraphics[width=1\linewidth]{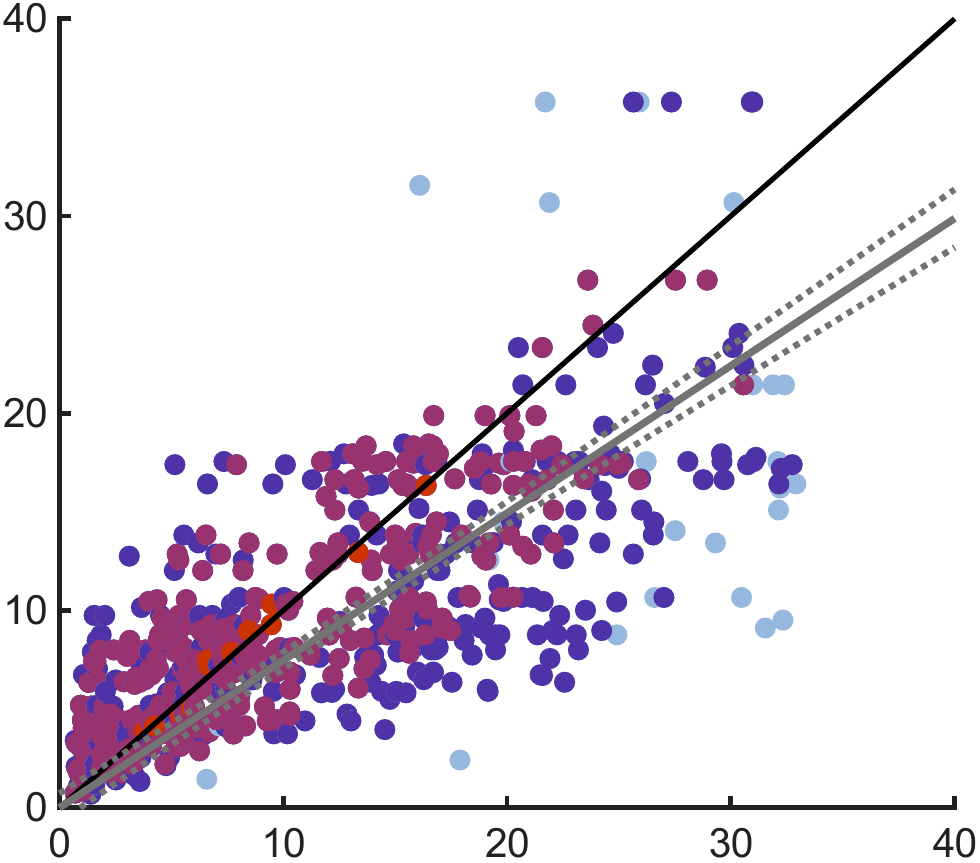}
     \end{minipage}\begin{minipage}{0.22\linewidth}
         \centering
         \tiny{wL}
         \includegraphics[width=1\linewidth]{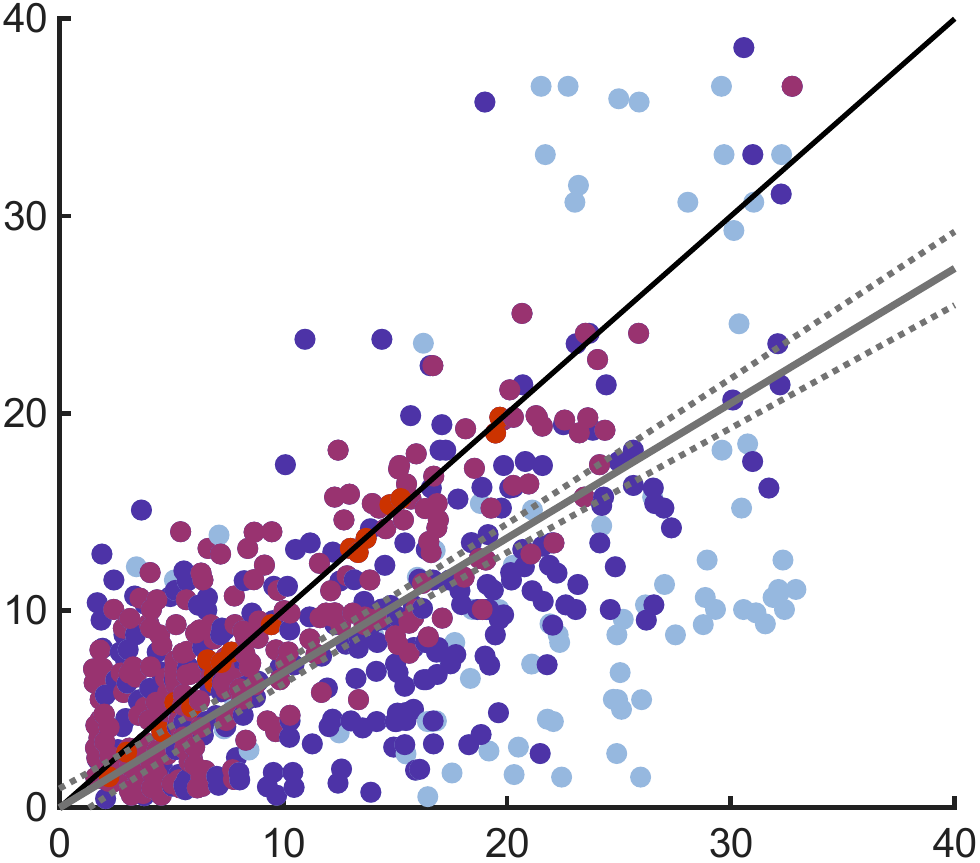}
     \end{minipage}\begin{minipage}{0.1\linewidth}
         \centering
         \includegraphics[width=1\linewidth]{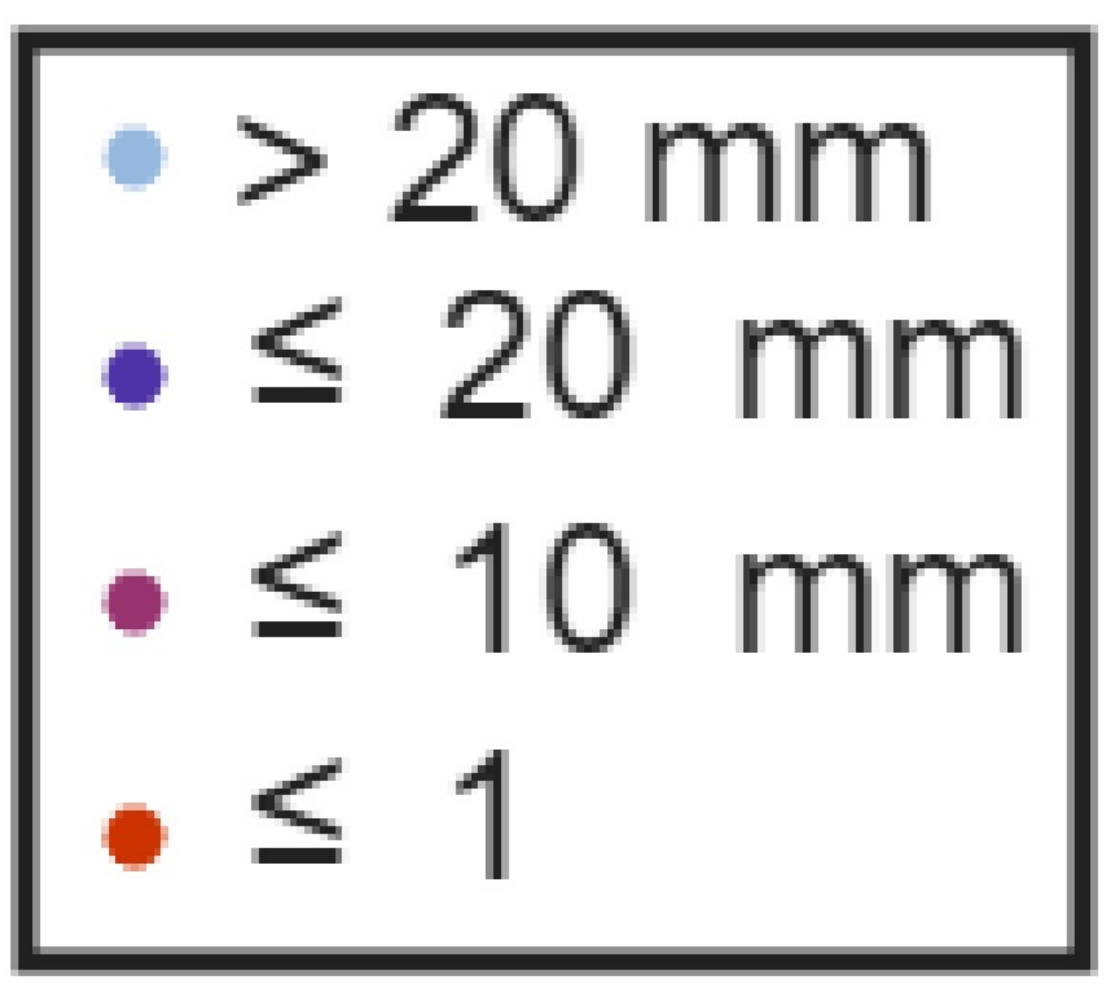}
     \end{minipage}
     \caption{Plots of reconstructed versus true source depths for 1000 simulated sources across different methods. The x-axis shows the depth of the true source maximum, and the y-axis shows the depth of the reconstructed maximum. The black diagonal line indicates perfect agreement (zero depth error). Point color encodes the distance between true and estimated source using the indicated color scale (orange, magenta, purple, blue), with larger magnitude indicating larger depth mismatch. The gray line shows the linear regression line obtained from the data with 95 \% confidence interval presented by the dashed gray lines.}
     \label{fig:depth_reconstruction}
 \end{figure}

Figure~\ref{fig:depth_reconstruction} presents the relationship between the reconstructed maximum depth and the true source depth for all considered methods. The diagonal line represents perfect agreement between the estimated and true depths.

The weighted conditionally Laplace models, particularly wCGL with EM updates, show the best agreement with the diagonal, indicating more accurate depth recovery across a wide range of source depths. Their regression lines closely follow the identity line, suggesting reduced depth bias compared to other methods.

In contrast, the classical conditionally Gaussian models (CG variants) exhibit a clear bias toward superficial regions, as reflected by the flatter slope of the regression lines and the clustering of reconstructed depths at lower values. This indicates that deeper sources are systematically reconstructed closer to the cortical surface. An exception to this is CG-Ga with the EM parameter update algorithm, which exhibits almost as much depth-bias reduction as weighted Laplace (wL).

The wMNE method shows the strongest depth bias, with reconstructed maxima largely concentrated at shallow depths regardless of the true source depth. Similarly, wGL and wL demonstrate intermediate behavior, with improved depth tracking compared to wMNE but still noticeable deviations from the ideal diagonal trend.

Overall, these results confirm that sensitivity weighting combined with EM-based hyperparameter estimation significantly reduces depth bias and improves the ability to recover deeper sources.

Table~\ref{tab:depth_error} provides a quantitative summary of these findings by reporting the percentage of reconstructed sources whose absolute depth error falls within predefined intervals. The results show that the weighted conditionally Laplace models, particularly wCGL (EM), achieve the highest proportion of small depth errors (e.g., within 1–5 mm and 5–10 mm), while also minimizing the occurrence of large errors (greater than 20 mm). In contrast, methods such as wMNE and the classical conditionally Gaussian variants exhibit a higher proportion of large depth errors, consistent with their observed bias toward superficial reconstructions. Overall, the table reinforces the trends observed in Figure~\ref{fig:depth_reconstruction}, confirming the improved depth localization accuracy of the weighted models with EM-based hyperparameter estimation.

Overall, these findings are consistent with the EMD-based analysis presented earlier and further emphasize the importance of algorithm selection in reducing depth-related effects in source reconstruction.

\subsection{Discussion}
In this work, we reviewed principal Bayesian methods for EEG source imaging, focusing on how each approach addresses the ill-posed nature of the inverse problem. Particular emphasis was placed on the role of sensitivity weighting, especially within the Bayesian formalism proposed by \citet{Calvetti2019AutomaticDepthWeighting}. We examined two primary families of Bayesian models for EEG source localization: Gaussian and Conditionally Gaussian models \citep{HamalainenMNE,Calvetti2009}, and Laplace-type models \citep{Uutela1999,Lahtinen2022}. For each class, we provided detailed derivations and algorithmic formulations, highlighting how sensitivity weighting \citep{Calvetti2019SensitivityWeight} can be used to automatically determine model parameters based on the signal-to-noise ratio (SNR) and lead field properties.


Bayesian models employing Laplace-type priors are well known for promoting sparse solutions. Our comparisons using the Earth Mover’s Distance (EMD), however, indicate that these models exhibit increased spatial spread under higher noise levels. Sensitivity weighting partially mitigates this effect, resulting in comparable performance between Laplace and conditional Laplace priors, particularly in terms of noise robustness and depth preservation.

Interestingly, even when SNR is used to modulate sensitivity weights, the anticipated improvement in noise robustness is not uniformly observed across all methods. Nevertheless, sensitivity weighting provides a principled and automatic mechanism for parameter tuning and can substantially reduce the risk of gross mislocalization, especially in cortical source estimation.

We further noted that sensitivity weighting shares conceptual similarities with the weighting employed in weighted Minimum Norm Estimates (wMNE) \citep{Lin2006wMNE}, which also aims to counteract depth bias. Our findings show that beyond approximately 5 mm depth, a pronounced bias emerges across most methods. However, sparsity-promoting approaches—particularly those using group Laplace priors—demonstrate improved localization within the cortical grey matter, suggesting that sensitivity weighting is most effective when combined with focal prior structures.

In addition, we investigated the performance of two hyperparameter updating algorithms, EM and IAS, within conditional models. Both approaches lead to analytically tractable update steps but differ subtly in their treatment of hyperpriors. A key finding is that for the Conditionally Gaussian model with a Gamma hyperprior, EM clearly outperforms IAS in terms of both sparsity and localization accuracy. For other models, EM consistently produces slightly more compact source estimates, albeit at the cost of increased derivational and implementation complexity.

Overall, the improved localization accuracy achievable with hierarchical Bayesian methods offers valuable opportunities for non-invasive diagnostics. For example, such approaches may enhance presurgical planning by improving the identification of functional cortical areas or epileptogenic zones \citep{Diamond2023}.

While sensitivity weighting mitigates the depth effect in EEG, causing the signal to decay further from the potential-measuring electrodes, the difficulty remains in distinguishing focal deep activity from widespread weak superficial activity that extends across the whole neocortex. Both alternative source distributions can yield the same measured observations, so it is up to the prior model to decide which is more likely to be true. Essentially, all the priors investigated here prefer the weaker, more widely spread activity to some degree, as evidenced by the increase in EMD with depth. However, unless the measurement data is almost noiseless, a highly focal estimate of a deep source should not be anticipated due to the subtlety of changes in the measurement data caused by sources at different locations in the deep regions.

Beyond measurement noise and the inherent ill-posedness of EEG source imaging, our results underscore the significant impact of numerical and discretization-induced model discrepancies between the forward simulations and the inverse model. In practice, the inverse solver relies on a discretized approximation of the forward operator, and mismatches arising from mesh resolution, source-space discretization, and assumed conductivity distributions introduce systematic biases that are not captured by noise models or sensitivity weighting alone. These effects manifest as increased EMD values and spatial spreading that persist even when sensitivity weighting is applied. While sensitivity weighting effectively compensates for depth-related sensitivity decay, it cannot correct structural inaccuracies in the forward model. Consequently, numerical bias and conductivity mismatch constitute an additional and distinct limitation, separate from depth bias and stochastic noise, which must be addressed through explicit model discrepancy handling rather than weighting schemes alone.

Thus, although sensitivity weighting can reduce inherent depth bias, improve noise robustness in low to moderate noise regimes, and provide a rational and automatic approach to parameter tuning—features that are particularly valuable in clinical settings where manual tuning is impractical—it cannot mitigate biases arising from numerical inaccuracies and forward-model uncertainty.

Finally, open questions remain regarding how to design priors that are both informative and non-restrictive, particularly in the presence of modeling uncertainty and inter-individual variability. One promising direction is the incorporation of multimodal constraints, such as fMRI or DTI-based connectivity priors \citep{Skudlarski2008DTI}, which can provide complementary structural or functional information to guide source localization without overly constraining the solution space. At the same time, the increased model complexity introduced by such priors highlights the need for scalable inference strategies. In this context, further work should explore efficient approximation techniques, including variational Bayesian and Monte Carlo methods, to make fully Bayesian treatments computationally feasible in clinical and real-time settings. Moreover, extending sensitivity-weighting schemes to time-varying or adaptive priors may offer a natural way to integrate dynamic information and more accurately model temporally evolving neural processes.

\section*{Data and Code Availability}
To be added a link upon acceptance


\section*{Funding}

The work of J.\ Lahtinen was supported by the Research Council of Finland (RCF) through the Flagship of Advanced Mathematics for Sensing, Imaging and modeling (FAME) (359185), and Doctoral Education Pilot on Advanced Mathematics for Modelling, Sensing, and Imaging (DREAM), Ministry of Education and Culture, Finland, VN/3137/2024. A.\ Koulouri was supported by the Institute for Mathematical Innovation, University of Bath, UK.




Supplementary Material (created during production as a web link to online material).

\printbibliography

\end{document}